\documentclass[11pt, letterpaper, reqno]{amsart}
\usepackage{amssymb}
\usepackage{comment}
\usepackage{bbm}
\usepackage{tikz}
\usepackage[all,cmtip]{xy}
\usepackage{hyperref}
\usepackage{stmaryrd}
\usepackage{mathtools}

\newcounter{mt}

\newtheorem{MainTheorem}[mt]{Theorem}

\newtheorem{Proposition}{Proposition}[section]
\newtheorem{Definition}[Proposition]{Definition}
\newtheorem{Lemma}[Proposition]{Lemma}
\newtheorem{Theorem}[Proposition]{Theorem}
\newtheorem{Corollary}[Proposition]{Corollary}
\newtheorem{Remark}[Proposition]{Remark}

\newtheorem{Example}[Proposition]{Example}

\DeclareMathOperator{\Val}{Val}

\DeclareMathOperator{\nc}{nc}

\DeclareMathOperator{\AGr}{AGr}
\DeclareMathOperator{\Gr}{Gr}

\DeclareMathOperator{\Res}{I}

\DeclareMathOperator{\res}{res}

\DeclareMathOperator{\proj}{pr}

\DeclareMathOperator{\push}{push}
\DeclareMathOperator{\id}{Id}

\DeclareMathOperator{\vol}{vol}

\DeclareMathOperator{\Dens}{Dens}
\DeclareMathOperator{\Span}{Span}

\DeclareMathOperator{\Image}{Im}
\DeclareMathOperator{\Supp}{Supp}

\DeclareMathOperator{\Stab}{Stab}
\DeclareMathOperator{\Sym}{Sym}

\DeclareMathOperator{\Kl}{Kl}

\DeclareMathOperator{\Sc}{Sc}

\DeclareMathOperator{\ori}{or}

\DeclareMathOperator{\GL}{GL}
\DeclareMathOperator{\SL}{SL}

\DeclareMathOperator{\OO}{O}
\DeclareMathOperator{\SO}{SO}
\DeclareMathOperator{\sign}{sign}
\DeclareMathOperator{\KS}{KS}

\newcommand{\R}{\mathbb{R}}
\newcommand{\N}{\mathbb{N}}
\newcommand{\Z}{\mathbb{Z}}

\makeatletter
\def\moverlay{\mathpalette\mov@rlay}
\def\mov@rlay#1#2{\leavevmode\vtop{%
		\baselineskip\z@skip \lineskiplimit-\maxdimen
		\ialign{\hfil$\m@th#1##$\hfil\cr#2\crcr}}}
\newcommand{\charfusion}[3][\mathord]{
	#1{\ifx#1\mathop\vphantom{#2}\fi
		\mathpalette\mov@rlay{#2\cr#3}
	}
	\ifx#1\mathop\expandafter\displaylimits\fi}
\makeatother

\newtheorem{note}{Remark} 
\def\note#1{\ifvmode\leavevmode\fi\vadjust{\vbox to0pt{\vss
			\hbox to 0pt{\hskip\hsize\hskip1em
				\vbox{\hsize3.5cm\small\raggedright\pretolerance10000
					\noindent #1\hfill}\hss}\vbox to8pt{\vfil}\vss}}}

\renewcommand{\emph}{\textsl}

\usepackage[abbrev]{amsrefs}
\BibSpec{article}{%
  +{}{\PrintAuthors}  		{author}
  +{,}{ \textit}     		{title}
  +{,}{ }             		{journal}
  +{}{ \textbf}       		{volume}
  +{}{ \parenthesize} 		{date}
  +{}{, no. } 		{number}
  +{,}{ }      	      		{conference}
  +{,}{ }      	      		{book}
  +{,}{ }            		{pages}
  +{,}{ }            	 	{note}
  +{,}{ }            	 	{status}
  +{,}{  \texttt } {eprint}
  +{.}{}              {transition}
}

\BibSpec{book}{%
  +{}{\PrintAuthors}  {author}
  +{,}{ \textit}      {title}
  +{,}{ }             {publisher}
  +{,}{ }             {place}
  +{,}{ }             {date}
  +{.}{}              {transition}
}

\BibSpec{incollection}{%
    +{}  {\PrintAuthors}                {author}
    +{,} { \textit}                     {title}
    +{.} { }                            {part}
    +{:} { \textit}                     {subtitle}
    +{,} { \PrintContributions}         {contribution}
    +{,} { \PrintConference}            {conference}
    +{,} { }                            {booktitle}
    +{,} { pp.~}                        {pages}
    +{,}{ }       		{series}
    +{}{ \textbf}       		{volume}
    +{,} { }                            {publisher}
    +{,} { }                 {date}
    +{,} { }                            {status}
    +{,} { \PrintDOI}                   {doi}
    +{,} { available at \eprint}        {eprint}
    +{}  { \parenthesize}               {language}
    +{}  { \PrintTranslation}           {translation}
    +{;} { \PrintReprint}               {reprint}
    +{.} { }                            {note}
    +{.} {}                             {transition}
}

\date{\today}

\begin{document}
	
\title{Convex valuations from Whitney to Nash}
	\author{Dmitry Faifman}
	\email{dmitry.faifman@umontreal.ca}
	\address{School of Mathematical Sciences, Tel Aviv University, Tel Aviv 6997801, Israel}
	\curraddr{D\'epartement de Math\'ematiques et de Statistique, Universit\'e de
		Montr\'eal, CP 6128 succ Centre-Ville, Montr\'eal, QC H3C 3J7, Canada}
	
	\author{Georg C.\ Hofst\"atter}
	\email{georg.hofstaetter@tuwien.ac.at}
	\address{Institute f. Mathematics, Friedrich-Schiller-University Jena, 07743 Jena, Germany}
	\curraddr{Institute of Discrete Mathematics and Geometry, TU Wien, 1040 Vienna, Austria}
	
	\begin{abstract}
		
		We consider the Whitney problem for valuations: does a smooth $j$-homogeneous translation-invariant valuation on $\R^n$ exist that has given restrictions to a fixed family $S$ of linear subspaces? A  necessary condition is compatibility: the given valuations must coincide on intersections. We show that for $S=\Gr_r(\R^n)$, the grassmannian of $r$-planes, this condition becomes sufficient once $r\geq j+2$. This complements the Klain and Schneider uniqueness theorems with an existence statement. Informally, the obstruction for a $j$-density to extend to a $j$-homogeneous valuation is localized in a single dimension, namely $j+2$. 
	
		We then look for conditions on $S$ when compatibility is also sufficient for extensibility, in two distinct regimes: finite arrangements of subspaces, and compact submanifolds of the grassmannian. In both regimes we find unexpected flexibility. As a consequence, we prove a Nash-type theorem for valuations on compact manifolds, from which in turn we deduce the existence of Crofton formulas for all smooth valuations on manifolds, answering a question of Fu. As an intermediate step of independent interest, we construct Crofton formulas for all odd translation-invariant valuations.
	\end{abstract}
\thanks{ 
	\\\indent DF was supported by the Israel Science Foundation grant No. 1750/20. GH was partially supported by the European Research Council 
	under the European Union’s Horizon 2020 research and innovation
	programme (grant agreement No 770127).
}
	
\maketitle

	\section{Introduction and results}\label{sec:intro}
	Valuations, that is finitely additive functionals on convex compact sets, belong to the most basic notions in (convex) geometry. While their emergence is often traced to Dehn's solution of Hilbert's third problem, many classical geometric quantities are instances of valuations. Some notable examples, such as intrinsic volumes (quermassintegrals) or affine surface area, can in fact be characterized as the unique valuations with certain symmetry and an analytic restriction such as continuity. 
	
	Consequently, valuation theory itself has gained a lot of interest and became a very active research subject, leading to numerous fruitful structural insights that could directly be applied to concrete problems in geometry, see, e.g., \cite{Dehn1901, Hadwiger1957, McMullen1977, Klain2000, Schneider1996, Ludwig2010, Alesker1999} or \cite{Klain1997, Schneider2014} for an overview.

	Notably, valuation theory was closely tied to integral geometry by Hadwiger \cite{Hadwiger1957}, who identified its role in the computation of Blaschke style kinematic formulas. Valuation theory has since played a pivotal role in the determination of kinematic formulas for various groups, see e.g. \cite{alesker_multiplicative,Fu2006,Bernig2006,Bernig2011, Bernig2014,Bernig2017b}. Further connections to integral geometry, this time Gelfand style, were uncovered in \cite{Alesker2009}.
	
	The proof of McMullen's conjecture by Alesker \cite{Alesker2001}, asserting the density of mixed volumes among continuous translation-invariant valuations, proved a turning point in valuation theory, with a rich algebraic structure subsequently uncovered therein. For an overview of this emergent theory, see, e.g., \cite{Alesker2014c}.
	The notion of smooth valuations on manifolds was subsequently introduced by Alesker \cite{Alesker2006b, Alesker2006} and studied by various authors  \cite{Bernig2007, Alesker2007,Alesker2008, Alesker2009, Alesker2017,Alesker2012}, in particular as an approach to integral geometry on more general spaces. For a survey of valuations on manifolds, see \cite{Alesker2007b}. 
	Development in valuation theory has since grown and accelerated, through the works of Alesker, Bernig, Fu and many others. For a further sample of works, see \cite{Alesker2003, Alesker2011b,  Alesker1999,  Alesker2014, Bernig2022, Bernig2022b, Faifman2021, Fu2019,Kotrbaty2020, Solanes2017,Wannerer2018}.  
		 
	A smooth valuation on an $n$-dimensional manifold $M$, which for simplicity we assume oriented, is a functional $\phi: \mathcal{P}(M) \to \R$ on the set $\mathcal{P}(M)$ of compact smooth submanifolds with corners of $M$, which is of the form 
	\begin{align*}
		\phi(A) = \int_{\nc(A)} \omega + \int_A \theta, \qquad A \in \mathcal{P}(M),
	\end{align*}
	where $\nc(A)$ is the conormal cycle of $A \in \mathcal{P}(M)$ (see, e.g., \cite{Zaehle1986}), which is a Lipschitz submanifold of the cosphere bundle $\mathbb P_M=\mathbb P_+(T^*M)$ of $M$, and $\omega \in \Omega^{n-1}(\mathbb P_M)$ and $\theta \in \Omega^n(M)$ are differential forms. The space of smooth valuations on $M$ is denoted by $\mathcal{V}^\infty(M)$. When $M=V$ is a linear space, then the subspace of translation-invariant valuations, denoted by $\Val^\infty(V)$, coincides with the space of smooth vectors of the $\GL(V)$-representation on $\Val(V)$, the space of continuous, translation-invariant valuations on $V$. Moreover, by a theorem of McMullen \cite{McMullen1977}, the space $\Val^\infty(V)$ is graded by the degree of homogeneity, defining the subspaces $\Val_k^\infty(V)$ of $k$-homogeneous valuations, $k=0, \dots, \dim V$. This grading can be further refined by parity. We denote by $\Val_k^{\pm, \infty}(V)$ the even/odd valuations, and similarly for $\Val(V)$. See Section~\ref{sec:bgval} for more background on  valuations.
	 
	The most important examples of valuations on a Riemannian manifold $M$ are the \emph{intrinsic volumes}, which extend the notion of intrinsic volumes on Euclidean space. Defined by integrals of certain invariant polynomials of the curvature tensor (the so-called Lipschitz-Killing curvatures), they appear as coefficients in Weyl's famous tube formula. Alternatively, they can be obtained as restrictions of the intrinsic volumes of $\R^N$ to $M$ by any isometric embedding, where it was Weyl's observation~\cite{Weyl1939}, now generally referred to as Weyl's principle, that the restrictions are independent of the embedding. This conceptually simple description allows transferring linear tools, such as Crofton formulas, from Euclidean space. Recently, generalizations of the Weyl principle were considered for contact~\cite{Faifman_contact}, pseudo-Riemannian~\cite{Bernig2022b,Bernig2021b}, Finsler~\cite{Faifman2021} and K\"ahler manifolds~\cite{Bernig2022}.	

	The smooth translation-invariant valuations in linear space enjoy a rich structure, with a large toolset available for their handling, notably Alesker's irreducibility theorem \cite{Alesker2001}, the Fourier transform \cite{Alesker2011b}, and explicit dense subspaces provided by mixed volumes or Crofton formulas. With this motivation in mind, we take a more general viewpoint and seek to represent a general $\phi\in\mathcal V^\infty(M)$ as a restriction $\phi = e^\ast \Phi$ of $\Phi\in \Val^\infty(\R^N)$ under some smooth embedding $e: M \hookrightarrow \R^N$.
	
	Clearly, the restriction of a translation-invariant valuation must be constant when evaluated at points $\{x\}$, $x \in M$. Recalling the distinguished subspace $\mathcal W_1^\infty(M)\subset \mathcal{V}^\infty(M)$ (see \cite{Alesker2006}) given by 
	\begin{align*}
		\mathcal W_1^\infty(M)=\{\phi\in \mathcal V^\infty(M): \phi(\{x\})=0,\,\,\forall x\in M\},
	\end{align*}
	we deduce that only valuations in $\mathcal W_1^\infty(M)\oplus \Span\{\chi\}$, where $\chi$ is the Euler characteristic, can possibly be the restriction of a translation-invariant valuation. Our first main result shows that no other obstructions exist, and moreover, a single embedding can be used to obtain all such smooth valuations. In analogy with the setting of Riemannian geometry, we refer to it as a Nash-type embedding theorem for valuations. Let us stress however that the Riemannian Nash embedding theorem is a substantially deeper result.
	\begin{MainTheorem}\label{mthm:nashThm}
		Suppose that $M$ is a compact smooth manifold. Then there exists an embedding $e:M\hookrightarrow \R^N$ such that $e^*\left(\Val^\infty(\R^N)\right)=\mathcal W_1^\infty(M)\oplus \Span\{\chi\}$.
	\end{MainTheorem}
	Theorem~\ref{mthm:nashThm} is a direct consequence of Theorem~\ref{mthm:nashThm2}, for which we require a definition. Let $\mathbb{P}_r$ be the bundle over $\Gr_r(\R^n)$ with fiber $\mathbb{P}_+(E)$ over $E \in \Gr_r(\R^n)$. Define the map $\theta_r: \mathbb{P}_r \to \mathbb{P}_+(\R^n)$ by $(E, [v]) \mapsto [v]$.
	
	\begin{Definition}\label{def:musical_chairs}
		A closed submanifold $Z \subset \Gr_r(\R^n)$ is called \emph{perfectly self-avoiding}, if $\theta_r$ restricted to $\mathbb{P}_r|_Z$ is an embedding.
	\end{Definition}
	
	In particular, if $Z$ is perfectly self-avoiding then $\theta_r(\mathbb{P}_r|_Z)$ is an embedded submanifold of $\mathbb{P}_+(\R^n)$, and $E\cap E'=\{0\}$ for all $E,E'\in Z$ .

	\begin{MainTheorem}\label{mthm:nashThm2}
		Suppose that $Z \subset \Gr_r(\R^n)$ is a perfectly self-avoiding compact submanifold, $1 \leq r \leq n-1$, and let $1 \leq j \leq r$. Then for every smooth assignment $E\mapsto \phi_E\in\Val^{\infty}_j(E)$, $E\in Z$, there exists a valuation $\psi \in \Val_j^\infty(\R^n)$ such that $\phi_E = \psi|_E$ for all $E \in Z$.
	\end{MainTheorem}

	For the proof of Theorem~\ref{mthm:nashThm}, we choose a generic embedding $e: M \hookrightarrow \R^N$ with $N$ large enough, which then induces an embedding of the tangent bundle of $M$ as a perfectly self-avoiding submanifold of $\Gr_n(\R^n)$. An alternative proof of Theorem~\ref{mthm:nashThm}, using the Alesker product on smooth valuations, the irreducibility theorem, and the nuclearity of Fr\'echet spaces of smooth valuations, is given in Appendix~\ref{app:NuclearNash}.
	
	As an application of Theorem~\ref{mthm:nashThm}, we deduce that all smooth valuations on $M$ are given by rather explicit Crofton formulas, and thereby answer a question of Fu \cite{Fu2016}. We do so by restricting the appropriate Crofton formulas on $\R^N$, which we obtain for odd translation-invariant valuation as an intermediate step of independent interest. In order to state the result, let $\AGr_{k}(\R^N)$ denote the grassmannian of $k$-dimensional affine subspaces, and $\mathrm{HGr}_{k}(\R^N)$ the space of $k$-dimensional affine half-spaces in $\R^N$.
	
	\begin{MainTheorem}\label{mthm:croftonFormulaMf}
		Suppose that $M$ is a compact smooth manifold of dimension $n$. Then there exist $N\in\N$, $R>0$ and an embedding $e: M \hookrightarrow \R^N$ such that, if $\phi \in  \mathcal{V}^{\infty}(M)$, then there exist $C \in \R$ and compactly supported smooth measures $\mu$ on $\R^N \times [0,R]$, $m_j$ on $\AGr_{N-j}(\R^N)$ and $\mu_j$ on $\mathrm{HGr}_{N-j+1}(\R^N)$, such that for $A \in \mathcal{P}(M)$,

		\begin{align}\label{eq:crofton}
			\phi(A) = \, &  C \chi(A) + \int_{\R^N \times [0,R]}\!\!\!\!\!\!\!\!\!\!\!\!\!\!\! \chi(B_\rho(y) \cap e(A))\, d\mu(y,\rho) \\ \nonumber
			&+ \sum_{j=1}^{n} \int_{\AGr_{N-j}(\R^N)}\!\!\!\!\!\!\!\!\!\!\!\!\!\!\!\!\!\!\!\!\!\! \chi(E\cap e(A))\,dm_j(E) + \sum_{j=1}^{n-1} \int_{\mathrm{HGr}_{N-j+1}(\R^N)}\!\!\!\!\!\!\!\!\!\!\!\!\!\!\!\!\!\!\!\!\!\!\!\!\!\!\!\chi(H\cap e(A)) d\mu_j(H),
		\end{align}
		where $B_\rho(y)$ denotes the Euclidean ball of radius $\rho>0$ centered at $y\in\R^N$.
	\end{MainTheorem}

	The technical nature of the Crofton formulas above is somewhat delicate. In particular, the intersections $B_\rho(y) \cap e(A)$, $E \cap e(A)$ and $H \cap e(A)$ above are transversal only for almost all values $(y, \rho), E$ and $H$, leading to questions of integrability. In Section~\ref{sec:crofton} we will further show that the formula above holds also as a Gelfand--Pettis (weak) integral in the space of generalized valuations, from which a simple formula for the Alesker product of valuations on manifolds can be derived.
	
	\bigskip
	
	 Let us now state our next main result. 
	\begin{MainTheorem}\label{mthm:resProperty_simple}
		Assume a smooth assignment $\Gr_r(\R^n)\ni E\mapsto \phi_E\in\Val^\infty_j(E)$ is given such that $\phi_E|_{E\cap E'}=\phi_{E'}|_{E\cap E'}$ for all pairs of subspaces $E, E'$. If $r\geq j+2$, then one can find $\phi\in\Val_j^\infty(\R^n)$ such that $\phi|_E=\phi_E$ for all $E\in\Gr_r(\R^n)$.
	\end{MainTheorem}
	Evidently the condition of compatibility on intersections is necessary for the existence of $\phi$.
	To facilitate the discussion, let us introduce some notation.
	
	Denote by $\Val_j^\infty(\Gr_r(\R^n))$ the Fr\'echet bundle with fiber $\Val_j^\infty(E)$ over $E\in\Gr_r(\R^n)$, $0 < r < n$, and let $V^\infty_j(r, \R^n)$ be the Fr\'echet space of smooth global sections. The restriction map $\res_r$ is defined for $j \geq 1$ by
	\begin{align*}
		\res_r:\Val_j^\infty(\R^n)\to V^\infty_j(r, \R^n), \qquad \res_r(\phi)(E)=\phi|_E.
	\end{align*}
	Observe that for $j=r$, $V^\infty_r(r, \R^n)$ is, by Hadwiger's theorem~\cite{Hadwiger1957}, the space of smooth assignments $\Gr_r(\R^n)\ni E\mapsto \phi_E^r \in \Dens(E)$, which can be identified with $C^\infty(\Gr_r(\R^n))$. 
	Thus $\res_r$, restricted to even valuations, is just the Klain map \begin{equation}\label{eq:klain}\Kl: \Val_r^{+,\infty}(\R^n)\to C^\infty(\Gr_r(\R^n)),\end{equation} which is injective by Klain's theorem \cite{Klain2000}. A theorem of Schneider \cite{Schneider1996} similarly asserts  that $\res_{r+1}:\Val_r^{-,\infty}(\R^n)\to V^\infty_{r}(r+1,\R^n)$ is injective. Thus we have
	\begin{Theorem}[\cite{Klain2000,Schneider1996}]\label{thm:KlSchnInjective}
		$\res_r:\Val_j^\infty(\R^n)\to V_j^\infty (r, \R^n)$ is injective if $r\geq j+1$. The restriction to even valuations is injective also for $r=j$.
	\end{Theorem}

	Theorem \ref{mthm:resProperty_simple} characterizes the image of $\res_r$, complementing Theorem~\ref{thm:KlSchnInjective}. Any element of $\mathrm{Image}(\res_r)$ lies in the subspace of compatible sections, defined as 
	\begin{align*}
		U_j^{\infty}(r, \R^n)=\{\psi\in V^{\infty}_j(r, \R^n): \psi_E|_{E\cap E'}=\psi_{E'}|_{E\cap E'}, \forall E, E' \in \Gr_r(\R^n) \}.
	\end{align*}
	Theorem \ref{mthm:resProperty_simple} can now be stated as follows.
	\begin{Theorem}\label{mthm:resProperty}
		If $r\geq j+2$, then $\mathrm{Image}(\res_r:\Val_j^\infty(\R^n)\to V_j^\infty(r, \R^n))=U_j^\infty(r, \R^n)$. 
	\end{Theorem}
	Note that for $r<j$ the image of $\res_r$ is zero by McMullen's decomposition theorem~\cite{McMullen1977}, while for $r=j,j+1$ Theorem~\ref{mthm:resProperty} is in general  false.
	Theorem~\ref{mthm:resProperty} implies that the only obstruction for a $j$-density to represent a valuation lies in dimension $j+2$. This appears to be a newly observed phenomenon.

	Let us briefly comment on the extension problem in the even case. The image of the Klain map is known \cite{Alesker2004c} to coincide with the image of the cosine transform $\mathcal C_r$ on $\Gr_r(\R^n)$. The latter is an isomorphism by Alexandrov's theorem for $r=1, n-1$, while for $2\leq r\leq n-2$ its image is a proper subspace of $C^\infty(\Gr_r(\R^n))$. Furthermore, Alesker and Bernstein described the irreducible $\SO(n)$-modules appearing in $\Image(\mathcal C_r)$. As $C^\infty(\Gr_r(\R^n))$ is a multiplicity-one $\SO(n)$-module, this description identifies $\Image(\mathcal C_r)$ uniquely inside $C^\infty(\Gr_r(\R^n))$. Recently it was shown in \cite{Faifman2022} that any nonzero $f\in \mathrm{Image}(\mathcal C_r)$ cannot be supported inside an open Schubert cell of $\Gr_r(\R^n)$, yielding a geometric obstruction. While an analytic description of the image of $\mathcal C_r$ is still lacking, it seems possible, in light of \cite{alesker_gourevitch_sahi} that it can also be described as the solution space of a PDE, as was done for the Radon transform \cite{john_radon, grinberg_radon, gonzalez_kakehi}.
	In the odd case, little appears to be known about the extension problem.

For even valuations, Theorem \ref{mthm:resProperty} is equivalent to the following.
\begin{Corollary}\label{mcor:imCosTransf}
	A function $f \in C^\infty(\Gr_k(\R^n))$ is in the image of the cosine transform $\mathcal C_k: C^\infty(\Gr_k(\R^n))\to  C^\infty(\Gr_k(\R^n))$, $1\leq k\leq n-2$, if and only if $f|_{Gr_k(E)}$ is in the image of $\mathcal C_k: C^\infty(\Gr_k(E))\to  C^\infty(\Gr_k(E))$ for all $k+2$-dimensional subspaces $E\subset\R^n$.
\end{Corollary} 

	In the proof of Theorem~\ref{mthm:resProperty}, we utilize the Alesker--Fourier transform on valuations to transform the problem into a dual statement concerning pushforwards of valuations, which can be solved using the representation by differential forms. An alternative proof for even valuations, using representation theory and composition series, appears in Appendix~\ref{app:repProofResThm}. 
	Theorem~\ref{mthm:resProperty} fails when the smooth valuations and sections in the assumption are replaced by continuous valuations and sections, as evidenced by the restrictions to $r$-subspaces of an $r$-homogeneous Klain--Schneider continuous valuation, a notion introduced in \cite{Bernig2017}. However, Theorem~\ref{mthm:resProperty} can be extended to Klain--Schneider continuous valuations, see Theorem \ref{thm:resPropertyCont}.

	\bigskip
	
	Theorems \ref{mthm:nashThm2} and \ref{mthm:resProperty_simple} are in fact special cases of the following general question. Given a set $S$ of linear subspaces of $\R^n$, and a family of valuations $\phi_E \in \Val_k^\infty(E)$, $E\in S$, of degree $0 < k < n$, does a globally defined valuation $\Phi \in \Val_k^\infty(\R^n)$ exist such that $\phi_E = \Phi|_E$ for all $E\in S$?
	
	As before, compatibility under restrictions, namely the requirement that for all $E, E' \in S$ one has $\phi_E|_{E \cap E'} = \phi_{E'}|_{E \cap E'}$, is a necessary condition. When $S$ is not discrete, we must further require the assignment $E\mapsto \phi_E$ to be smooth. We are interested in conditions on $S$ guaranteeing that those necessary conditions are also sufficient. 
	
	The problem can be viewed as an analogue of Whitney's extension problem~\cite{Whitney1934} for smooth valuations. Moreover, valuations are often easy or natural to describe in terms of their lower dimensional restrictions, e.g.\ through their Klain section. The problem gains further relevance if we recall that many important examples of valuations on linear spaces, such as the Euclidean intrinsic volumes, Hermitian intrinsic volumes~\cite{Bernig2011, Bernig2022} and Holmes--Thompson intrinsic volumes~\cite{bernig_ht}  satisfy a linear Weyl principle, namely behave naturally under restrictions.	

	\medskip
	Since the cases of $S=\Gr_k(\R^n)$, as well as $S\subset\Gr_k(\R^n)$ a perfectly self-avoiding submanifold, have been the subject of Theorems \ref{mthm:resProperty_simple} and \ref{mthm:nashThm2}, respectively, let us now consider finite arrangements of subspaces. Thus we are given finitely many subspaces $E_1, \dots, E_N$ of $\R^n$, possibly of different dimension, and compatible smooth valuations $\phi_1, \dots, \phi_N$ defined on these subspaces. 
	
Examples show that compatibility alone is in general too weak a condition for the existence of a global valuation (see Section~\ref{sec:finArr}). For an exception to this rule, asserting that the Klain section of an even valuation can assume arbitrary values on any finite subset of the grassmannian, see Proposition~\ref{prop:cosine_arbitrary}.

Thus we will be making geometric assumptions on the arrangement, all of which amount to restrictions on the possible dimensions of the various intersections. These will allow us to solve an extension problem for linear forms, similar to the one under consideration for valuations, that comes up routinely in the proof. Not less important, those assumptions allow to interpolate a family of compatible smooth functions defined on the spheres in $E_i$, $i=1, \dots, N$, by a smooth function on the sphere $S^{n-1}$ in $\R^n$.

 In order to state the assumptions, we introduce some definitions.
	\begin{Definition}\label{def:MinIntersSemiGeneric}
		Suppose that $E_1, \dots, E_N \subset \R^n$ is an arrangement of subspaces of $\R^n$. Then the arrangement is called
		\begin{itemize}
			\item \emph{minimally intersecting}, if for every non-empty $I \subseteq \{1, \dots, N\}$
			\begin{align*}
				\mathrm{codim} \bigcap_{i \in I} E_i = \sum_{i \in I} \mathrm{codim} E_i;
			\end{align*}
			\item \emph{semi-generic}, if $\{E_i\}_{i\in I}$ is minimally intersecting within $\sum_{i\in I}E_i$, whenever $\cap_{i\in I}E_i\neq \{0\}$ for some non-empty $I\subseteq\{1,\dots, N\}$.
		\end{itemize}
	\end{Definition}
	For example, any $\{E_1, E_2, E_3\}$ with $E_1\cap E_2\cap E_3=\{0\}$ is semi-generic. Let us point out that an arrangement $E_1, \dots, E_N$ is minimally intersecting if and only if the sum of the annihilators $E_1^\perp + \dots + E_N^\perp \subset (\R^n)^\ast$ is direct. 
	We will usually need the property of minimal intersection for tangent spaces of intersecting subspheres. This is where the notion of semi-genericity comes into play: 
	The arrangement $\{E_i\}_{i=1}^N$ is semi-generic, if and only if for any common point $\xi$ of a subcollection of intersecting spheres $S(E_i)\subset S^{n-1}$, $i\in I$, their tangent spaces $T_\xi S(E_i)$ intersect minimally.
	
	Our last main result now reads
	\begin{MainTheorem}\label{mthm:extensionFinArrang}
		Let $S=\{E_i\}_{i=1}^N$ be a semi-generic arrangement of proper subspaces in $\R^n$, and $k\geq 1$. Assume also one of the following:
		\begin{enumerate}
			\item Any $(k+1)$ subspaces from $S$ have at most $(k-1)$-dimensional intersection.
			\item Any  $2(k+1)$ subspaces or less from $S$ are minimally intersecting inside $E_1+\dots+E_N$.
			\item $k\geq\frac{n}{2}$.
			\item $N=3$.
		\end{enumerate}
		Suppose that $\phi_i\in\Val_{k}^{\infty}(E_i)$, $1\leq i  \leq N$, are given such that $\phi_i|_{E_i\cap E_j}=\phi_j|_{E_i\cap E_j}$, for all $1 \leq i,j \leq N$. Then there exists $\phi\in\Val^\infty_k(\R^n)$ with $\phi|_{E_i} = \phi_i$, $1 \leq i \leq N$.
	\end{MainTheorem}
The two main cases of the theorem are of complementary nature. The first condition is increasingly lax as $k$ increases: if satisfied for $k$, it is also satisfied for $k'>k$, in particular it is always satisfied for $k\geq\frac n2$; it is most easily satisfied by low dimensional subspaces. In contrast, the second condition is increasingly stringent as $k$ increases, and is never satisfied for $k\geq \frac{n-1}{2}$; all but one subspace must have dimension at least $n/2$. Taken together, the two cases imply the following.
	\begin{Corollary}\label{cor:hyperplanes}
		When $k\neq \frac{n-1}{2}$, the conclusion of Theorem \ref{mthm:extensionFinArrang} holds for any set $S$ of hyperplanes in general position.
	\end{Corollary}The conditions in i) and ii) are in some sense tight, as evidenced by Example~\ref{example:counterexample} showing that Corollary \ref{cor:hyperplanes} fails for $k=\frac{n-1}{2}$.
	
	In the proof of Theorem~\ref{mthm:extensionFinArrang}, we will again work with an Alesker--Fourier-dual formulation and differential forms. The problem then reduces to solving an extension problem for differential forms. However, as the representation of a smooth valuation by differential forms is not unique, an additional alignment-step is necessary, which is in fact the main step in the proof. 
	
	It would be interesting to extend the results and methods of Theorem \ref{mthm:extensionFinArrang} to submanifolds $Z\subset\Gr_k(\R^n)$ more general than those in Theorem \ref{mthm:nashThm2}. In fact, the proof ot the latter amounts to an adaptation of the the proof of Theorem \ref{mthm:extensionFinArrang} in the simplest case of i) when every two subspaces intersect trivially.
	
	
	\subsection*{Proof ideas and organization of the article}
	We give a short overview of the proof ideas for the main theorems and comment on the structure of the article. The proofs in the three main sections~\ref{sec:extFullGrass} to \ref{sec:Nash} are formally independent and can be read separately. The proof of Theorem~\ref{mthm:croftonFormulaMf} in Section~\ref{sec:crofton} depends on Theorem~\ref{mthm:nashThm}, but not on the techniques of the proof. In Section~\ref{sec:bg}, we give some background and preparatory steps on double forms (needed in Section~\ref{sec:extFullGrass}) and valuations (needed throughout).
	A key ingredient in most proofs is a preparatory application of the Alesker--Fourier transform. This is particularly useful when representing valuations by differential forms. While the Alesker--Fourier transform interchanges the restriction of  valuations with the conceptually more involved push-forward under projection, the effect is reversed on the corresponding differential forms, and so instead of solving a complicated system of integral equations, we find ourselves in a Whitney-type extension problem for differential forms.
	A similar effect occurs also when representing valuations using the surface area measure -- this is exploited in the proof of Theorem \ref{mthm:nashThm} in Appendix~\ref{app:NuclearNash}.

	\subsubsection*{Nash-type embedding theorem -- proof of Theorems~\ref{mthm:nashThm} and \ref{mthm:nashThm2}}
	We give two different proofs of the Nash-type embedding Theorem~\ref{mthm:nashThm}, both depending on the existence of a suitable generic embedding, which we establish in Lemma~\ref{lem:exPerfNonParEmbed} using Thom's transversality theorem.
	
	In Section~\ref{sec:Nash}, we prove Theorem~\ref{mthm:nashThm2} by first applying the Alesker--Fourier transform and then choosing a smooth family of differential forms that represent the given assignment of valuations. As there are no compatibility constraints, this family can be extended to a globally defined differential form, yielding the claimed valuation. Theorem~\ref{mthm:nashThm} is then deduced as a straightforward corollary of Theorem \ref{mthm:nashThm2}. 
	
	In Appendix~\ref{app:NuclearNash}, we prove Theorem~\ref{mthm:nashThm} using the Alesker product. We use a natural filtration of $\mathcal{V}^\infty(M)$ where quotients of subsequent subspaces are isomorphic to fields of smooth homogeneous valuations on the tangent spaces of $M$. We then solve the extension problem for such fields in the $1$-homogeneous case and lift this solution to general degrees of homogeneity by taking the Alesker product of the extensions, and utilizing the nuclearity of the Fr\'echet spaces involved.

	\subsubsection*{Crofton formulas -- proof of Theorem~\ref{mthm:croftonFormulaMf}}
	 
	To establish the Crofton formula from Theorem~\ref{mthm:croftonFormulaMf}, we first prove a Crofton formula for all translation-invariant valuations on a linear space. We then restrict it to the manifold by Theorem~\ref{mthm:nashThm}. In the linear construction (proved in Proposition~\ref{prop:exCroftMeasOddValonR} by a representation theoretic argument), we extend a representation for $1$-homogeneous valuations to arbitrary degrees of homogeneity by tracking its behavior under Alesker product.
	
	\subsubsection*{Extension from the full grassmanian -- proof of Theorems~\ref{mthm:resProperty_simple} and \ref{mthm:resProperty}}
	To prove Theorems~\ref{mthm:resProperty_simple} and \ref{mthm:resProperty}, we rewrite the problem using the Alesker--Fourier transform to an extension problem for double forms that uniquely determine the valuation. After extending the double forms (Section~\ref{sec:extFullGrassDoubleForms}) simultaneously in a smooth way, we then show that the extension defines a valuation.
	
	In Appendix~\ref{app:repProofResThm}, we give a different proof of Theorem~\ref{mthm:resProperty} for even valuations utilizing representation theory.
	
	\subsubsection*{Extension from finite arrangements -- proof of Theorem~\ref{mthm:extensionFinArrang}}
	Here, we first study extension problems for linear forms in Section~\ref{sec:finArrLinForms}. Dualizing the problem again by the Alesker--Fourier transform and representing the valuations by (non-unique) differential forms, we first modify the forms in order to make them compatible on intersections by adding suitably chosen closed forms (Proposition~\ref{prop:compatible_primitives}) and then construct a common extension, defining the claimed extension of the valuations. In the remainder of the section, we give Example~\ref{example:counterexample} showing the sharpness of the conditions of Corollary~\ref{cor:hyperplanes}. We also present a more general extensibility statement for densities, which has a simple proof.

	\section{Background and preparations}\label{sec:bg}
	In this section we review the needed background on double forms and valuations, and do some basic preparatory steps. As a general reference on valuation theory we recommend \cite{Klain1997} and \cite{Schneider2014}, as well as \cite{Alesker2014c}, \cite{Alesker2007b} and \cite{Bernig2010} for the modern theory of valuations. Throughout the paper we denote by $\Gamma(X, \mathcal E)$ the $C^\infty$-smooth sections of a bundle $\mathcal E$ over $X$.

	\subsection{Double forms}\label{sec:prepDoubleforms}
	We start with a short introduction on double forms, which we will need (only) in the proof of Theorem~\ref{mthm:resProperty} in Section~\ref{sec:extFullGrass}. The techniques and statements of this section are not new. However, as we did not find references for the exact statements, we will also sometimes give proofs for the reader's convenience.
		
	 Let in the following $V$ be a linear space, and denote by $D_{p, q}(V)=\wedge^p V^*\otimes \wedge^q V^*$ the space of double forms of bi-degree $(p,q)$, $p,q \geq 0$, on $V$.
	
	Following the presentation in \cite{Gray1969}, the wedge product on $V$ and $V^\ast$ induces a pairing $\wedge: D_{p,q}(V)\otimes D_{p',q'}(V)\to D_{p+p',q+q'}(V)$, which is given on decomposable vectors by
	\begin{align*}
		(v \otimes w) \wedge (v' \otimes w') = (v \wedge v') \otimes (w \wedge w'),\quad v,v',w,w' \in \wedge^{\bullet} V^\ast.
	\end{align*}

	For $p=q$, we have the subspace $Y_p(V) = \Sym^2(\wedge^p V^*)\subset D_{p,p}(V)$ of symmetric forms, that is, forms satisfying $Q(v,w) = Q(w,v)$, for every $v,w \in \wedge^p V$. Clearly, the wedge product restricts to a pairing $\wedge: Y_{p}(V)\otimes Y_q(V)\to Y_{p+q}(V)$.
	
	For every $Q \in Y_p(V)$, there exists a naturally associated quadratic form on $\wedge^p V$,  namely $v\mapsto Q(v,v)$, $v \in \wedge^p V$. In the following, we will often abuse notation and will not distinguish between $Q \in Y_p(V)$ and its quadratic form.
	
	\medskip
	
	Next, consider the subspace $Z_p(V) \subseteq Y_p(V)$, consisting of all forms $Q \in Y_p(V)$ whose associated quadratic forms vanish on all decomposable vectors $v \in \wedge^p V$, that is, those $Q$ satisfying
	\begin{align*}
		Q(v_1 \wedge \dots \wedge v_p, v_1 \wedge \dots \wedge v_p) = 0, \quad \forall v_1, \dots, v_p \in V. 
	\end{align*}
	The generators of $Z_p(V)$ are well-known (see, e.g., \cite{Griffiths1994}*{p.209ff}) and can be described using the contraction $i_v : \wedge^p V^\ast \to \wedge^{p-k}V^\ast$, for $p \geq k$ and $v\in \wedge^k V$, 
	\begin{align*}
		\langle  i_v \xi, w \rangle=\langle \xi, v\wedge w\rangle,\qquad\forall \xi\in\wedge^p V^\ast, w\in \wedge^{p-k}V.
	\end{align*}

	\begin{Proposition}\label{prop:zpSpannedPluckerRelations}
		$Z_p(V)$ is spanned by the quadratic forms $P_{\xi,\eta}(v):=\langle i_v\xi \wedge \eta, v\rangle$, $v \in \wedge^p V$, where $\xi\in\wedge^{p+1}V^*$ and $\eta\in\wedge^{p-1}V^*$ are decomposable.
	\end{Proposition}

	\medskip

	Using that $Y_p(V)^\ast = Y_p(V^\ast)$, we may define $A_p(V):=Z_p(V^*)^\perp\subseteq Y_p(V)$. 
	
	\begin{Lemma}\label{lem:1d_to_Ak}
		Let $p \geq 1$ and $Q\in\Sym^2 V^*$. Then $Q^{\wedge p}\in A_p(V)$.
	\end{Lemma}
	\proof
	We may choose functionals $\eta_i\in V^*$, $1\leq i\leq r$, such that $Q=\sum_{i=1}^r \epsilon_i \eta_i \otimes \eta_i $ for some $\epsilon_i\in\{\pm 1\}$. Hence there are constants $c_I \in \Z$ such that $Q^{\wedge p}=\sum_I c_I \eta_I\otimes \eta_I$, where we sum over all $I=\{1\leq i_1<\dots<i_p\leq r\}$, and $\eta_I=\eta_{i_1}\wedge\dots\wedge \eta_{i_p}$. 
	Let $Q^*\in Z_p(V^*)$ be arbitrary. Since $\eta_I \in \wedge^p V$ is a decomposable vector, $Q^*(\eta_I, \eta_I) = 0$ for every $I$, and so $\langle Q^*, Q^{\wedge p}\rangle =\sum_I c_I Q^*(\eta_I, \eta_I)=0,$
	that is, $Q^{\wedge p}\in A_p(V)$.
	\endproof
	
	Note that if $P\in \Sym^2 V^*$ is a Euclidean structure, then $P^{\wedge p}\in A_p(V)$ is the induced Euclidean structure on $\wedge^p V$. 
	
	\medskip
	
	Next, observe that both $Z_p(V)$ and $A_p(V)$ are invariant subspaces with respect to the standard $\GL(V)$-representation on $Y_p(V)$, that is,
	\begin{align*}
		g \cdot Q(v,w) = Q(g^{-1}(v), g^{-1}(w)),
	\end{align*}
	for all $g \in \GL(V), Q \in Y_p(V)$ and $v,w \in \wedge^p V$. Indeed, $Y_p(V)$ decomposes into a direct sum of $Z_p(V)$ and $A_p(V)$, as the following lemma shows. Moreover, it is well-known that $Y_p(V)/Z_p(V)$ is an irreducible $\GL(V)$-module (see, e.g., \cite[Ex.~15.43]{Fulton1991}), and, consequently, $A_p(V) \cong Y_p(V)/Z_p(V)$ is irreducible.
	
	\begin{Lemma}\label{lem:plucker_decomposition}
		Let $p \geq 1$. Then the $GL(V)$-module $Y_p(V)$ decomposes into 
		\begin{align*}
			Y_p(V) = Z_p(V)\oplus A_p(V).
		\end{align*}
		\noindent Moreover, $A_p(V)$ is the only $\GL(V)$-invariant complement of $Z_p(V)$ in $Y_p(V)$.
	\end{Lemma}
	\proof
	Since $\left(A_p(V)+Z_p(V)\right)/Z_p(V)$ is an invariant submodule of the irreducible $\GL(V)$-module $Y_p(V)/Z_p(V)$, it suffices to show that $A_p(V)$ is not contained in $Z_p(V)$ in order to prove $Y_p(V) = A_p(V) + Z_p(V)$.
	
	To this end, let $P$ be a Euclidean structure on $V$. By Lemma~\ref{lem:1d_to_Ak} and the remark below it, $P^{\wedge p}\in A_p(V)$ and $P^{\wedge p}$ is positive definite on $\wedge^p V$, that is, $P^{\wedge p}(v,v) > 0$ for all $v \in \wedge^p V$. Hence, $P^{\wedge p} \in A_p(V) \setminus Z_p(V)$, and so $Y_p(V) = A_p(V) + Z_p(V)$.
	
	Noting that $\dim A_p(V) = \dim Y_p(V^*) - \dim Z_p(V^*)$ and that $\dim Y_p(V)$ and $\dim Z_p(V)$ only depend on $\dim V$, we deduce $\dim Y_p(V) = \dim Z_p(V) + \dim A_p(V)$ and, therefore, $Y_p(V) = A_p(V) \oplus Z_p(V)$.
	
	\medskip

	For the second statement, we recall that $\SL(V)$ is semi-simple and $A_p(V)$ and $Z_p(V)$ are its representations. Thus $Z_p(V)$ decomposes into a direct sum of irreducible $\SL(V)$ modules. It then suffices to see that $Z_p(V)$ does not contain a representation of $\SL(V)$ isomorphic to $A_p(V)$.  But this is clear since $Y_p=\Sym^2(\wedge ^p V^*)$ has multiplicity one under the action of $\SL(V)$ (see, e.g., \cite{Fulton1991}*{Ex.~15.32}).
	\endproof

	The space $A_p(V)$ can be described more directly using the operation $\omega \mapsto \omega'$, $D_{p,q}(V)\mapsto D_{p+1,q-1}(V)$, which is defined, following \cite{Gray1969}*{Eq.~(2.6)}, by
	\begin{align*}
		\omega'(v_1\wedge\dots\wedge v_{p+1}, w_2\wedge\dots\wedge w_q)=\sum_{j=1}^{p+1}(-1)^{j+1}\omega(\wedge_{i\neq j}v_i, v_j\wedge w_2\wedge\dots\wedge w_q),
	\end{align*}
	where $v_1, \dots, v_{p+1}, w_2, \dots, w_q \in V$. Note that, by \cite{Gray1969}*{Prop.~2.1}, this operation satisfies a Leibnitz rule, that is, 
	\begin{align}\label{eq:Leibnitz}
		(\omega \wedge \eta)' = \omega' \wedge \eta + (-1)^{p+q} \omega \wedge \eta',
	\end{align}
	where $\omega \in D_{p,q}(V)$ and $\eta \in D_{r,s}(V)$.
	
	The following description essentially appeared in \cite{Gray1969}.
	\begin{Proposition}\label{prop:description_Ap} Let $p \geq 1$. Then
		\begin{align*}
			A_p(V)&=\Span\{Q_1\wedge\dots\wedge Q_p: Q_1,\dots, Q_p\in \Sym^2(V^*)\}\\
			&=\{Q\in Y_p(V): Q'=0\}.
		\end{align*}
	\end{Proposition}
	\proof
	Writing $B_p(V)=\Span\{Q_1\wedge\dots\wedge Q_p: Q_1,\dots, Q_p\in \Sym^2(V^*)\}$ and $C_p(V)=\{Q\in Y_p(V): Q'=0\}$, we will show that $A_p(V) \subseteq B_p(V) \subseteq C_p(V) \subseteq A_p(V)$, which will yield the claim.
	
	First, fix a Euclidean structure $P \in \Sym^2V^*$. Then, $P^{\wedge p}$ is non-zero (indeed, it is a Euclidean structure on $\wedge^p V$) and, by Lemma~\ref{lem:1d_to_Ak}, $P^{\wedge p} \in A_p(V) \cap B_p(V)$. Thus, $A_p(V) \cap B_p(V)$ is a non-zero $\GL(V)$-submodule of the irreducible module $A_p(V)$, that is, $A_p(V) \cap B_p(V) = A_p(V)$, which implies $A_p(V) \subseteq B_p(V)$.
	
	Second, since clearly $Y_1(V)=C_1(V)$, \eqref{eq:Leibnitz} implies that $B_p(V) \subseteq C_p(V)$.
	
	Finally for the inclusion $C_p(V) \subseteq A_p(V)$, assume that $Q \in C_p(V)$, that is $Q\in Y_p(V)$ with $Q' = 0$
	and let $Z \in Z_p(V^\ast)$. We need to show that $\langle Q, Z\rangle = 0$. By Proposition~\ref{prop:zpSpannedPluckerRelations}, it suffices to consider $Z=P_{v,w}$, where $v = v_1 \wedge\dots\wedge v_{p+1} \in \wedge^{p+1} V$ and $w=w_2\wedge \dots\wedge w_p \in \wedge^{p-1} V$ are decomposable.
	
	Writing $Q = \sum_{i \in I} c_i \eta_i \otimes \xi_i$, for some $\eta_i, \xi_i \in \wedge^p V^\ast$, we first need to calculate
	\begin{align*}
		\langle \eta_i \otimes \xi_i, P_{v,w}\rangle = \frac{1}{2}\left(\langle \eta_i, i_{\xi_i} v \wedge w\rangle + \langle \xi_i, i_{\eta_i} v \wedge w\rangle\right).
	\end{align*}
	The first term becomes
	\begin{align}
		\langle \eta_i, i_{\xi_i} v \wedge w\rangle &= \sum_{l=1}^{p+1}(-1)^{l+1}  \eta_i(v_l \wedge w) \xi_i(v_1 \wedge \dots\wedge \widehat v_l \wedge\dots \wedge v_{p+1}) \nonumber\\&= (\xi_i \otimes \eta_i)'(v, w),\label{eq:prfLemApBpCp1}
	\end{align}
	where we implicitly used a Laplace expansion to calculate $i_{\xi_i} v$. The second term is analogous, yielding $(\eta_i \otimes \xi_i)'(v, w)$. Summing up all equations~\eqref{eq:prfLemApBpCp1} and using that by symmetry, $Q =\sum_{i \in I} c_i \eta_i \otimes \xi_i =  \sum_{i \in I} c_i \xi_i \otimes \eta_i$, we obtain
	\begin{align*}
		\langle Q, P_{v,w}\rangle = \sum_{i \in I} c_i \langle \eta_i \otimes \xi_i, P_{v,w}\rangle = \frac{1}{2}\left(Q'(v,w) + Q'(v,w)\right) = 0.
	\end{align*}
	Consequently, $Q \in Z_p(V^\ast)^\perp = A_p(V)$, which concludes the proof.
	\endproof

	\subsection{Valuations}\label{sec:bgval}
		We turn now to valuations, ways to represent them and operations on them. While we are mostly interested in smooth valuations, let us first recall the classical notion of valuations on convex bodies.
		\begin{Definition}
			 A functional $\phi: \mathcal{K}(V) \to \R$ on the space $\mathcal{K}(V)$ of compact convex subsets (convex bodies) of a linear space $V$ is called a \emph{valuation} if
			 \begin{align*}
			 	\phi(K \cup L) + \phi(K \cap L) = \phi(K) + \phi(L), \quad \forall K, L, K \cup L \in \mathcal{K}(V).
			 \end{align*}  
		\end{Definition}
		We denote by $\Val(V)$ the space of all valuations on $V$ that are translation-invariant and continuous with respect to the Hausdorff metric. The subspace of smooth vectors with respect to the action of $\GL(V)$ on $\Val(V)$, defined by $(\eta \cdot \phi)(K) = \phi(\eta^{-1}K)$, $\eta \in \GL(V)$, $\phi \in \Val(V)$, coincides with the space $\Val^\infty(V)$, defined in the introduction (see \cite{Alesker2006b}).
			
		\medskip
		
		In the following, we will mostly restrict our presentation to smooth trans\-lation-invariant valuations on a linear space $V$. We remark that most notions apply equally to valuations on manifolds. The cosphere bundle of $V$ is $\mathbb P_V = V \times \mathbb{P}_+(V^\ast)$, where the oriented projectivization $\mathbb{P}_+ (W)$ is the space of oriented $1$-dimensional linear subspaces of a linear space $W$.
			By $\Dens(W)$ we denote the $1$-dimensional space of densities, that is Lebesgue measures, on $W$.
		
		\subsubsection{Smooth valuations and differential forms} \label{sec:bgValDiffForm}
		Recall that, fixing an orientation on a linear space $V$, a valuation $\phi \in \Val^\infty(V)$ is given by a pair of translation-invariant differential forms $\omega \in \Omega^{n-1}(\mathbb P_V)^{tr}$ and $\theta = c \vol_V \in \Omega^n(V)^{tr}$, where $\vol_V$ is a volume form on $V$, and
		\begin{align}\label{eq:smoothValVec}
			\phi(A) = \int_{\nc(A)} \omega + c \vol_V(A), \qquad A \in \mathcal{P}(V).
		\end{align}
		Specifying an orientation at this point is actually not necessary if we instead consider $\omega \in\Omega^{n-1}(\mathbb P_V)^{tr} \otimes \ori(V)$ and $\theta \in \Omega^n(V)^{tr}\otimes \ori(V)$, where $\ori(V)$ is the orientation bundle of $V$, that is, the space of all functions $\rho: \wedge^n V \to \R$, where $n = \dim V$, satisfying $\rho(\lambda \tau) = \sign(\lambda) \rho(\tau), \lambda \in \R, \tau \in \wedge^n V$.  We remark that \eqref{eq:smoothValVec} can serve as a definition for all sets $A$ admitting a conormal cycle, see \cite{Alesker2009, Federer1959,Fu1994}.
		
		The product structure of $\mathbb P_V= V \times \mathbb{P}_+(V^\ast)$ induces a bi-grading on $\Omega(\mathbb P_V)$, and we denote by $\Omega^{k,l}(V \times \mathbb{P}_+(V^\ast))$ the subspace of forms of bi-degree $(k,l)$. A smooth valuation $\phi$ is $k$-homogeneous exactly if it can be represented by $(\omega, \theta)$, where $\omega$ has bi-degree $(k, n-1-k)$ and $\theta = 0$, for $0 \leq k < n$, or where $\omega = 0$, for $k=n$.
		
		Next, note that the representation \eqref{eq:smoothValVec} is not unique and the kernel was described in \cite{Bernig2007} using the Rumin differential $D$, giving rise to a unique representation by an $n$-form on the cosphere bundle. To state the result, observe that the cosphere bundle $\mathbb P_V$ is a contact manifold, and fix a contact form $\alpha$. A form $\omega \in \Omega^\bullet(\mathbb P_V)$ is called vertical if $\alpha \wedge \omega = 0$, which is equivalent to $\omega = \alpha \wedge \tau$ for some form $\tau$.

		\begin{Theorem}[Bernig--Br\"ocker \cite{Bernig2007}]\label{thm:bernigBroecker}
			Suppose that $1 \leq k \leq n-1$. Then the map
			\begin{align*}
				\Val_k^\infty(V) \to \Omega^{k,n-k}(V\times\mathbb P_+(V^*))^{tr}\otimes \ori(V), \quad \phi \mapsto \tau_\phi,
			\end{align*}
			given by $\tau_\phi = D\omega$, where $\phi(A)=\int_{\nc(A)}\omega$ with $\omega\in\Omega^{k,n-1-k}(V\times\mathbb P_+(V^*))^{tr}\otimes \ori(V)$, is injective. For $2 \leq k \leq n-1$, its image is given by the subspace of closed and vertical forms, whereas for $k=1$, its image consists of all closed and vertical $\tau$ satisfying $\pi_*\tau=0$, where $\pi: V\times\mathbb P_+(V^*)\to V$ is the natural projection.
		\end{Theorem}

		We will call $\tau_\phi$ the defining form of $\phi$.

		In the proof of Theorem~\ref{mthm:resProperty}, we will extensively use the following natural identifications induced by translation-invariance and verticality of the forms and the Hodge-star operator (see also \cite[Sec.~4]{Bernig2016}). First, by translation-invariance
		\begin{align*}
			\Omega^{k,n-k}(V\times\mathbb P_+(V^*))^{tr}\otimes\ori(V)\cong\Omega^{n-k}(\mathbb P_+(V^*), \wedge^k V^*\otimes\ori(V))\\
		= \Gamma(\mathbb P_+(V^*), \wedge^{n-k}T_\xi^*\mathbb P_+(V^*)\otimes \wedge^k V^*\otimes\ori(V)),
		\end{align*}
		and, since 
		\begin{align*}
			\wedge^{n-k}T_\xi^*\mathbb P_+(V^*)\cong\wedge^{n-k}(\xi^*\otimes V^*/\xi)^*\cong \xi^{\otimes (n-k)}\otimes\wedge^{n-k}\xi^\perp,
		\end{align*}
		and
		\begin{align*}
			\wedge^k V^*\otimes\ori(V)\cong\wedge^{n-k}V\otimes\Dens(V),
		\end{align*}
		through the Hodge star operator, we obtain
		\begin{align}\label{eq:forms_sections}
			&\Omega^{k,n-k}(V\times\mathbb P_+(V^*))^{tr}\otimes\ori(V) \nonumber\\
			&\qquad\qquad\cong\Gamma(\mathbb P_+(V^*), \wedge^{n-k}\xi^\perp\otimes  \wedge^{n-k}V\otimes \xi^{\otimes (n-k)}\otimes\Dens(V)).
		\end{align}
		We will slightly abuse notation and switch between those spaces as needed. In particular, when a Euclidean structure and orientation are fixed on $V$, these spaces reduce to $\Omega^{n-k}(S(V^*), \wedge^k V^*)\cong \Omega^{n-k}(S(V^*), \wedge^{n-k}V)$. Here and throughout, we denote by $S(W)$ the sphere in a linear space $W$, whenever $W$ is endowed with a Euclidean structure.
		
		Note that the subspace of vertical forms corresponds to the subspace of sections
		\begin{align*}
			\Gamma(\mathbb P_+(V^*), \wedge^{n-k}\xi^\perp\otimes \wedge^{n-k}\xi^\perp\otimes  \xi^{\otimes(n-k)}\otimes \Dens(V)),
		\end{align*} 
		and so we will call such sections \emph{vertical} as well.

	A closed vertical form must satisfy an additional pointwise constraint, as follows, which was observed in \cite{Faifman2023} and \cite{bernig_personal}. 
	\begin{Proposition}\label{prop:pwConstraintApSecVal}
		Suppose that $\tau\in \Omega^{k,n-k}(V\times\mathbb P_+(V^*))^{tr}\otimes\ori(V)$, and let $Q\in \Gamma(\mathbb P_+(V^*), \wedge^{n-k}\xi^\perp\otimes  \wedge^{n-k}V\otimes \xi^{\otimes(n-k)}\otimes \Dens(V))$ be the corresponding section under \eqref{eq:forms_sections}. If $\tau$ is vertical and closed, then 
		$$Q\in  \Gamma(\mathbb P_+(V^*), A_{n-k}(V^*/\xi)\otimes\xi^{\otimes(n-k)}\otimes \Dens(V)).$$
	\end{Proposition}
	\proof
	The fact that $Q$ is symmetric, namely $Q|_\xi\in  Y_{n-k}(V^*/\xi)\otimes\xi^{\otimes(n-k)}\otimes \Dens(V)$, appears in \cite{Faifman2023}, where it is proved by showing that $(Q|_\xi)'=0$. Proposition~\ref{prop:description_Ap} then immediately implies the statement.
	\endproof
	
	\medskip
	
	By the representation of a smooth valuation from \eqref{eq:forms_sections}, we deduce from Proposition~\ref{prop:pwConstraintApSecVal} that every $\phi\in \Val^\infty_k(V)$ defines a section of quadratic forms 
	\begin{align*}
		Q(\phi) \in \Gamma(\mathbb P_+(V^*), A_{n-k}(V^*/\xi)\otimes\xi^{\otimes(n-k)}\otimes \Dens(V)).
	\end{align*}
	We will use the same notation also for twisted valuations in $\Val^\infty_k(V)\otimes\Dens(V)^*$.

	\subsubsection{Klain and Schneider embedding}
	
	According to a theorem of Hadwiger, $\Val_n(\R^n)$ is the one dimensional space $\Dens(\R^n)$ of Lebesgue measures. This fact is used to define, through restrictions, the Klain map $\Kl:\Val_j^{+,\infty}\to \Gamma(\Gr_j(\R^n), \Dens(E))$ on even valuations, see eq. \eqref{eq:klain}.  Klain \cite{Klain1997, Klain2000} has shown it to be injective.

	For general valuations, restrictions to $(j+1)$-dimensional subspaces yield an equivariant injective map $\Sc: \Val_j^{\infty}(V) \to \Gamma(\widehat{\mathcal{F}}_{j}^{j+1}, \Dens(E))/L$, where $\widehat{\mathcal{F}}_{j}^{j+1}$ is the partial flag manifold of cooriented pairs $(E,F)$, $F \in \Gr_{j+1}(V)$ and $E \in \Gr_j(F)$, and $L$ denotes the (closed) subspace corresponding to linear functions. More precisely, $L$ consists of all sections $\phi \in \Gamma(\widehat{\mathcal{F}}_{j}^{j+1}, \Dens(E))$ such that after trivializing $\Gamma(\widehat{\mathcal{F}}_{j}^{j+1}, \Dens(E)) \cong C^\infty(\widehat{\mathcal{F}}_j^{j+1})$ by a Euclidean structure on $V$, the function $S(F) \ni E \mapsto \phi(E, F)$ is the restriction to $S(F)$ of a linear function, for every fixed $F \in \Gr_{j+1}(V)$.
	
	Note that the part of $\Sc_ \phi$, $\phi \in \Val_j^\infty(V)$, which is even with respect to the change of coorientation, is essentially the pullback of the Klain section of the even part of $\phi$, under the natural projection $\widehat{\mathcal{F}}_{j}^{j+1}\to\Gr_j(V)$, divided by a factor of $2$. Moreover, when restricted to odd valuations, $\Sc$ maps into the subspace $\Gamma^-(\widehat{\mathcal{F}}_{j}^{j+1}, \Dens(E))/L$ of odd sections.
	
	\medskip
		
	\subsubsection{Operations on valuations}\label{sec:bgopVal} Next, we recall some useful operations on smooth valuations, which we will need later on.

	First, we review some of the functorial properties of smooth valuations on manifolds, see \cite{Faifman2017, Alesker2011b, Alesker2009}, starting with the notion of pullback by an embedding.

	\begin{Theorem}[{\cite[Claim~3.1.1]{Alesker2009}}]
		Suppose that $M,N$ are smooth manifolds and let $f:M \to N$ be a smooth embedding. Then there exists a continuous linear map
		\begin{align*}
			f^\ast: \mathcal{V}^\infty(N) \to \mathcal{V}^\infty(M),
		\end{align*}
		called the \emph{pullback by $f$}, given by $(f^\ast \phi)(A) = \phi(f(A))$, $A \in \mathcal{P}(M)$.
	\end{Theorem}

 When $f$ is an inclusion map, we often call $f^\ast \phi$ the restriction of $\phi$ to $M$, and denote it by $\phi|_M$. We remark that pullbacks can be defined for larger classes of maps and valuations. 
	
	Let us also recall the definition of pullback by surjections, but only in the linear setting.
	\begin{Definition}
		Given an epimorphism $\pi:V\to W$, the \emph{pullback} $\pi^*:\Val(W)\to\Val(V)$ is given by $\pi^*\phi(K)=\phi(\pi(K))$ for every convex body $K \subset V$.
	\end{Definition}
	
		The space of smooth valuations (on a manifold) can be endowed with a product structure, called the Alesker product.
	
	\begin{Theorem}[{\cite{Alesker2008, Alesker2012}}]
		Suppose that $M$ is a smooth manifold. Then there exists a natural commutative product $\mathcal{V}^\infty(M) \times \mathcal{V}^\infty(M) \to \mathcal{V}^\infty(M)$, such that the pullback is an algebra homomorphism, that is,
		\begin{align*}
			f^\ast (\phi \cdot \psi) = (f^\ast \phi) \cdot (f^\ast \psi), \quad \phi, \psi \in \mathcal{V}^\infty(N),
		\end{align*}
		for every smooth embedding $f: M \to N$. Moreover, the Euler characteristic $\chi$ is the multiplicative identity, that is $\chi \cdot \phi = \phi$ for all $\phi \in \mathcal{V}^\infty(M)$.
	\end{Theorem}
	Denoting by $\mathcal{V}_c^\infty(M)$ the subspace of $\mathcal{V}^\infty(M)$ of valuations with compact support, the product composed with evaluating at the manifold itself yields a perfect pairing $\mathcal{V}^\infty(M) \times \mathcal{V}_c^\infty(M) \to \R$, that is, an injective map $\mathcal{V}^\infty(M) \to \mathcal{V}_c^\infty(M)^\ast$ with dense image, called Poincar\'e duality (see \cite{Alesker2014c}). The space of generalized valuations is defined by $\mathcal{V}^{-\infty}(M) = \mathcal{V}_c^\infty(M)^\ast$, and by the above naturally contains $\mathcal{V}^{\infty}(M)$.
	
	For translation-invariant, smooth valuations on a linear space $V=\R^n$ we similarly have a Poincar\'e duality \cite{alesker_multiplicative} with values in $\Val_n^\infty(V) = \Dens(V)$, namely \[\Val^\infty(V)\times \Val^\infty(V)\to\Dens(V), \quad \langle \phi, \psi\rangle=(\phi\cdot\psi)_n,\]
	where $\phi_n$ denotes the $n$-homogeneous component of $\phi \in \Val^\infty(V)$. Moreover, it extends to a non-degenerate pairing $\Val^\infty(V)\times \Val(V)\to\Dens(V)$ (see \cite{Alesker2011b}). 
	 Hence we obtain an injective map $\Val_j(V) \to \Val_{n-j}^\infty(V)^\ast \otimes \Dens(V)$ with $\Val_j^\infty(V)$ dense in the image, and we define $\Val_j^{-\infty}(V) = (\Val_{n-j}^\infty(V)^\ast \otimes \Dens(V))^\ast$.
	
	The pushforward of valuations is in general dual to the pullback. We will only define it for surjective linear maps between linear spaces and translation-invariant valuations, since this is the only case we will use. We remark that there are two closely related but distinct constructions of pushforward, the other being the pushforward of valuations on manifolds by a proper map.
	\begin{Theorem}[{\cite[Prop.~3.2.1]{Alesker2011b}}]
		Suppose that $V,W$ are linear spaces and let $f:V \to W$ be linear and surjective, and $d=\dim V-\dim W$. Then there exists a continuous linear map
		\begin{align*}
			f_\ast: \Val^\infty_k(V) \otimes \Dens(V)^\ast \to \Val^\infty_{k-d}(W) \otimes \Dens(W)^\ast,
		\end{align*}
		called the \emph{pushforward by $f$}, which is formally adjoint to pull-back: 
		\[ \langle f_*\phi, \psi\rangle =\langle \phi, f^*\psi\rangle, \quad \forall\phi\in\Val^\infty(V)\otimes\Dens(V)^*, \psi\in\Val^\infty(W).\]
	\end{Theorem}

	The pushforward under a linear surjection can easily be expressed using the quadratic form of  valuations, as the following proposition shows.
	
	\begin{Proposition}\label{prop:ConvPushfOnQuadrForms}
		Let $\phi \in \Val_k^\infty(V)\otimes \Dens(V)^*$ be a twisted $k$-homogeneous valuation, and denote by $Q(\phi)\in\Gamma(\mathbb P_+(V^*), A_{n-k}(V^*/\xi)\otimes\xi^{\otimes(n-k)})$ its corresponding section of quadratic forms. Let $\pi: V \to W$ be  a linear surjection. Then
		\begin{align*}
			Q(\pi_\ast \phi)|_\xi = (\pi^\vee)^\ast \left( Q(\phi)|_{\pi^\vee (\xi)} \right) , \quad \xi \in \mathbb{P}_+(W^\ast),
		\end{align*}
	 where $\pi^\vee$ denotes the adjoint map.
	\end{Proposition}
	\proof
	 The statement follows by the argument in \cite[Prop.~3.2.3]{Alesker2009}, adapted to the linear setting.  
	\endproof

	Next, denote by $\mathbb{F}_V$, or just $\mathbb{F}$, the Alesker--Fourier transform, introduced in \cite{Alesker2003} for even smooth valuations and in \cite{Alesker2011b} for general smooth valuations.
	\begin{Theorem}[\cite{Alesker2011b}]
		Suppose that $V$ is a vector space. Then there exists an isomorphism of linear topological spaces
		\begin{align*}
			\mathbb{F}_V: \Val^\infty(V) \to \Val^\infty(V^\ast) \otimes \Dens(V),
		\end{align*}
		which commutes with the natural action of the group $\GL(V)$ on both spaces.
	\end{Theorem}
	Note that $\mathbb{F}$ interchanges degree and codegree, that is, if $\phi \in \Val^\infty_j(V)$, then $\mathbb{F} \phi \in \Val_{n-j}^\infty(V^\ast)\otimes \Dens(V)$. Moreover, it interchanges pushforward and pullback.
	\begin{Theorem}[{\cite[Thm.~6.2.1]{Alesker2011b}}]\label{thm:fourier_interchange}
		Let $i: L \hookrightarrow V$ be an injection of linear spaces and let $\phi \in \Val^\infty(V)$. Then $i^\ast \phi \in \Val^\infty(L)$ and 
		\begin{align*}
			\mathbb{F}_L (i^\ast \phi) = (i^\vee)_\ast \mathbb{F}_V(\phi),
		\end{align*}
		where $i^\vee$ denotes the adjoint map.
	\end{Theorem}

	The space $\Val_j^{-\infty}(V)$ contains the subspaces of Klain--Schneider (KS)-continuous valuations $\Val_j^{\KS}(V)$, defined in \cite{Bernig2017} as the completion of $\Val_j^{\infty}(V)$ in the norm $\|\cdot\|_{\Sc}$, which depends on a choice of Euclidean structure and is given by $\|\psi\|_{\Sc} = \|\Sc_ \psi \|_{\infty}, \psi \in \Val_j^{\infty}(V)$. Note that a generalized valuation $\phi$ is KS-continuous precisely if its (generalized) Schneider section is continuous. Moreover, the operation of pullback by linear injections extends by continuity to $\Val_j^{\KS}(V)$, making it a natural setting for a generalization of Theorem~\ref{mthm:resProperty}.

\section{Extension from the full grassmannian}\label{sec:extFullGrass}
	In this section we will give a proof of Theorem~\ref{mthm:resProperty}. To this end, we reformulate Theorem~\ref{mthm:resProperty} using the Alesker--Fourier transform, yielding a dual formulation where restrictions are replaced by pushforwards. As every smooth valuation of degree $1 \leq k \leq n-1$ can be represented by a section over $\mathbb{P}_+ (V^\ast)$ with values in the subspace $A_{n-k}$ of symmetric double forms (see Sections~\ref{sec:prepDoubleforms} and \ref{sec:bgValDiffForm}), twisted by a line bundle, and since taking the pushforward of a valuation corresponds to restricting these double forms (see Proposition~\ref{prop:ConvPushfOnQuadrForms}), we need to solve an extension problem for such double forms (using representation theory) and then show that the (unique) extension represents a globally defined valuation.

	\subsection{Extending double forms}\label{sec:extFullGrassDoubleForms}
	We start with the question of extending double forms. For this reason, let $k \geq 0$ and denote by
	\begin{align*}
		\res_H: Y_k(V) \to Y_k(H), \quad H \in \Gr_{n-1}(V),
	\end{align*}
	the restriction map. As a first step, we show in the following lemma that $\res_H$ respects the decomposition $Y_k(V) = Z_k(V) \oplus A_k(V)$.
	
	\begin{Lemma}\label{lem:restrict_Ak}
		Let $H \in \Gr_{n-1}(V)$. Then the restriction map $\res_H$ satisfies
		\begin{align*}
			\res_H(Z_k(V))=Z_k(H) \quad \text{ and } \quad \res_H(A_k(V))=A_k(H). 
		\end{align*}
	\end{Lemma}
	\proof
	It is clear that $\res_H(Z_k(V))\subset Z_k(H)$, as decomposable vectors in $\wedge^k H$ are also decomposable in $\wedge^k V$. 
	
	If $n\leq 3$, or $k\in\{1,n-1\}$, then $Z_k(V)=\{0\}$, $Z_k(H)=\{0\}$ for all $H$, and the claim is trivial as $\res_H$ is clearly surjective.
	
	When $n=4$ or $k=n-2$, we have $Z_k(H)=\{0\}$, and so $\res_H(Z_k(V))=\{0\}=Z_k(H)$. Hence, as $\res_H$ is surjective and linear,
	\begin{align*}
		A_k(H)&=Y_k(H)=\res_H(Y_k(V)) \\&= \res_H(A_k(V)) + \res_H(Z_k(V)) =\res_H(A_k(V)).
	\end{align*}
	
	\medskip
	
	Now assume $n\geq 5$, $2\leq k\leq n-3$, and consider the subspace $X\subset Y_k(V)$ of all quadratic forms $Q$ such that $\res_E(Q)\in A_k(E)$ for all $E\in\Gr_{n-1}(V)$. This is evidently a $\GL(V)$-invariant linear subspace.
	
	Fixing a Euclidean structure $P$ on $V$ and letting $P_E=(P|_E)^{\wedge k}\in Y_k(E)$, $P_V=P^{\wedge k}\in Y_k(V)$ be the induced Euclidean structures, we have $\res_E(P_V)=P_E$. By Lemma \ref{lem:1d_to_Ak}, $P_E \in A_k(E)$ for all $E$, and therefore $P_V\in X$. As $P_V\in A_k(V)$, while $A_k(V)$ is an irreducible $\GL(V)$-module, it follows that $A_k(V)\subset X$.
	
	Hence $\res_H(A_k(V)) \subseteq A_k(H)$. Recalling that also $\res_H(Z_k(V))\subset Z_k(H)$, while $\res_H(Y_k(V))=Y_k(H)$ by Lemma \ref{lem:plucker_decomposition}, we must have equality in both inclusions.
\endproof

Next we prove that under the obviously necessary compatibility conditions, every continuous assignment $H\mapsto Q_H \in A_k(H)$, $H\in\Gr_{n-1}(V)$, can be extended to a global quadratic form.
\begin{Proposition}\label{prop:lift_quadrform_exists}
	Assume $1\leq k\leq n-2$, and suppose that $Q_H\in A_k(H)$ is given for all hyperplanes $H\in\Gr_{n-1}(V)$ and
	\begin{align*}
		\res_{H\cap H'}Q_{H}=\res_{H\cap H'}Q_{H'}, \quad \forall H \neq H' \in \Gr_{n-1}(V).
	\end{align*}
	If, moreover, $H \mapsto Q_H$ is continuous, then there exists a unique quadratic form $Q\in A_k(V)$ such that $\res_H(Q)=Q_H$ for all $H \in \Gr_{n-1}(V)$.
\end{Proposition}

\noindent
In the proof of this statement, we make essential use of the following fact proved in \cite[Prop.~4.4]{Wannerer2018}. Here, we denote by $\wedge^k_sV$ the subset of decomposable vectors in $\wedge^kV$.

\begin{Proposition}[\cite{Wannerer2018}]\label{prop:wannerer}
	Let $f:\wedge^k_sV\to \R$ be a continuous even function, such that for any hyperplane $H\subset V$ there exists a quadratic form $Q_H\in Y_k(H)$ such that $f|_{\wedge^k_s H}=Q_H|_{\wedge^k_s H}$. Then there exists $Q\in Y_k(V)$ such that $f=Q|_{\wedge^k_sV}$. 
\end{Proposition}

\noindent We are now ready to prove Proposition~\ref{prop:lift_quadrform_exists}.
\begin{proof}[Proof of Proposition~\ref{prop:lift_quadrform_exists}]
The compatibility condition allows to define a continuous function $f:\wedge^k_s V\to \R$ on the set of all decomposable vectors by setting $f(X)=Q_H(X)$ for any $H$ with $X\in\wedge^kH$. The continuity of $f$ follows from the continuity of $H \mapsto Q_H$. Also $f$ is even, as $Q_H$ is a quadratic form.

By Proposition \ref{prop:wannerer}, we can find $Q\in Y_k(V)$ such that $Q|_{\wedge^k_sV}=f$. We may moreover assume that $Q\in A_k(V)$: writing $Q=Q_A + Q_Z$ with $Q_A \in A_k(V)$ and $Q_Z \in Z_k(V)$ as in Lemma~\ref{lem:plucker_decomposition}, we may replace $Q$ with $Q_A$ since $Q_Z|_{\wedge^k_sV} = 0$.

\noindent As $\res_H(Q)$ and $Q_H$ coincide on decomposable vectors, for all $H\in\Gr_{n-1}(V)$,
\begin{align*}
	\res_H(Q)-Q_H\in Z_k(H).
\end{align*}
By Lemma \ref{lem:restrict_Ak}, $\res_H(Q)\in A_k(H)$, hence $\res_H(Q)-Q_H\in Z_k(H)\cap A_k(H) = \{0\}$, that is, $\res_H(Q)=Q_H$ for all $H$.

For the uniqueness of $Q$, merely observe that any lift $\widetilde Q\in A_k(V)$ must satisfy $\widetilde Q-Q\in Z_k(V)$, as $\wedge^k_s V\subset \cup_{H\in\Gr_{n-1}(V)}\wedge^k H$. As $A_k(V)\cap Z_k(V)=\{0\}$ by Lemma~\ref{lem:plucker_decomposition}, it follows that $\widetilde Q=Q$.
\end{proof}

\medskip

In the application of Proposition~\ref{prop:lift_quadrform_exists}, we will need to ensure that the lift $Q$ depends smoothly on the data $(Q_H)_{H \in \Gr_{n-1}(V)}$. The rest of the section is devoted to establishing this fact.

Let $M$ be a smooth manifold of dimension $m$, and let $G$ be a Lie group acting transitively on $M$. Consider the Fr\'echet bundle $\mathcal E_k$ over $M$, with fiber
\begin{align*}
	\mathcal E_k|_x=\Gamma(\Gr_{m-1}(T_xM), A_k(H))
\end{align*} 
over $x \in M$. It is naturally a $G$-bundle. The subspace
\begin{align*}
	E_k|_x=\{q\in \mathcal E_k|_x: \res_{H\cap H'}(q_H)=\res_{H\cap H'}(q_{H'}), \,\, \forall H,H'\in\Gr_{m-1}(T_xM)\}
\end{align*}
of compatible quadratic forms is finite-dimensional by Proposition \ref{prop:lift_quadrform_exists}, and is evidently $\Stab(x)$-invariant. It therefore defines a $G$-subbundle $E_k$ of finite rank of $\mathcal E_k$, which by Proposition \ref{prop:lift_quadrform_exists} is isomorphic to the bundle $\mathcal A_k$ over $M$ with fiber $\mathcal A_k|_x=A_k(T_xM)$ over $x\in M$. Let $\res:\mathcal A_k\to \mathcal E_k$ denote the natural inclusion, given by restrictions to the various hyperplanes. Applying those remarks with $M=S^{n-1}$ and $G=\SO(n)$, we find

\begin{Proposition}\label{prop:smooth_quad_extension}
Let $1 \leq k \leq n-3$. A smooth section $s\in\Gamma(S^{n-1}, \mathcal E_k)$ lying in $E_k$ is given by $\res(q)$ for $q\in \Gamma(S^{n-1}, A_k(T_xS^{n-1}))$. 
\end{Proposition}

\subsection{Smoothness of compatible sections}
	
	Natural extensions of the embeddings by Klain and Schneider yield useful characterizations of smooth compatible sections in $V_j^{(\pm), \infty}(r,\R^n)$.
	
	Indeed, consider first the case of even valuations. Let $\Gr_r(\R^n) \ni E \mapsto \phi(E) \in \Val_j^{+,\infty}(E)$ be a compatible assignment (not necessarily smooth), that is $\phi_E|_{E\cap E'}=\phi_{E'}|_{E\cap E'}$, $\forall E, E' \in \Gr_r(\R^n)$. Denote, abusing notation, by $\Kl_\phi$ the assignment of an element in $\Dens(E)$ to every $E \in \Gr_j(\R^n)$, defined by $\Kl_ \phi(E) = \phi(F)|_E$, where $F \in \Gr_r(\R^n)$ with $E \subset F$ is arbitrary. Note that $\Kl_ \phi$ is well-defined by the compatibility condition and Theorem~\ref{thm:KlSchnInjective}.
	
	In the following, we will give a criterion on when such a compatible assignment $\phi$ is smooth, that is, when $\phi$ is an element of $U_j^{+, \infty}(r, \R^n)$. 
	
	\begin{Lemma}\label{lem:charSmoothCompSecEven}
		A compatible assignment $\phi$ as above satisfies
		\begin{align*}
			\phi \in U_j^{+, \infty}(r, \R^n) \quad \iff \quad  \Kl_ \phi \in \Gamma(\Gr_j(\R^n), \Dens(E)).
		\end{align*}
	\end{Lemma}
	\begin{proof} Note that the statement is trivial for $r < j$ as then both $\phi$ and $\Kl_\phi$ are zero, while for $r=j$  it is tautological. Therefore, let $r \geq j+1$.
		
		Assume first that $\phi$ is a smooth section and let $U'$ be a neighborhood of $F_0 \in \Gr_j(\R^n)$. Let $g: U' \to \OO(n)$ be a smooth map with $g(F)F_0 = F$ and pick $E_0 \in \Gr_r(\R^n)$ with $E_0 \supset F_0$. Next, fix a neighborhood $U$ of $E_0 \in \Gr_r(\R^n)$ such that $g(F)E_0 \in U$ for all $F \in U'$. By possibly making $U'$ and $U$ smaller, we may assume that $\pi^{-1}U$ is trivializable and consider $\phi$ as a map $U \to \Val_j^{+,\infty}(E_0)$. The map
		\begin{align*}
			F \mapsto g(F)^{-1} \cdot \phi(g(F)E_0)|_F = \left.\left(g(F)^{-1} \cdot \phi(g(F)E_0)\right)\right|_{F_0}
		\end{align*}
		is smooth as composition of a smooth map (since $\phi$ is smooth) and the linear and continuous restriction map.  Hence, $\Kl_ \phi$ is smooth.
		
		Next, assume that $\Kl_ \phi$ is smooth and pick a Euclidean structure, trivializing $\Gamma(\Gr_j(\R^n), \Dens(E)) \cong C^\infty(\Gr_j(\R^n))$. Fix a trivializable neighborhood $U$ of $E_0 \in \Gr_r(\R^n)$, as before, and consider $\phi$ as a map $U \to \Val_j^{+,\infty}(E_0)$. We need to show that this map is smooth. To this end, let $\gamma = E_t \in U$ be a smooth curve and consider $\phi \circ \gamma: \R \to \Val_j^{+,\infty}(E_0)$. As the Klain embedding $\Kl:\Val_j^{+,\infty}(E_0)\to C^\infty(\Gr_j(E_0))$ is an isomorphism of Fr\'echet spaces (onto its closed image), $\phi \circ \gamma$ is smooth exactly if $\Kl \circ \phi \circ \gamma: \R \to C^\infty(\Gr_j(E_0))$ is smooth. This is clear, however, as $\Kl \circ \phi \circ \gamma$ is the smooth map $t \mapsto \Kl _\phi|_{\Gr_j(E_t)}$ composed with the trivialization map. 
	\end{proof}
	
	Next, consider the case of odd valuations. For a compatible assignment \[\Gr_r(\R^n) \ni E \mapsto \phi(E) \in \Val_j^{-,\infty}(E),\] we denote by $\Sc_ \phi$, abusing notation, the assignment $(E,F) \mapsto \Sc_{\phi(G)}(E,F)$, where $G \in \Gr_r(\R^n)$ is arbitrary with $(E,F) \in \widehat{\mathcal{F}}_j^{j+1}$, $F \subset G$, and $\Sc_{\phi(G)}$ is the Schneider section of $\phi(G) \in \Val_j^{-,\infty}(G)$.
	
	\begin{Lemma}\label{lem:charSmoothCompSecOdd}
		A compatible assignment $\phi$ as above satisfies
		\begin{align*}
			\phi \in U_j^{-, \infty}(r, \R^n) \quad \iff \quad  \Sc_ \phi \in \Gamma^-(\widehat{\mathcal{F}}_{j}^{j+1}, \Dens(F))/L.
		\end{align*}
	\end{Lemma}
	\begin{proof}
		The claim follows by a similar argument as in the proof of Lemma~\ref{lem:charSmoothCompSecEven} with the Schneider embedding replacing the Klain embedding.
	\end{proof}
	
	We note a direct consequence of Lemmas~\ref{lem:charSmoothCompSecEven} and \ref{lem:charSmoothCompSecOdd} for later reference.
	
	\begin{Lemma}\label{lem:smooth_restriction}
		Let $r \geq j+2$. A compatible assignment
		\begin{align*}
			\Gr_r(\R^n) \ni E \mapsto \phi(E) \in \Val_j^\infty(E)
		\end{align*}
		is in $U_j^\infty(r, \R^n)$ if and only if the assignment 
		\begin{align*}
			\Gr_{j+2}(\R^n) \ni F \mapsto \phi(E)|_F \in \Val_j^\infty(F), \quad \text{ for arbitrary $F \subset E \in \Gr_r(\R^n)$,}
		\end{align*}
		is in $U_j^\infty(j+2, \R^n)$.
	\end{Lemma}

\subsection{Proof of Theorem~\ref{mthm:resProperty}}
We now have all ingredients needed for the proof of Theorem~\ref{mthm:resProperty}. However, before proceeding to the actual proof, we do two further preparatory steps.

\medskip

First, we observe that it suffices to prove Theorem~\ref{mthm:resProperty} for $r=n-1$. Indeed, assume this case is established for all $n$. 
For an arbitrary $j+2\leq r\leq n-2$, take $\psi\in U^\infty_j(r, V)$, and fix $F\in \Gr_{r+1}(V)$. Consider the natural restriction $\psi|_F\in U^\infty_j(r, F)$. By assumption, one may find $\phi_F\in\Val_j^\infty(F)$ with $\res_r(\phi_F)=\psi|_F$, and by Theorem~\ref{thm:KlSchnInjective} such $\phi_F$ is uniquely defined. The assignment $F\mapsto \phi_F$ is smooth by Lemma \ref{lem:smooth_restriction}. We now define $\psi^{r+1}\in U_j^\infty(r+1, V)$ by $\psi^{r+1}(F)=\phi_F$.
Repeating the argument for $r+1, r+2, \dots, n-1$ we finally arrive at a valuation $\phi\in \Val_j^\infty(V)$ with the property that $\phi|_E=\psi(E)$ for all $E\in\Gr_r(V)$, that is $\res_r(\phi)=\psi$.

\medskip

The second (and main) preparatory step is to replace restrictions with push-forwards. Writing $\Gr^q_r(V)$ for the grassmannian of $r$-dimensional quotient spaces of $V$, we put
\begin{align*}
	\widetilde V_k^\infty(r, V)=\Gamma(\Gr^q_{r}(V), \Val_{k}^\infty(F)\otimes \Dens(F)^*)
\end{align*}
as the Alesker--Fourier dual analogue of $V_k^\infty(r,V)$. By Theorem \ref{thm:fourier_interchange}, the Alesker--Fourier transform interchanges pullback (restriction) and pushforward. Consequently, we define $\widetilde U_k^\infty(r, V)$ as the subspace of all sections $\psi\in \widetilde V_k^\infty(r, V)$ satisfying
\begin{align*}
(\pi^{E_1}_{E_1 + E_2})_*\psi(F_1)=(\pi^{E_2}_{E_1 + E_2})_*\psi(F_2), \quad \forall F_i = V/E_i\in \Gr^q_{r}(V), i=1,2,
\end{align*}
where for $E \subset E'$, $\pi^E_{E'} : V/E \to V/E'$ denotes the canonical quotient map.
	
Defining the (simultaneous) pushforward $\mathrm{push}_r(\phi)\in \widetilde V_k(r, V)$ of a twisted valuation $\phi\in\Val_{k+n-r}^\infty(V)\otimes \Dens(V)^*$ by
\begin{align*}
	\mathrm{push}_r(\phi)(F)= (\pi_F)_*\phi, \quad F \in \Gr^q_r(V),
\end{align*}
where $\pi_F:V\to F$ is the quotient map, and applying the Alesker--Fourier transform, Theorem \ref{mthm:resProperty} is easily seen to be equivalent to the following
	
\begin{Theorem}
	Assume $n-1 \geq r > k \geq 2$. The image of 
	\begin{align*}
		\mathrm{push}_r:\Val_{n+k-r}^\infty(V)\otimes\Dens(V)^*\to \widetilde V_k(r, V)
	\end{align*}
	coincides with $\widetilde U^{\infty}_k(r, V)$.
\end{Theorem}
\proof
By the first preparatory step described above, we may assume that $r=n-1$.

Next, take $ \psi\in \widetilde U^{\infty}_k(n-1, V)$. For $F\in\Gr_{n-1}^q(V)$ and $\xi\in\mathbb P_+(F^*)$, write 
\begin{align*}
	Q_{F, \xi}=Q(\psi(F))|_\xi\in  A_{n-k-1}(F^*/\xi)\otimes \xi^{\otimes(n-k-1)} .
\end{align*}
By the compatibility assumption, the collection of forms $Q_{F, \xi}$ defines a smooth section of $\mathcal E_{n-k-1}$ over $\mathbb P_+(V^*)$ which in fact belongs to $E_{n-k-1}$, by Proposition~\ref{prop:ConvPushfOnQuadrForms}. By Proposition~\ref{prop:smooth_quad_extension}, we may find a smooth section of quadratic forms $Q\in\Gamma (\mathbb P_+(V^*),  A_{n-k-1}(V^*/\xi)\otimes  \xi^{\otimes(n-k-1)})$ with $\res_{F^*/\xi}(Q|_\xi)=Q_{F, \xi}$.

Let us fix a Euclidean structure and orientation on $V$. Let us write $\tau_Q\in\Omega^{n-k-1}(S(V^*), \wedge^{n-k-1}V)$ for the form corresponding to $Q$, and we denote by $\widetilde \tau_Q \in\Omega^{n-k-1}(S(V^*), \wedge^{k+1}V^*)$ the corresponding form under the natural identification $\wedge^{n-k-1}V\cong\wedge^{k+1}V^*$. By construction, $\widetilde{\tau}_Q$ is a vertical form. It remains to see that, if necessary, we can modify $\tau_Q$ by a vertical form so as to make it closed, without altering its restrictions to hyperplanes.

For a quotient space $H\in \Gr_{n-1}^q(V)$, write
\begin{align*}
	\res_{H^\ast}\tau_Q\in\Omega^{n-k-1}(S(H^\ast), \wedge^{n-k-1}H)
\end{align*}
for the restriction of $\tau_Q$, obtained by restricting the form $\tau_Q$ to $S(H^\ast)$, and subsequently applying to its values the canonical extension of the quotient map $V \to H$ to $\wedge^{n-k-1}V \to \wedge^{n-k-1}H$. 

By the construction of $Q$ and Proposition~\ref{prop:ConvPushfOnQuadrForms}, $ \res_{H^\ast} \tau_Q$ corresponds to the defining form of $\psi(H)$ under the identification $\wedge^{n-k-1} H \cong \wedge^k H^\ast\otimes \ori(H)$. Consequently, $d(\res_{H^\ast}\tau_Q)=0$ for all $H\in\Gr_{n-1}^q(V)$, that is,
\begin{equation}
	\label{eq:hyperplane_pairing_zero}
	\langle d\tau_Q|_\xi(u_1,\dots, u_{n-k}), v \rangle =0
\end{equation}
whenever there exists $H \in \Gr_{n-1}^q(V)$ such that $\xi \in S(H^\ast)$, $u_1, \dots, u_{n-k} \in T_\xi S(H^\ast)$ and $v \in \wedge^{n-k-1}H^\ast$. We claim that this implies that
\begin{align}\label{eq:prfExtFullDTau}
	d\widetilde\tau_Q=\beta\wedge\widetilde \psi_k+\gamma\wedge\psi_{k+1},
\end{align} 
where $\widetilde \psi_k\in\Omega^k(S(V^*), \wedge^{k+1}V^*)$ is defined by
\begin{align*}
	\widetilde \psi_k|_\xi(u_1, \dots, u_{k})=\xi\wedge u_1\wedge\dots\wedge u_k,
\end{align*}
and $\psi_{k+1}\in\Omega^{k+1}(S(V^*), \wedge^{k+1}V^*)$ by
\begin{align*}
	\psi_{k+1}|_\xi(u_1, \dots, u_{k+1})= u_1\wedge\dots\wedge u_{k+1},
\end{align*} 
for $u_1, \dots, u_{k+1} \in T_\xi S(V^\ast)$ and $\xi \in S(V^\ast)$, while $\beta\in\Omega^{n-2k}(S(V^*))$ and $\gamma\in \Omega^{n-2k-1}(S(V^*))$ are some scalar-valued forms.

Indeed, fix $\xi \in S(V^\ast)$ and assume that $u_1,\dots, u_{n-k}\in T_\xi S(V^*)$ are linearly independent. Choose $w_1,\dots, w_{k-1}\in T_\xi  S(V^*)$ that together with $\xi$ complement them to a basis $B$ of $V^*$, and decompose $d\widetilde \tau_Q|_\xi(u_1,\dots, u_{n-k})$ in the corresponding basis of $\wedge^{k+1}V^*$. If $u_I\wedge w_J$ or $u_I\wedge w_J\wedge \xi$ appears in the sum with non-zero coefficient and $|J|\geq 1$, then choose $v_1,\dots  v_{n-k-1}$ as the complement in $B$ of $\{u_i\}_{i\in I}\cup \{w_j\}_{j\in J}$, resp. $\{u_i\}_{i\in I}\cup \{w_j\}_{j\in J}\cup\{\xi\}$. Then
\begin{align*}
	0 &\neq d\widetilde{\tau}_Q|_{\xi}(u_1, \dots, u_{n-k}) \wedge v_1 \wedge \dots \wedge v_{n-k-1} \\&= \langle d \tau_Q|_{\xi}(u_1, \dots, u_{n-k}), v_1 \wedge \dots \wedge v_{n-k-1} \rangle,
\end{align*}
which contradicts \eqref{eq:hyperplane_pairing_zero} since $u_1,\dots, u_{n-k}, v_1, \dots, v_{n-k-1}, \xi$ lie in $\Span\{ B\setminus \{w_j\}_{j\in J} \}$. We conclude that 
$$d\widetilde \tau_Q|_\xi(u_1,\dots, u_{n-k})=\sum_{|I|=k}a_I \xi\wedge\bigwedge_{i\in I} u_i + \sum_{|J|=k+1}b_J \bigwedge_{j\in J} u_j$$ 
for some coefficients $a_I, b_J$. The claim now follows since $\widetilde \tau_Q|_\xi$ is multi-linear and anti-symmetric in the vectors $u_i$.

Next, observe that, if $n-2k<0$, $\beta$ and $\gamma$ vanish and, by \eqref{eq:prfExtFullDTau}, $d\widetilde \tau_Q|_\xi=0$. For the general case note that $d\widetilde \psi_k=\psi_{k+1}$ and, consequently, the equation
$$0=d(d\widetilde\tau_Q)=d\beta\wedge\widetilde \psi_k+(-1)^{n-2k}\beta\wedge \psi_{k+1}+d\gamma\wedge\psi_{k+1},$$
implies that $\beta=(-1)^{n+1}d\gamma$ and $d\beta=0$. Thus
$$d\widetilde\tau_Q=(-1)^{n+1}d\gamma\wedge\widetilde \psi_{k}+\gamma\wedge d\widetilde \psi_{k}=(-1)^{n+1} d(\gamma\wedge \widetilde \psi_{k}).$$
and, therefore, the form 
$$\widetilde\zeta:=\widetilde \tau_Q+(-1)^{n}\gamma\wedge \widetilde \psi_{k}\in\Omega^{n-k-1}(S(V^*), \wedge^{k+1}V^*)$$
is closed.

 As $ \gamma\wedge\widetilde{\psi}_k$ is vertical, $\widetilde \zeta$ defines a vertical and closed form in $\Omega^{k+1, n-k-1}(V\times S(V^*))^{tr}$, which, by Theorem~\ref{thm:bernigBroecker}, corresponds to a valuation $\phi\in\Val^\infty_{k+1}(V)$.
Moreover, notice that
\begin{align*}
	(\gamma\wedge \widetilde \psi_{k})(u_1\wedge\dots\wedge u_{n-k-1}) \wedge v_1\wedge\dots\wedge v_{n-k-1}=0
\end{align*}
whenever $u_1,\dots, u_{n-k-1}, v_1, \dots, v_{n-k-1}, \xi \in H^\ast$ for some $H \in \Gr_{n-1}^q(V)$, since $\dim H=n-1$, while the wedge product is a linear combination of wedges of $n$ vectors all belonging to $H^\ast$. Consequently, 
\begin{align*}
	\res_{H^\ast/\Span\xi}(Q(\phi)|_{\xi})=\res_{H^\ast/\Span\xi}(Q|_{\xi}), \quad \xi \in S(H^\ast),
\end{align*}
for any hyperplane $H\in \Gr_{n-1}^q(V)$, that is, by Proposition~\ref{prop:ConvPushfOnQuadrForms}, we obtain $\push_{n-1}(\phi)=\psi$, concluding the proof.
\endproof

\begin{Remark}
	Note that in the proof above, the extension $Q$ is actually unique by Proposition~\ref{prop:lift_quadrform_exists}. Consequently, the form $\gamma$ must be zero, that is, our proof shows that already $\tau_Q$ is closed. 
\end{Remark}

A standard approximation argument extends Theorem \ref{mthm:resProperty} to the KS-continuous valuations $\Val^{\mathrm{KS}}_k(V)$ (see Section~\ref{sec:bgopVal}).
Let $V^{\KS}_j(r, \R^n)$ denote the continuous global sections of the Banach bundle $\Val_j^{\KS}(\Gr_r(\R^n))$, whose fiber over $E\in\Gr_r(\R^n)$ is the space of KS-continuous valuations $\Val^{\KS}_j(E)$. We equip $V^{\KS}_j(r, \R^n)$ with the Banach norm 
$$\| \phi\|_{\Sc}:=\sup_{E\in\Gr_r(\R^n)}\|\phi(E)\|_{\Sc}.$$
The subspace  $U^{\KS}_j(r, \R^n)\subset V^{\KS}_j(r, \R^n)$ of compatible sections is defined as in the smooth case.

\noindent
The direct analogue of Theorem~\ref{mthm:resProperty} for KS-continuous valuation now reads
\begin{Theorem}\label{thm:resPropertyCont}
	If $r\geq j+2$, $\mathrm{Image}(\res_j:\Val^{\KS}_j(\R^n)\to V^{\KS}_j(r, \R^n))=U^{\KS}_j(r, \R^n)$. 
\end{Theorem}

Let us finally deduce Corollary~\ref{mcor:imCosTransf} from Theorem~\ref{mthm:resProperty} for even valuations.
\begin{proof}[Proof of Corollary~\ref{mcor:imCosTransf}]
	By \cite{Alesker2004c}, the image of the $i$-cosine transform coincides with the image of the Klain map on $i$-homogeneous even valuations.
	The corollary is then easily seen to follow from Theorem~\ref{mthm:resProperty}.
\end{proof}

\section{Extension from finite arrangements of subspaces}\label{sec:finArr}
In this section we give a proof of Theorem~\ref{mthm:extensionFinArrang} using the representation of valuations by differential forms (see Section~\ref{sec:bgValDiffForm}). Working again in the dual setting (that is, after applying the Alesker--Fourier transform), we are thus given differential forms on finitely many subspaces, which we aim to extend to a globally defined differential form. In contrast to Section~\ref{sec:extFullGrass}, we will work with the $(n-1)$-form associated to a valuation and not the $n$-form from Theorem~\ref{thm:bernigBroecker}. This approach has the advantage that, once we find a form that restricts to the given data, it defines a valuation without any further assumptions such as being closed and vertical. However, as the representation by an $(n-1)$-form is not unique, the compatibility assumptions for the given valuations do not imply compatibility of differential forms. We therefore need an additional "alignment" step in which we replace the given differential forms by compatible (aligned) forms representing the same valuations (see Proposition~\ref{prop:compatible_primitives}).

In the end of the section, we will discuss the necessity of the conditions of Theorem~\ref{mthm:extensionFinArrang} as well as some related results.


\subsection{Preparations and extensions of linear forms}\label{sec:finArrLinForms} By an arrangement we mean a finite collection of linear subspaces. The notion of \emph{minimally intersecting arrangement} introduced in Definition \ref{def:MinIntersSemiGeneric} will play a key role throughout the section. We first collect some facts about such arrangements, and then proceed to study extension problems for linear forms on various arrangements. 

The following basic fact from linear algebra will often be used implicitly. We use the convention $\wedge^k E=0$ when $k\notin [0,\dim E]\cap \Z$. We will often deal with collections $\{E_i\}$ of linear subspaces, and will henceforth use the notation $E_I=\cap_{i\in I}E_i$ without comment. We typically omit the curly brackets in the index, and list the elements of $I$ in increasing order. If a Euclidean structure is given in the ambient space, we similarly write $S_I=S(E_I)$ for the unit sphere.

\begin{Lemma}
	Let $X_1,\dots, X_m\subset \R^n$ be subspaces. Then $\wedge^k(X_1\cap \dots\cap X_m)=\wedge^k X_1\cap \dots\cap \wedge^k X_m$.
\end{Lemma}
\begin{proof}
	Induction immediately reduces the statement to $m=2$. Denote $X=X_1$, $Y=X_2$. Choose a basis $e=\{e_1, \dots, e_a\}$ of $X\cap Y$, and let $f=\{f_1,\dots, f_b\}\subset X$ and $g=\{g_1, \dots, g_c\}\subset Y$ complete $e$ to a basis of $X$, respectively $Y$. Then the $k$-wedges of vectors from $e$, $e \cup f$ and $e \cup g$ form bases of $\wedge^k (X \cap Y)$, $\wedge^k X$ and $\wedge^k Y$, respectively. As the $k$-wedges of vectors from $e\cup f\cup g$ are all linearly independent, the claim follows.
\end{proof}
	
\noindent 
Next, we spell out a few basic examples and easily checked properties of minimally intersecting and semi-generic arrangements.
\begin{itemize}
	\item A pair of subspaces $E_1, E_2$ is always minimally intersecting within $E_1+E_2$.
	\item Any finite arrangement of subspaces in general position is semi-generic.
	\item A subset of a minimally intersecting arrangement is minimally intersecting.
	\item A subset of a semi-generic arrangement is semi-generic.
	\item A minimally intersecting arrangement of proper subspaces in $\R^n$ has size at most $n$.
\end{itemize}

We will further need the following properties, for which we give a proof for the reader's convenience. Recall for the following that $E_I$ denotes $\cap_{i\in I}E_i$.
\begin{Lemma}\label{lem:min_intersecting}
	Let $\{E_i\}_{i=1}^N$ in $\R^n$ be a minimally intersecting arrangement. Then 
	\begin{enumerate}
		\item $E_i+E_j=\R^n$ for all $i\neq j$.
		\item If $F\subset\cap_{i=1}^NE_i$, then $\{E_i/F\}_{i=1}^N$ are minimally intersecting in $\R^n/F$.
		\item For any $\emptyset\neq I\subset\{1,\dots, N\}$, $\{E_j\cap E_I\}_{j\notin I}$ is minimally intersecting in $E_I$.
	\end{enumerate}
\end{Lemma}
\proof
To show (i), note that, by assumption, $\dim (E_i\cap E_j)=\dim E_i+\dim E_j -n$, and since $\dim (E_i\cap E_j)=\dim E_i+\dim E_j-\dim (E_i+E_j)$, we conclude that $\dim (E_i+E_j)=n$, that is $E_i+E_j=\R^n$.

The second claim follows directly, as the codimensions of $E_i\subset \R^n$ and $E_i/F\subset\R^n/F$ are equal, while $(\cap_{i=1}^N E_i)/F=\cap_{i=1}^N (E_i/F)$.

For the final statement, it suffices by induction to check that $\{E_N\cap E_i\}_{i=1}^{N-1}$ is minimally intersecting in $E_N$, which is straightforward. 
\endproof

Minimally intersecting arrangements are best described by well chosen bases.

\begin{Lemma}\label{lem:BaseGenericFam}
	Suppose that $E_1, \dots, E_N$ is a minimally intersecting arrangement of subspaces in $\R^n$, and $N\geq 2$. Then there exist unordered bases $\underline e_i$ of $E_i$ such that $\underline e = \bigcup_{i=1}^N \underline e_i$ is a basis of $\R^n$.
\end{Lemma}
\proof
As the arrangement $E_1,\dots, E_N$ is minimally intersecting, the sum $(\cap_{i=1}^N E_i)^\perp = E_1^\perp + \dots + E_N^\perp \subset (\R^n)^\ast$ is direct. Choose an inner product on $(\R^n)^*$ for which $E_i^\perp$, $i=1, \dots, N$, are pairwise orthogonal, and use it to  identify $\R^n \cong (\R^n)^\ast$. We then choose orthonormal bases $\underline c_i$ of $E_i^\perp$, $i=1, \dots, N$. Picking an orthonormal basis $\underline{b}$ of $\cap_{i=1}^N E_i$, we obtain an orthonormal basis $\underline{b} \cup \bigcup_{i=1}^N \underline c_i$ of $\R^n$, so that $E_i = \Span\{\underline{b} \cup \bigcup_{j \neq i} \underline c_j\}$. 
\endproof

The following property of collections of subspaces will play a major role.
\begin{Definition}
	 Let $S$ be a collection of subspaces in $V$, and $k\geq 0$. We will say that $S$ has the \emph{$k$-extension property} if, whenever $\omega_E\in \wedge ^k E^*$ is given for all $E\in S$ such that $\omega_E|_{\wedge^k (E\cap E')}=\omega_{E'}|_{\wedge^k(E\cap E')}$ for all $E\neq E'$ in $S$, then one can find a form $\omega\in\wedge^kV^*$ such that $\omega_E=\omega|_{\wedge^k E}$, for all $E\in S$.
\end{Definition}

Next we show that minimally intersecting arrangements have the $k$-ex\-tension property for all $k$.

\begin{Corollary}\label{cor:lift_form_min_intersection}
	Let $S=\{E_1,\dots, E_N\}$ be an arrangement in $V$ which is minimally intersecting within its span. Then $S$ has the $k$-extension property for all $k$. 
\end{Corollary}
\proof
Assume that forms $\omega_i\in \wedge^kE_i^*$ are given such that $\omega_i|_{\wedge^k(E_{i}\cap E_j)}=\omega_j|_{\wedge^k(E_i\cap E_j)}$. We look for $\omega\in\wedge^kV^*$ such that $\omega_i=\omega|_{\wedge^k E_i}$, for all $i$. By the surjectivity of $\wedge^k V^*\to \wedge^k(E_1+\dots+E_N)^*$, we may assume that $E_1+\dots+E_N=V$.

By Lemma~\ref{lem:BaseGenericFam} there exist bases $\underline e_i=\{e^i_1,\dots, e^i_{n_i}\}$ of $E_i$ such that $\underline e:=\cup_{i=1}^N\underline e_i$ is a basis of $V$. The statement then readily follows: the value of $\omega$ on any $k$-wedge of vectors from $\underline e$ is determined by $\omega_i$ if all vectors belong to $\underline e_i$. If such $i$ is not unique, compatibility ensures the choice of $i$ does not matter. If no such $i$ exists, the value can be assigned arbitrarily.
\endproof

More generally, minimally intersecting subspaces allow simultaneous extension in the sense of the following corollary. Here, the extension problem is reformulated in terms of an exact sequence.

\begin{Corollary}\label{cor:lift_double_forms}
	Let $\{E_i\}_{i=1}^N$ be minimally intersecting subspaces in $\R^m$, and $\{F_j\}_{j=1}^N$ be minimally intersecting subspaces in $\R^n$. Then the sequence
	$$ \wedge^p(\R^m)^*\otimes\wedge^q(\R^n)^*\xrightarrow{r} \bigoplus_{i=1}^N \wedge^p E_i^*\otimes\wedge^q F_i^*\xrightarrow{d} \bigoplus_{i<j} \wedge^p(E_i\cap E_j)^*\otimes \wedge^q(F_i\cap F_j)^*$$
	is exact, where $r(\omega \otimes \eta) = \sum_{i=1}^N \omega|_{\wedge^p E_i} \otimes \eta|_{\wedge^q F_i}$ and 

\begin{align*}
	d\left(\sum_{i=1}^N \omega_i \otimes \eta_i\right) = \sum_{i<j} \left(\omega_i|_{\wedge^p(E_i \cap E_j)}\otimes \eta_i|_{\wedge^q(F_i \cap F_j)}  \qquad\right. \\ \left.- \omega_j|_{\wedge^p(E_i \cap E_j)}\otimes \eta_j|_{\wedge^q(F_i \cap F_j)}\right).
\end{align*}

\end{Corollary}
\proof
Choose bases $\underline e_i$ of $E_i$ and $\underline f_i$ of $F_i$ as in Lemma~\ref{lem:BaseGenericFam}, so that $\underline e=\bigcup\underline e_i$ is a basis of $\R^m$, and $\underline f=\bigcup\underline f_i$ is a basis of $\R^n$.
Take $\omega_i\in \wedge^p E_i^*\otimes\wedge^q F_i^*$ in the kernel of $d$, that is, 
$$\omega_i|_{\wedge^p(E_i\cap E_j)\otimes \wedge^q(F_i\cap F_j)}=\omega_j|_{\wedge^p(E_i\cap E_j)\otimes \wedge^q(F_i\cap F_j)}.$$

Let $e_I$ be a $p$-wedge from $\underline e$, and $f_J$ a $q$-wedge from $\underline f$.
Define $\omega(e_I, f_J)$ by $\omega_i(e_I, f_J)$ if for some $i$ it holds that $e_I$ is the wedge of elements from $\underline  e_i$, while $f_J$ is the wedge of elements from $\underline f_i$. By the above, the value of $\omega(e_I, f_J)$ does not depend on the choice of $i$. In all other cases, set $\omega(e_I, f_J)=0$. Then $\omega\in \wedge^p (\R^m)^*\otimes\wedge^q (\R^n)^*$ is the desired lift, that is, $r(\omega) = \sum_{i=1}^N \omega_i$.
\endproof

Corollary~\ref{cor:lift_double_forms} directly implies the following statement:
\begin{Corollary}\label{cor:smoothDependenceLift}
	 Let $\{E_i\}_{i=1}^N$ and $\{F_i\}_{i=1}^N$ be as in Corollary~\ref{cor:lift_double_forms} and suppose that $\bigcap_{i=1}^N E_i \neq \{0\}$. Then there exist linear maps $L_\xi$, $\xi \in S(\bigcap_{i=1}^N E_i)$, mapping compatible forms $\omega_i \in \wedge^p (E_i/\Span\xi)^\ast \otimes \wedge^q F_i^\ast$, $i=1, \dots, N$, to
	 \begin{align*}
	 	\omega = L_\xi\left(\sum_{i=1}^{N}\omega_i \right) \in \wedge^p (\R^m/\Span\xi)^\ast \otimes \wedge^q (\R^n)^\ast,
	 \end{align*}
	 lifting all $\omega_i$, that is, $r(\omega) = \sum_i \omega_i$, and the maps $L_\xi$ depend smoothly on $\xi$.
\end{Corollary}
\proof
Indeed, let $\xi_0 \in S(\bigcap_{i=1}^N E_i)$ be arbitrary. By Lemma~\ref{lem:BaseGenericFam}, there exist bases for the minimally intersecting arrangement $(E_i)_{i=1}^N$, which contain $\xi_0$. Removing $\xi_0$ from the bases and projecting the bases to the factor space then yields bases for the arrangements $(E_i/\Span\xi)_{i=1}^N$ for all $\xi$ belonging to a neighborhood of $\xi_0$. These bases depend smoothly on $\xi$. Consequently, by the construction of the lifts in Corollary~\ref{cor:lift_double_forms}, the lifts are smooth in this neighborhood of $\xi_0$. By a partition of unity argument, finally, the lifts can be chosen to be smooth everywhere.
\endproof

We now consider another class with the extension property.

\begin{Lemma}\label{lem:multiple_hyperplanes}
	Let $S$ be a collection of subspaces in $V$ such that any $S'\subset S$ with $|S'|\leq k+2$ is minimally intersecting. Then $S$ has the $k$-extension property. 
\end{Lemma}
Observe that if $|S|\leq k+2$, the assumption simply reads that $S$ is minimally intersecting, while for $|S|>k+2$ we need only assume that all $S'\subset S$ with $|S'|=k+2$ are minimally intersecting.
\proof
Let $\omega_E\in\wedge ^k E^*$, $E\in S$ be a compatible family of forms. 
When $|S|\leq k+1$, $S$ forms a minimally intersecting arrangement, and we conclude by Corollary \ref{cor:lift_form_min_intersection}.

Next, observe that when $|S|=k+1$, a form $\omega$ such that $\omega|_{\wedge ^k E}=\omega_E$ is uniquely determined by all $\omega_E$, $E\in S$. Indeed, use Lemma \ref{lem:BaseGenericFam} to choose a basis $e_1, \dots, e_n$ of $V$ such that all $E\in S$ are coordinate subspaces. By assumption $E+E'=V$ for all $E\neq E'$ in $S$, and it follows that for each $e_i$, there is at most one $E\in S$ such that $e_i\notin E$. 

Note that the values $\omega(e_{i_1}\wedge \dots\wedge e_{i_k})$, where $i_1<\dots<i_k$, determine $\omega$ uniquely, and since $|S|=k+1$, we can find $E\in S$ such that $e_{i_1},\dots, e_{i_k}\in E$, wherefore $$\omega(e_{i_1}\wedge\dots \wedge e_{i_k})=\omega_E(e_{i_1}\wedge\dots\wedge e_{i_k}).$$

Now assume $|S|\geq k+2$. Fix $S_0\subset S$ with $|S_0|=k+1$, and let $\omega\in\wedge ^k V^*$ be the uniquely defined common lift of $(E, \omega_E)$, $E\in S_0$, which exists by the previous step. Then for each $F\in S\setminus S_0$, the collection $\{E\cap F\}$, $E\in S_0$, is, by assumption and invoking Lemma \ref{lem:min_intersecting}, a minimally intersecting  arrangement in $F$. By the previous step, an extension of the forms $\omega|_{\wedge^k(E\cap F)}$, $E\in S_0$, is unique. As both $\omega_F$ and $\omega|_{\wedge^k F}$ are extensions, they must coincide. This completes the proof.  
\endproof
\begin{Example}
Any collection of hyperplanes in general position satisfies the assumption of Lemma \ref{lem:multiple_hyperplanes} as long as $k\leq \dim V-2$.
\end{Example}

We close the section by verifying the vanishing of a certain cohomology.

\begin{Proposition}\label{prop:first_cohomology}
	Fix $p\geq 1$, and suppose that $S=\{E_i\}_{i=1}^N$ is an arrangement of subspaces in $\R^n$.
	Assume that for all $1\leq i\leq N-3$, $S_i:=\{E_j\cap E_i\}_{j>i}$ has the $p$-extension property in $E_i$.
	Then the sequence
\[ \bigoplus_{i=1}^N \wedge^p E_i^*\xrightarrow{d_1}\bigoplus_{i<j} \wedge^p(E_i\cap E_j)^*\xrightarrow{d_2} \bigoplus_{i<j<k} \wedge^p(E_i\cap E_j\cap E_k)^*\]
is exact, where \[d_1((L_i)_{i=1}^N) = (L_j-L_i)_{i<j}, \qquad d_2((R_{ij})_{i<j}) = (R_{ij} - R_{ik} + R_{jk})_{i<j<k},\] omitting the restrictions to $\wedge^p(E_i \cap E_j)$ and $\wedge^p(E_i \cap E_j \cap E_k)$ everywhere.
\end{Proposition}
\proof

Let $R_{ij}\in \wedge^p (E_i\cap E_j)^*$ be such that $R_{ij}-R_{ik}+R_{jk}=0$ on $E_i\cap E_j\cap E_k$. We have to construct $L_i\in \wedge^p E_i^*$ such that $R_{ij}=L_j|_{E_i\cap E_j}-L_i|_{E_i\cap E_j}$.

Define $L_N=0$, and let $L_{N-1}$ be any extension of $-R_{N-1, N}$. 

For $i\leq N-2$, assume that $L_j\in \wedge^p E_j^*$ is constructed for all $j>i$ such that on $E_j\cap E_k$ one has $R_{jk}=L_k-L_j$ for all $i<j<k$. Let us construct $L_i$.

Define $\lambda_j\in \wedge^p (E_i\cap E_j)^*$ by $\lambda_j=L_j|_{\wedge^p (E_i\cap E_j)}-R_{ij}$ for all $j>i$.
Then $\lambda_j|_{\wedge^p(E_i\cap E_j\cap E_k)}=\lambda_k|_{\wedge^p(E_i\cap E_j\cap E_k)}$ for all $k>j>i$, as this is equivalent to having on $E_i\cap E_j\cap E_k$ the identity
$$ L_j-R_{ij}= L_k-R_{ik}\iff R_{ij}-R_{ik}+R_{jk}=0,$$
which is given. 

For $i=N-2$, observe that $E_N\cap E_{N-2},E_{N-1}\cap E_{N-2}$ are minimally intersecting within their span, and apply Corollary \ref{cor:lift_form_min_intersection} to obtain the desired $L_i$ satisfying $L_i|_{\wedge^p(E_i\cap E_j)}=\lambda_j$ for $j>i$. 

For $i\leq N-3$, $\mathcal F_i=\{E_j\cap E_i\}_{j>i}$ has the $p$-extension property by assumption, and thus one can choose $L_i\in \wedge^p E_i^*$ such that $L_i|_{\wedge^p(E_i\cap E_j)}=L_j|_{\wedge^p(E_i\cap E_j)}-R_{ij}$, as required. This completes the inductive construction of $L_i$, $1\leq i\leq N$.
\endproof
\begin{Corollary}\label{cor:first_cohomology}
	The following satisfy the assumptions and thus the conclusion of Proposition \ref{prop:first_cohomology}:
	\begin{itemize}
		\item Any minimally intersecting arrangement $S$.
		\item Any arrangement $S$ such that any $p+3$ or less subspaces from $S$ are minimally intersecting in $\R^n$.
		\item Any set $S$ with $N=|S|=3$.
	\end{itemize}
\end{Corollary}
\proof
The first case follows by Lemma \ref{lem:min_intersecting} and Corollary \ref{cor:lift_form_min_intersection}. 
The second case follows from Lemma~\ref{lem:multiple_hyperplanes}. In the last case there is nothing to check.
\endproof 

\subsection{Proof of Theorem~\ref{mthm:extensionFinArrang}}
In order to extend differential forms (rather than linear forms), we need not only be able to lift linear forms on tangent subspaces over a fixed point of the sphere, but also extend partially defined sections that respect the previously chosen lift. This is the content of the next statements.

\begin{Proposition}\label{prop:smooth_extension}
	Let $W$ be a finite-dimensional affine bundle over $S^{n-1}\subset\R^n$, and suppose that $\{E_i\}_{i=1}^N$ is a semi-generic arrangement of subspaces in $\R^n$. Let $f_i\in \Gamma (S(E_i),W)$ be smooth sections such that $f_i|_{S_{ij}}=f_j|_{S_{ij}}$ for all $1\leq i<j\leq N$. Then there exists $f\in \Gamma(S^{n-1}, W)$ such that $f|_{S_i}=f_i$ for all $1\leq i\leq N$.
\end{Proposition}
\proof

The idea is to use a partition of unity argument to reduce the statement to the case of minimally intersecting arrangements, which we prove separately at the end of the proof. Indeed, to find a local extension near $x \in S^{n-1}$, we only need to consider those functions $f_i$, $i \in I \subset \{1, \dots, N\}$, whose domains $S(E_i)$ contain $x$. As $x \in \cap_{i\in I} E_i$, the arrangement $(E_i)_{i \in I}$ is minimally intersecting in $\sum_{i\in I} E_i$.

Formally, take an open cover of $S^{n-1}$ consisting of small balls $B_x$, with the further property that if the center $x$ of $B_x$ lies in $S_I= S(\cap_{i\in I} E_i)$ for some maximal $I\subset \{1,\dots, N\}$, then $B_x\cap E_j=\emptyset$ for all $j\notin I$. Choose a finite subcover $B_{x_\nu}$ of $S^{n-1}$.

Fix a subordinate partition of unity $\rho_{\nu} \in C_c^\infty(B_{x_\nu})$.  If we can find functions $g_{\nu} \in \Gamma(B_{x_\nu}, W)$ with $g_{\nu}|_{S_i\cap B_{x_\nu}} = f_i|_{S_i \cap B_{x_\nu}}$ for all $i$, then setting $f = \sum_{\nu} \rho_{\nu} g_{\nu}$ completes the proof.

By choosing a local trivialization of $W$ in $B_{x_\nu}$, choosing an origin and fixing coordinates, we can restrict our attention to $W = \R$.

Thus assume $x=x_\nu\in S_I$ for some maximal $I\subset \{1,\dots, N\}$. If $I = \emptyset$ we set $g_\nu = 0$, so we will assume $I \neq \emptyset$ in the following. Using the exponential map at $x$, we identify $B_x$ with an open neighborhood of the origin in $T_xS^{n-1}$, so that $S(E_i)$ corresponds to $F_i=T_xS(E_i)$ for all $i\in I$. Identifying $T_xS^{n-1}=\R^n/\Span\{x\}$, we have $F_i = E_i / \Span \{x\}$ for all $i\in I$. As by assumption $x\in\cap_{i\in I}E_i$, the arrangement $\{E_i\}_{i \in I}$ is minimally intersecting in its span, and by Lemma~\ref{lem:min_intersecting}(ii), the arrangement $\{F_i\}_{i\in I}$ is minimally intersecting in $\sum_{i \in I} F_i$.

Putting $m=n-1$, we have reduced the problem to constructing a common extension $g \in C^\infty(\R^{m})$ of compatible functions $f_i \in C^\infty(F_i)$, $i \in I$, from a minimally intersecting arrangement $\{F_i\}_{i \in I}$  of subspaces of $\R^{m}$.

If $\sum_{i \in I} F_i$ is a proper subspace of $\R^m$, we will first extend to $\sum_{i \in I} F_i$ and then choose an arbitrary extension to $\R^{m}$. Thus assume $\sum_{i \in I} F_i = \R^{m}$.

As the arrangement is minimally intersecting, the sum $\sum_{i \in I} F_i^\perp$ is direct. It will now be convenient to fix a Euclidean structure on $\R^m$ for which this sum is orthogonal.

For $J \subset I$, $J \neq \emptyset$, consider the orthogonal projection
$
	\pi_{J}: \R^{m} \to \cap_{j\in J}F_j
$
 Set $f_J := f_j \circ \pi_J$, where $j \in J$ can be chosen arbitrarily by the compatibility assumption on $f_j$. The extension $g$ is now defined by the following inclusion-exclusion formula
\begin{align*}
	g := \sum_{\emptyset \neq J \subset I} (-1)^{|J| + 1} f_J.
\end{align*}
One easily verifies that $g|_{F_i}$ coincides with $f_i$, concluding the proof. 
\endproof
\noindent Next, define the restriction operator $\res_E$ by
\begin{align*}
	(\res_E \omega)_p (v_1, \dots, v_s) = \omega_p(v_1, \dots, v_s)|_{\wedge^t E}, p \in S(E), v_1, \dots, v_s \in T_p(S(E)),
\end{align*}
where $\omega \in \Omega^s(S(V), \wedge^t V^\ast)$ and $E\subset V$ is a subspace, that is, the differential form $\res_E \omega \in \Omega^s(S(E), \wedge^t E^\ast)$ is the restriction of $\omega$ as a form to $S(E)$, composed with the restriction $\wedge^t V^\ast \to \wedge^t E^\ast$. We will use the same notation for the analogous operation on sections of various vector bundles.
\begin{Proposition}\label{prop:common_lift_extension}
	Suppose that $\{E_i\}_{i=0}^N$ is a semi-generic arrangement of subspaces in $V=\R^n$ and let $s,t \geq 0$. Let $\mathcal E_I$ denote the vector bundle $\wedge^s T_\xi^*S_I\otimes \wedge^t E_I^*$ over $S_I$, and $\mathcal E$ the vector bundle $\wedge^s T_\xi^*S^{n-1}\otimes \wedge^t V^*$ over $S^{n-1}$. Assume $g_i\in\Gamma(S_i, \mathcal E_i)$ are given for all $1\leq i\leq N$, as well as $f_0\in\Gamma(S_0, \mathcal E)$, which are all pairwise compatible on intersections: $\res_{E_i\cap E_0}f_0=\res_{E_i\cap E_0}g_i$, $\res_{E_i\cap E_j}g_i=\res_{E_i\cap E_j} g_j$. Then there exists $f\in \Gamma(S^{n-1}, \mathcal E)$ such that $\res_{E_i}f=g_i$, and $f|_{S_0}=f_0$.
\end{Proposition}
\proof

We will construct elements $f_i\in\Gamma(S_i, \mathcal E)$, $1\leq i\leq N$, that lift the corresponding $g_i$, such that $f_i|_{S_i\cap S_j}=f_j|_{S_i\cap S_j}$ for all $0\leq i<j\leq N$. 

If all $S_i$ are pairwise disjoint for $i=0,\dots, N$, we may choose arbitrary lifts $f_i$ of $g_i$ for all $1\leq i\leq N$, and then choose any $f\in \Gamma(S^{n-1}, \mathcal E)$ restricting to $f_i$ over $S_i$. Thus we assume henceforth that the spheres $S_i$ are not pairwise disjoint.

Now, consider the collections $\mathcal{B}_j$ of subsets $I \subset \{0, \dots, N\}$ of size $j$ such that $E_I\neq \{0\}$, and note that for each $I\in \mathcal B_j$, the family $\{E_{i}\}_{i \in I}$ is minimally intersecting in its span.

Let $m$ be the maximal $j$ such that $\mathcal{B}_j \neq \emptyset$. By assumption, $m\geq 2$. We will define, for all $I\in \mathcal B_j$, sections $f^I\in \Gamma(S_I,\mathcal E)$ that lift $g_i$ on $S_i\cap S_I$, for all $i\geq1$, and coincide with $f_0$ over $S_I\cap S_0$. We will do so by induction on $j$, starting with $j=m$.

Choose $I \in \mathcal{B}_m$. We define $f^I\in \Gamma(S_I, \mathcal E)$ as follows. If $0\in I$, set $f^I=f_0|_{S_I}$. 

If $0\notin I$, using that the arrangement of subspaces $T_\xi S_{i} \cong E_{i}/\Span\{\xi\}$, $i \in I$, is minimally intersecting by Lemma~\ref{lem:min_intersecting}(ii), we can apply Corollary~\ref{cor:smoothDependenceLift} to obtain linear maps $L_\xi$ that depend smoothly on $\xi \in S_{I}$ and map compatible forms in $\wedge^{s}T_\xi^\ast S_{i} \otimes \wedge^t E_{i}^\ast$, $i \in I$, to a form in $\wedge^{s}T_\xi^\ast S^{n-1} \otimes \wedge^t V^\ast$ simultaneously lifting them. We define $f^I\in \Gamma(S_I, \mathcal E)$ by $f^I(\xi) = L_\xi((g_{i}(\xi))_{i \in I})$, $\xi \in S_{I}$.

Assume for $j\geq 1$ that $f^I \in \Gamma(S_{I }, \mathcal E)$ is constructed for all $I \in \mathcal{B}_{j+1}$. We need to construct $f^I$ for all $I \in \mathcal{B}_j$.

Take $I \in \mathcal{B}_j$. If $0\in I$, we set as before $f^I=f_0|_{S_I}$. Otherwise, consider the restriction map
\begin{align*}
	\pi_{I}: \wedge^{s}T^*S^{n-1}\otimes \wedge^t V^* \to \oplus_{i\in I}\wedge^{s}T^*S_{i}\otimes \wedge^t E_{i}^*,
\end{align*}
and the associated  affine sub-bundle $\mathcal F_I$ of the linear bundle $\mathcal E$ over $S_{I}$, with fiber
\begin{align*}
	\mathcal F_I|_{\xi}=\pi_I^{-1}((g_{i}(\xi))_{i \in I}).
\end{align*}
By Corollary \ref{cor:lift_double_forms}, $\mathcal{F}_I$ is non-empty. It holds by construction that $f^{I'}\in\Gamma(S_{I'}, \mathcal F_I)$ for all  $I'\in \mathcal B_{j+1}$ such that $I'\supset I$. Moreover, $(E_{I'})_{I\subset I'\in \mathcal B_{j+1}}$ is semi-generic by Lemma \ref{lem:min_intersecting}(iii). Now use Proposition \ref{prop:smooth_extension} to extend the sections $(f^{I'})_{I\subset I'\in \mathcal B_{j+1}}$ to a global section $f^I$ of $\mathcal F_I$ over $S_I$. 

This completes the induction. The sections $f_i=f^{\{i\}}$ then have the desired properties.
An application of Proposition \ref{prop:smooth_extension} now completes the proof.
\endproof
By letting $E_0$ be a line intersecting $\cup_{i=1}^N E_i$ trivially, we arrive at the following more basic statement. We use the notation of Proposition \ref{prop:common_lift_extension}. 
\begin{Proposition}\label{prop:common_lift_extension0}
	Suppose that $\{E_i\}_{i=1}^N$ is a semi-generic arrangement of subspaces in $V=\R^n$. Assume $g_i\in\Gamma(S_i, \mathcal E_i)$ are given for all $1\leq i\leq N$ which are pairwise compatible on intersections: $\res_{E_i\cap E_j} g_i=\res_{E_i\cap E_j} g_j$. Then there exists $f\in \Gamma(S^{n-1}, \mathcal E)$ such that $\res_{E_i}f=g_i$ for all $1\leq i\leq N$.
\end{Proposition}

\medskip

We are now ready to prove that a family of valuations over an admissible arrangement, compatible under pushforwards, can be represented by an aligned family of $(n-1)$-forms. In order to state this without reference to valuations, we need the following definition.

\begin{Definition}
	For an exact form $\omega\in\Omega^k(S^{n-1}, W)$ and a $k$-dimensional oriented subspace $E$, we define $\Res_E\omega\in W$ by $\int_{S(E)}\eta$, where $\eta$ is any form with $d\eta=\omega$. 
\end{Definition}
Let us point out that $\Res_E \omega$ of an exact form $\omega$ is well-defined. Indeed, suppose that $\eta_1, \eta_2 \in \Omega^{k-1}(S^{n-1}, W)$ are forms with $d \eta_{i} = \omega$, $i=1,2$, that is, $\eta_1 - \eta_2$ is closed. As $S(E)$ is a boundary, it follows that $\int_{S(E)} \eta_1 = \int_{S(E)} \eta_2$.

Note also that the equality $\Res_E \omega= \Res_E \omega'$ for two forms $\omega$, $\omega'$ is meaningful without the need to specify an orientation on $E$.

\medskip

\begin{Proposition}\label{prop:compatible_primitives}
	Fix integers $t\geq s\geq 1$. Let $S=\{E_i\}_{i=1}^N$ be a semi-generic arrangement of distinct subspaces in $V=\R^n$. Let $\omega_i\in \Omega^s(S_i, \wedge^t E_i^*)$, $1\leq i\leq N$, be exact forms such that $\res_{E_{ij}}\omega_i=\res_{E_{ij}}\omega_j$, for all $i,j$. 
	Assume also that whenever $\dim E_{ij}=s$, one has
		\begin{equation}\label{eq:extra_aligned}
			(\Res_{E_{ij}}\omega_i)|_{\wedge^t E_{ij}}=(\Res_{E_{ij}}\omega_j)|_{\wedge^t E_{ij}}.
		\end{equation}
		Furthermore, assume one of the following
		\begin{enumerate}
			\item \textbf{Bounded intersections.} For any subset $S'\subset S$ of $s+1$ subspaces one has $\dim(\cap_{E\in S'} E)\leq t-1$.
			\item\textbf{Minimal intersections.} Any subset $S'\subset S$ with $|S'|\leq s+t+2$ is minimally intersecting in $E_1+\dots+E_N$.
			\item{$\mathbf{N=3}$.} $N=3$.
		\end{enumerate} 
	Then one can choose $\eta_i\in \Omega^{s-1}(S(E_i), \wedge^tE_i^*)$ such that $\res_{E_{ij}}\eta_i=\res_{E_{ij}}\eta_j$ and $d\eta_i=\omega_i$ for all $i,j$.
\end{Proposition}
\noindent Whenever two forms $\eta\in\Omega(S(E), \wedge^t E^*)$ and $\eta'\in \Omega(S(F), \wedge^t F^*)$ satisfy $\res_{E\cap F}\eta=\res_{E\cap F}\eta'$, we will call them aligned, or compatible. We will denote $d_I=\dim E_I$. 
\proof
The proof follows the same lines in all cases, with a few key differences. In the case $N=3$, we may assume that $E_1\cap E_2\cap E_3 = \{0\}$, or else this is a special case of the minimal intersections setting. Furthermore, in this case, the result follows from the bounded intersection case for $s \geq 2$, so we only consider it below as a separate case when $s=1$.
\medskip

\textsl{Step 0.} Fix arbitrary forms $\eta_i$ with $d\eta_i=\omega_i$. We will see how to modify the $\eta_i$ by adding closed forms to make them aligned. By assumption,
\begin{equation}\label{eq:closed_difference}
	d(\res_{E_{ij}}\eta_i-\res_{E_{ij}}\eta_j)=\res_{E_{ij}}\omega_i-\res_{E_{ij}}\omega_j=0.
\end{equation}

Consider first $s=1$. In the case of bounded intersections, there is nothing to prove since $\wedge^tE_{ij}^*=\{0\}$, and we pass to the other cases.

Put $R_{ij}=\res_{E_{ij}}\eta_i-\res_{E_{ij}}\eta_j\in C^\infty(S_{ij}, \wedge^tE_{ij}^*)$. We then have $dR_{ij}=0$ by eq.~\eqref{eq:closed_difference}. Thus if $d_{ij}\geq 2$, $S_{ij}$ is connected so $R_{ij}$ is constant. Also if $d_{ij}=1$, eq.~\eqref{eq:extra_aligned} implies that $R_{ij}$ assumes equal values at both points of $S_{ij}$. In the following, we will therefore tacitly consider $R_{ij}$ as an element of $\wedge^tE_{ij}^*$.

Observe that one has 
\begin{align*}
	&R_{ij}|_{\wedge^t E_{ijk}}-R_{ik}|_{\wedge^t E_{ijk}}+R_{jk}|_{\wedge^t E_{ijk}}\\
	&\equiv\res_{E_{ijk}}\eta_i-\res_{E_{ijk}}\eta_j-(\res_{E_{ijk}}\eta_i-\res_{E_{ijk}}\eta_k)+\res_{E_{ijk}}\eta_j-\res_{E_{ijk}}\eta_k=0.
\end{align*}
In both the minimal intersections case and the $N=3$ case, we may by Corollary~\ref{cor:first_cohomology} find $L_i\in \wedge^t E_i^*$ such that $R_{ij}=L_j|_{\wedge^t E_{ij}}-L_i|_{\wedge^t E_{ij}}$. Now replace $\eta_i$ with $\eta_i+L_i$ for all $i$, concluding this case.

\medskip

\textsl{Step 1.} Assume henceforth $s>1$. By eq. \eqref{eq:closed_difference} we may write
$$\res_{E_{ij}}\eta_i-\res_{E_{ij}}\eta_j=df_{ij},\quad f_{ij}\in\Omega^{s-2}(S_{ij}, \wedge^t E_{ij}^*),$$
which for $1<s<d_{ij}$ follows from the vanishing of the de Rham cohomology, and for $d_{ij}=s$ follows from the additional assumption on $\Res_{E_{ij}}\omega_i$ and $\Res_{E_{ij}}\omega_j$. Note that $df_{ij}=0$ for dimensional reasons if $s>d_{ij}$.

Assume $\eta_1,\dots, \eta_{p-1}$ are already pairwise aligned, so that we may assume $f_{ij}=0$ for $1\leq i<j< p$. We will align $\eta_{p}$ with $\eta_q$ for all $1\leq q< p$, by adding a closed form to $\eta_p$. If $t>d_{qp}$ then $\wedge^tE_{qp}^*=0$, so there is nothing to prove. Thus assume $t\leq d_{qp}$.

Assume that $\eta_p$ is already aligned with $\eta_1,\dots, \eta_{q-1}$. We will align $\eta_p$ with $\eta_q$, by replacing $\eta_p$ with $\eta_p+df$, where $f\in\Omega^{s-2}(S_{p}, \wedge^t E_{p}^*)$ satisfies 
\begin{itemize}
	\item $\res_{E_{ip}} df = 0 $ for $1 \leq i < q$ with $E_{ip} \neq \{0\}$, so that alignment of $\eta_p$ with $\eta_i$, $i<q$ is maintained; and
	\item $\res_{E_{qp}}df=df_{qp}$, aligning $\eta_p$ with $\eta_q$.
\end{itemize}
We will in fact produce closed forms $h_{ip}\in\Omega^{s-2}(S_{ip}, \wedge^t E_{ip}^*)$, which are aligned among themselves and with $f_{qp}$, and then use Proposition \ref{prop:common_lift_extension0} to construct $f$ such that
\begin{itemize}
	\item $\res_{E_{ip}}f=h_{ip}$ for $1 \leq i < q$ with $E_{ip} \neq \{0\}$, and
	\item $\res_{E_{qp}}f=f_{qp}$.
\end{itemize}

Denote $B_{qp}=\{1\leq i\leq q-1:  E_{iqp}\neq \{0\} \}$. For each $i\in B_{qp}$, we define $f_{iqp}=\res_{E_{iqp}}f_{qp} \in \Omega^{s-2}(S_{iqp}, \wedge^t E_{iqp}^*)$. Since $\eta_q$ and $\eta_p$ are aligned with all $\eta_i$, $1\leq i\leq q-1$, it holds that
\begin{align*}
	df_{iqp}=\res_{E_{iqp}}df_{qp}=\res_{E_{iqp}}\eta_q-\res_{E_{iqp}}\eta_p=\res_{E_{iqp}}\eta_i-\res_{E_{iqp}}\eta_i=0,
\end{align*}
and so $f_{iqp}$ is closed. Note that $f_{iqp}=0$ if $t>d_{iqp}$ since  $\wedge^tE_{iqp}^*=0$.

\medskip

\textsl{Step 2.} Assume now $s=2$. It follows that $f_{iqp}\in C^\infty(S_{iqp}, \wedge^t E_{iqp}^*)$ is constant, since if $d_{iqp}\geq 2$ then $S_{iqp}$ is connected, while for $d_{iqp}=1$ we have $d_{iqp}<s\leq t$ and so $f_{iqp}=0$. In the following, we will therefore consider $f_{iqp}$ as element of $\wedge^t E_{iqp}^*$.

We next define an element $h_p\in \wedge^t E_p^*$ such that $h_p|_{\wedge^t E_{iqp}}=f_{iqp}$ for all $i\in B_{qp}$, as follows. In the case of bounded intersections, we have $f_{iqp}=0$ for all $i\in B_{qp}$ since $d_{iqp}\leq t-1$, and we set $h_p=0$.

Consider now the case of minimal intersections. We note that for all $i,j\in B_{qp}$
\begin{equation}\label{eq:quadrupatible}
	f_{iqp}|_{\wedge^t E_{ijqp}}\equiv\res_{E_{ijqp}}f_{qp}\equiv f_{jqp}|_{\wedge^t E_{ijqp}}.
\end{equation}
Observe that $\{E_{iqp}\}_{i\in B_{qp}}$ have the property that any $t+2$ of them are minimally intersecting by Lemma~\ref{lem:min_intersecting}(iii), and it follows from Lemma \ref{lem:multiple_hyperplanes} and \eqref{eq:quadrupatible} that $h_{p}$ exists with $f_{iqp}=h_p|_{\wedge^t E_{iqp}}$ for all $i \in B_{qp}$. 
For trivial reasons, $f_{iqp} = 0 = h_p|_{\wedge^t E_{iqp}}$ also for $i \notin B_{qp}$.

We then set $h_{ip}:=h_p|_{\wedge^t E_{ip}}$ for all $1\leq i\leq q-1$. It holds by construction that
$h_{ip}|_{\wedge^t E_{ijp}}=h_{p}|_{\wedge^t E_{ijp}}=h_{jp}|_{\wedge^t E_{ijp}}$ for $1\leq i<j\leq q-1$, while
\begin{align}\label{eq:prfPropAlignFormsHipComp}
	h_{ip}|_{\wedge^t E_{iqp}} = h_p|_{\wedge^t E_{iqp}} = f_{iqp} \equiv \res_{E_{iqp}}f_{qp},
\end{align}
for all $1 \leq i \leq q-1$. This completes the construction of $h_{ip} \in \wedge^t E_{ip}^*$. Note that when considered as (constant) element of $\Omega^0(S_{ip}, \wedge^t E_{ip}^\ast)$, $h_{ip}$ is closed.

We may now construct the desired $f \in C^\infty(S_p, \wedge^t E_p^\ast)$ using Proposition \ref{prop:common_lift_extension0}, and thus the claim for $s=2$ follows. 

\medskip

\textsl {Step 3.} We next assume $s\geq 3$, and proceed by induction on $s$.

If $ 3\leq s\leq d_{iqp}$, $df_{iqp}=0$ implies that $f_{iqp}$ is exact, and we may write $f_{iqp}=dg_{iqp}$, $g_{iqp}\in\Omega^{s-3}(S_{iqp}, \wedge^tE_{iqp}^*)$. 

For $s\geq d_{iqp}+1$, we have $t>d_{iqp}$ and so $\wedge^t E_{iqp}^*=0$. In that case write $f_{iqp}=dg_{iqp}$ with $g_{iqp}=0$.

As the $f_{iqp}\in\Omega^{s-2}(S_{iqp}, \wedge^t E_{iqp}^*)$, $i\in B_{qp}$, are all restrictions of $f_{qp}$, they are aligned. Note that eq. \eqref{eq:extra_aligned} is trivially satisfied, since $d_{ijqp}=s-2$ implies $d_{ijqp}<t$, and so $\Res_{E_{ijqp}}f_{iqp}|_{\wedge^t E_{ijqp}}=\Res_{E_{ijqp}}f_{jqp}|_{\wedge^t E_{ijqp}}=0$. 

In the case of bounded intersections, the subspaces $(E_{iqp})_{i\in B_{qp}}$ are semi-generic by Lemma~\ref{lem:min_intersecting}(iii), and any $s-1$ of them have at most $(t-1)$-dimensional intersection. 
In the case of minimal intersections, any $s+t$ subspaces among $(E_{iqp})_{i\in B_{qp}}$ are minimally intersecting in $E_{qp}$ by Lemma~\ref{lem:min_intersecting}(iii).
Thus in both cases we may use the induction assumption to choose $g_{iqp}$ that are aligned.

For $j\geq 1$, let $\mathcal B_j$ be the collection of subsets $I\subset\{1,\dots, q\}$ of size $|I|=j$ such that $E_I\cap E_p\neq \{0\}$. Let $m$ be the maximal $j$ with $\mathcal B_j\neq \emptyset$.

If $m=1$, we may find a neighborhood $U$ of $S_{qp}$ that is disjoint from all $S_i$ with $1\leq i\leq q-1$. We then take $h_{ip}=0$ for $1 \leq i < q$, and $f$ an arbitrary lift and extension of $f_{qp}$ that is supported inside $U$.

Now assume $m\geq 2$. We will inductively define a family of forms $g_{I \cup p} \in \Omega^{s-3}(S_{I \cup p}, \wedge^t E_{I \cup p}^\ast)$, indexed by $\{q\}\neq I \in \mathcal{B}_j$, $j\geq 1$, such that the $g_{I \cup p}$ are aligned, and $g_{I \cup p} = \res_{E_{I \cup p}} g_{iqp}$, whenever $i, q \in I, i \neq q$, leading to a definition of $g_{ip}$, $1\leq i<q$ that lift and extend the forms $g_{iqp}$.

For all $I\in \mathcal B_m$, if $q\in I$, choose any $1\leq i\leq q-1$ with $i\in I$ and define $g_{I\cup p}=\res_{E_{I\cup p}} g_{iqp}$. If $q\notin I$, set $g_{I\cup p}=0$. Note that $(g_{I\cup p})_{I\in\mathcal B_m}$ is an aligned family of forms. 

Assume $g_{I\cup p}$ is constructed for all $I\in \mathcal B_r$ with $r\geq j+1$ and $j \geq 1$. For $I\in \mathcal B_j$, $I \neq \{q\}$, if $q\in I$, set $g_{I\cup p}=\res_{E_{I\cup p}}g_{iqp}$ for any $i\in I$ with $1\leq i\leq q-1$.

If $q\notin I$, observe that $(E_{I'\cup p})_{I\subset I'\in\mathcal B_{j+1}}$ is a semi-generic arrangement in $E_{I\cup p}$ by Lemma~\ref{lem:min_intersecting}(iii). If $I \cup \{q\} \not \in \mathcal{B}_{j+1}$, we set $g_{I \cup p} = 0$. Otherwise we have $I \cup \{q\} \in \mathcal{B}_{j+1}$, and choosing an arbitrary $i \in I$ we define
\begin{align*}
	g^{I \cup q} = (\res_{E_{I \cup p}}g_{iqp})|_{S_{I \cup \{q,p\}}} \in \Gamma(S_{I \cup \{q,p\}}, \wedge^{s-3}T_\xi^\ast S_{I \cup p} \otimes \wedge^t E_{I \cup p}^\ast).
\end{align*}
Note that for every  $I'$ such that $I \subset I' \in \mathcal{B}_{j+1}$, $g^{I \cup q}$ is compatible with $g_{I' \cup p}$ by the inductive assumption. Indeed, for $\xi \in S_{I' \cup \{q,p\}}$,
\begin{align*}
	g^{I \cup q}(\xi)|_{\wedge^{s-3}T_\xi S_{I' \cup p} \otimes \wedge^t E_{I' \cup p}} = \res_{E_{I' \cup \{p\}}}g_{iqp}(\xi) = g_{I' \cup p}(\xi).
\end{align*}

By Proposition \ref{prop:common_lift_extension}, we may find \[g_{I\cup p}\in \Gamma (S_{I\cup p}, \wedge^{s-3}T_\xi^\ast S_{I \cup p} \otimes \wedge^t E_{I \cup p}^\ast)=\Omega^{s-3}(S_{I \cup p}, \wedge^tE_{I \cup p}^*)\] which extends $g^{I\cup q}$ and lifts all $g_{I'\cup p}$, $I\subset I'\in\mathcal B_{j+1}$ with $q\notin I'$.

Moreover, by construction, the elements $g_{I \cup p}$, $I \in \mathcal{B}_j$, are aligned. Indeed, let $I_1, I_2 \in \mathcal{B}_j$ and note that if $S_{I_1 \cup I_2 \cup p} = \emptyset$, there is nothing to show. If $S_{I_1 \cup I_2 \cup p} \neq \emptyset$, then there exist sets $J_1, J_2 \in \mathcal{B}_{j+1}$ such that $I_1 \cup I_2 \supset J_k \supset I_k$, $k=1,2$. Thus,
\begin{align*}
	\res_{E_{I_1 \cup I_2 \cup p}} g_{I_1 \cup p} =  \res_{E_{I_1 \cup I_2 \cup p}} g_{J_1 \cup p},
\end{align*}
and, by the inductive assumption, $g_{J_1 \cup p}$ and $g_{J_2 \cup p}$ are aligned. Repeating the steps for $g_{I_2 \cup p}$ then yields the claim.

Arriving at $j=1$, we have our desired elements $g_{ip}$ for $1 \leq i < q$, and we set $h_{ip}=dg_{ip} \in \Omega^{s-2}(S_{ip}, \wedge^t E_{ip}^\ast)$. It holds by construction that
\begin{align*}
	\res_{E_{ijp}}h_{ip}=d\res_{E_{ijp}}g_{ip}=d\res_{E_{ijp}}g_{jp}=\res_{E_{ijp}}h_{jp}
\end{align*}
for $1\leq i<j\leq q-1$, while
\begin{align*}
	\res_{E_{iqp}}h_{ip}=d \res_{E_{iqp}}g_{ip} = dg_{iqp} = f_{iqp}=\res_{E_{iqp}}f_{qp}
\end{align*}
for all $1\leq i\leq q-1$.

As before, we can now construct $f \in \Omega^{s-2}(S_p, \wedge^t E_p^\ast)$ with the aforementioned properties from the sections $h_{ip}$ and $f_{qp}$ using  Proposition \ref{prop:common_lift_extension0}. 
\endproof

\begin{proof}[Proof of Theorem \ref{mthm:extensionFinArrang}]
We will follow the notation from the proof of Proposition~\ref{prop:compatible_primitives}. Moreover, we will write $\phi_{E_i}$ and $\phi_{i}$ interchangeably, and similarly also for various families indexed by $S$.
Observe that when $k\geq \frac{n}{2}$, any $(k+1)$ subspaces would have interesection of dimension at most $n-(k+1)\leq k-1$, thus we need not consider this case separately.

\textsl{Step 0.} Fix a Euclidean structure in $V$. Then the Alesker--Fourier transform allows to define $\Phi_i=\mathbb F\phi_i\in\Val_{d_i-k}^\infty(E_i^*)$. The compatibility condition on intersections then corresponds to $(\pi^i_{ij})_*\Phi_i=(\pi^j_{ij})_{*}\Phi_j$ for all pairs of indices $i,j$, where $\pi^i_{ij}:E_i^*\to (E_i\cap E_j)^*$ is the natural restriction map. We must find $\Phi\in\Val^\infty_{n-k}(V^*)$ such that $\Phi_i=(\pi_i)_*\Phi$ for all $i$, where $\pi_i:V^*\to E_i^*$ is the restriction map. Denote by $\tau_{E_i}=\tau_{\Phi_i}\in\Omega^{k}(S_i, \wedge^k E_i^*)$ the defining form of $\Phi_i$, whenever $d_i>k$.

\textsl{Step 1.} Consider $S_k=\{E\in S: \dim E=k\}$  and $S'=S\setminus S_k$. For all $E\in S'$, we use Proposition \ref{prop:compatible_primitives} to choose $\eta_E\in\Omega^{k-1}(S(E), \wedge^{k}E^*)$ such that $\tau_E=d\eta_E$, and $\res_{E\cap E'}\eta_E=\res_{E\cap E'}\eta_{E'}$ for all $E, E'\in S'$. Note here, that the $\tau_{E}$ are aligned since the valuations $\Phi_i$ are compatible, whereas \eqref{eq:extra_aligned} follows from
\begin{align*}
	\Res_{E_{ij}} \tau_{E_i} = \int_{S_{ij}} \res_{E_{ij}}\eta_{E_i} = ((\pi_{ij}^i)_\ast \Phi_i)(\{0\}) = ((\pi_{ij}^j)_\ast \Phi_j)(\{0\})=\Res_{E_{ij}} \tau_{E_j} ,
\end{align*}
where $E_{i}, E_j \in S'$ with $\dim E_{ij} = k$, and the second equality follows from Proposition~\ref{prop:ConvPushfOnQuadrForms}, and the third from the compatibility condition for the $\Phi_i$.

For all $E\in S_k$, $\Phi_E$ is just a multiple of the Euler characteristic on $E$, and we simply choose arbitrary forms $\eta_E$ representing $\Phi_E$. In all cases under consideration it holds that $E\cap E'$ is a proper subspace of $E$, for all $E,E'\in S$. In particular, $\res_{E\cap E'}\eta_E=\res_{E\cap E'}\eta_{E'}=0$ whenever $E\in S_k$. Thus all forms $\eta_E$, $E\in S$ are aligned.

\textsl{Step 2.} 
Thus we are given compatible forms $\eta_E \in \Omega^{k-1}(S_E, \wedge^kE^\ast)$, $E \in S$, and we may by Proposition  ~\ref{prop:common_lift_extension0} find a form $\eta \in \Omega^{k-1}(S(V), \wedge^k V^\ast)$ with $\res_E \eta = \eta_E$ for all $E\in S$. The corresponding form in $\Omega^{n-k, k-1}(V^*\times S(V))$ then defines the desired valuation $\Phi$. 
\end{proof}
\begin{proof}[Proof of Corollary \ref{cor:hyperplanes}]
	For $k>\frac{n-1}{2}$ we apply the first case of Theorem \ref{mthm:extensionFinArrang}. For $k<\frac{n-1}{2}$, apply the second case.
\end{proof}

Let us show that the additional constraints in the first two cases of the theorem, on top of semi-genericity, are in fact necessary. The following example applies in both cases. We will use the equivalent Alesker-Fourier dual formulation, and fix a Euclidean structure.

\begin{Example}\label{example:counterexample}
	Consider an arrangement $\{E_i\}_{i=1}^N$ of planes in $\R^3$ in general position, that is $E_{ij}:=E_i\cap E_j$ is a line for all $i<j$, and each triple has trivial intersection. Observe that for $N\geq 4$, the dimensions of the sequence
	
	$$ (\R^3)^*	\xrightarrow{d_0} \bigoplus_{i=1}^N E_i^*	\xrightarrow{d_1} \bigoplus_{i<j} E_{ij}^*$$ 
	satisfy
	$$\dim (\R^3)^*-\dim \bigoplus_{i=1}^N E_i^*+\dim \bigoplus_{i<j} E_{ij}^*=3-2N+{N\choose 2}>0.$$
	Thus for $N\geq 4$, the map $d_1$, which factors through $ \bigoplus_{i=1}^N E_i^*/(\R^3)^*$, is not onto. 
	
	Now fix arbitrary elements $R_{ij}\in (E_i\cap E_j)^*$ not in the image of $d_1$, and construct $\Phi_i\in\Val_1^{-,\infty}(E_i)$ given by 
	\begin{align}\label{eq:counterExRepPhii}
		\Phi_i(K)=\int_{S(E_i)}f_i(\theta)dS(K;\theta)
	\end{align}
	inductively as follows. Set $f_1=0$, so $\Phi_1=0$. Assume $f_1,\dots, f_{i-1}$ are constructed. Then, since the $0$-dimensional circles $S(E_{ij})$ do not intersect, we may define $f_{i}$ on $S(E_{ij})$ by $f_j|_{S_{ij}}+R_{ji}$ for all $1\leq j<i$, and then extend it arbitrarily to an odd smooth function on $S(E_i)$.
	
	Assume to the contrary that $\Phi\in \Val_2^{-,\infty}(\R^3)$ exists with $(\pi_i)_*\Phi=\Phi_i$. By a classical result by McMullen~\cite{McMullen1980}, $\Phi$ is given by
	\begin{align*}
		\Phi(K) = \int_{S^2} f(\theta) dS(K, \theta)
	\end{align*} 
	for some odd, continuous function $f$ on $S^2$, and $(\pi_i)_* \Phi$ is given by a formula as in \eqref{eq:counterExRepPhii} with $f|_{S(E_i)}$ replacing $f_i$. It then follows that
	$L_i:=f|_{S(E_{i})}-f_i\in E_i^*$ for all $i$. But then 
	$$L_i|_{S_{ij}}-L_j|_{S_{ij}}=f_j|_{S_{ij}}-f_i|_{S_{ij}}=R_{ij}, $$
	which is a contradiction.
\end{Example}
\begin{Remark}
For $1$-homogeneous even valuations, the conclusion of the theorem holds for any semi-generic arrangement $S$. This follows at once from the Klain embedding, which in this case is surjective, and Proposition \ref{prop:smooth_extension}.
\end{Remark}

	It follows from Theorem \ref{mthm:extensionFinArrang} that for any finite semi-generic subset $S=\{E_i\}_{i=1}^N\subset\Gr_k(\R^n)$, and any choice of densities $\mu_i\in\Dens(E_i)$, one can find a valuation $\phi\in\Val_k^\infty(\R^n)$ with $\phi|_{E_i}=\mu_i$.  The assumption of semi-genericity is in fact superfluous in this case, as the following proposition shows.

\begin{Proposition}\label{prop:cosine_arbitrary}
	For any subset $S=\{E_i\}_{i=1}^N\subset\Gr_k(\R^n)$, and any choice of densities $\mu_i\in\Dens(E_i)$, there is a valuation $\phi\in\Val_k^{+,\infty}(\R^n)$ with $\phi|_{E_i}=\mu_i$. 
\end{Proposition}
\proof
Fix a Euclidean structure on $\R^n$, identifying $\Dens(E)=\R$ for all $E\in\Gr_k(\R^n)$. Denote by $W\subset C^\infty(\Gr_k(\R^n))$ the image of the Klain embedding, which is the image of the cosine transform $\mathcal C_k$ by \cite{Alesker2004c}. 
It then suffices to show that the map $W\to \R^N$ given by $f\mapsto (f(E_i))_{i=1}^N$ is onto. Assuming the contrary, there exist constants $c_i$, $i=1,\dots, N$, not all zero, such that $\sum_{i=1}^N c_i f(E_i)=0$ for all $f\in W$. 

Thus for any $h\in C^\infty$, taking 
\begin{align*}
	f(E)=(\mathcal C_k h) (E)=\int_{\Gr_k(\R^n)} |\cos(E, F)|h(F)dF
\end{align*}
we find that 
\[\int_{\Gr_k(\R^n)}\left(  \sum c_i |\cos(E_i, F)| \right)h(F)dF=0,\qquad\forall h\in C^\infty(\Gr_k(\R^n)),\]
that is
\begin{equation}\label{eq:cosine_combination}
	\sum c_i |\cos(E_i, F)| =0.
\end{equation}

Note that the function $F \mapsto |\cos(E_i, F)|\in C(\Gr_k(\R^n))$ is smooth precisely at the complement of its zero set, which is $\Xi_i=\{F: F\cap E_i^\perp\neq\{0\}\}$. Assuming without loss of generality that $c_N\neq 0$, it follows from eq. \eqref{eq:cosine_combination} that $|\cos(E_n, F)|$ is smooth on the complement of $\Xi_1\cup\dots\cup\Xi_{N-1}$.
That is, $\Xi_N\subset \Xi_1\cup\dots\cup\Xi_{N-1}$.

However $\Xi_i$ has Hausdorff codimension $1$ in $\Gr_k(\R^n)$ for all $i$, while the Hausdorff codimension of $\Xi_i\cap \Xi_j$ in $\Gr_k(\R^n)$ is at least $2$, for all $i\neq j$.  Thus the equality $\Xi_N= (\Xi_1\cap \Xi_N)\cup\dots\cup(\Xi_{N-1}\cap \Xi_N)$ cannot hold. 
\endproof

\section{Extension from compact submanifolds and the Nash theorem}\label{sec:Nash}

In this section, we give a proof of Theorem~\ref{mthm:nashThm}. For this reason, we first show by methods similar to Section~\ref{sec:finArr} that a smooth assignment of valuations to subspaces from a compact submanifold $Z \subseteq \Gr_r(\R^n)$ can be extended to a globally defined valuation if $Z$ has the property that pairs of subspaces from $Z$ intersect trivially, in a slightly stronger sense described in Definition~\ref{def:musical_chairs}.

\medskip

Restricting the bundle $\Val_j^\infty(\Gr_r(\R^n))$ over $\Gr_r(\R^n)$ to a bundle over $Z$, we can consider the space of smooth sections over $Z$, denoted by $V_j^\infty(Z, \R^n)$. Our next result is a version of Theorems~\ref{mthm:resProperty} and \ref{mthm:extensionFinArrang} for perfectly self-avoiding compact submanifolds $Z$, where the proof follows the ideas used in the proof of Theorem~\ref{mthm:extensionFinArrang}. Since if $Z$ is perfectly self-avoiding, there are only trivial intersections, we do not need further compatibility conditions and we can skip the step of aligning the differential forms (corresponding to Proposition~\ref{prop:compatible_primitives}).
\begin{proof}[Proof of Theorem \ref{mthm:nashThm2}]
Let $\phi \in V_j^\infty(Z, \R^n)$. Fix a Euclidean structure on $\R^n$. We first apply the Alesker--Fourier transform $\mathbb{F}$, to obtain a section $\Phi_E=\mathbb F\phi_E\in\Val^{\infty}_{r-j}(E^\ast)$.

Next, we are going to define forms $$\eta_E \in \Omega^{j-1}(S(E), \wedge^j E^\ast)=\Gamma(S(E), \wedge^{j-1}(\xi^\perp)^*\otimes \wedge^j E^\ast), \quad E \in Z,$$ which represent $\Phi_E$ and depend smoothly on $E$.

To this end, first consider the case $j = r$. Then $\Phi_E = f(E) \chi_E$, where  $f: Z \to \R$ is smooth, and we may take $\eta_E$ given by $\eta_E|_\xi=f(E)\vol_{\xi^\perp}\otimes (\vol_{\xi^\perp}\wedge \xi)$, where the orientation of $\xi^\perp$ defining the Euclidean volume form $\vol_{\xi^\perp}$ can be arbitrary.

Now assume $j<r$, and let $\tau_E \in \Omega^j(S(E), \wedge^j E^\ast)$ be the defining forms of $\Phi_E$, $E \in Z$, which depend smoothly on $E$. Let $G_E$ denote the Green operator on $\Omega^{j}(S(E))$, so that $\id=\mathcal H_E+\Delta_E G_E$, where $\Delta_E$ is the Laplace-Beltrami operator on $S(E)$, and $\mathcal H_E$ the orthogonal projection to the space of harmonic forms. Since $\tau_E$ is the defining form of a valuation $\Phi_E$, it is an exact form. In particular, $\mathcal{H}_E \tau_E = 0$, by the Hodge decomposition.

We then set $\eta_E=d^*G_E\tau_E$. Recalling that $d\tau_E=0$ while the image of $d^*$ is $L^2$-orthogonal to the space of closed forms, we have 
\begin{align*}
	\tau_E=\Delta_E G_E\tau_E=dd^*G_E\tau_E+d^*dG_E\tau_E=d\eta_E+d^*(dG_E\tau_E)
\end{align*}
and consequently $d^*dG_E\tau_E=0$ and $\tau_E=d\eta_E$.

It remains to find a form $\eta \in \Omega^{j-1}(S^{n-1}, \wedge^j (\R^n)^\ast)$ with $\res_E \eta = \eta_E$. For this reason, consider $\eta_E$ as a section in $\Gamma(S(E), \wedge^{j-1} T_\xi^\ast S(E) \otimes \wedge^j E^\ast)$ and let
\begin{align*}
	b_\xi: \wedge^{j-1} T_\xi^\ast S(E) \otimes \wedge^j E^\ast \hookrightarrow \wedge^{j-1} T_\xi^\ast S^{n-1} \otimes \wedge^j (\R^n)^\ast,  
\end{align*}
be the natural map induced by the Euclidean structure, $E^\ast \cong E$ and $(\R^n)^\ast \cong \R^n$, and the inclusion $E \subset \R^n$. We may then define $\eta|_\xi = b_\xi(\eta_E|_\xi)$ for every $\xi \in S^{n-1}$ such that $\xi \in S(E)$ for some $E \in Z$. Since $Z$ is perfectly self-avoiding and $\eta_E$ depends smoothly on $E$, $\eta$ is a well-defined smooth section, defined on the embedded submanifold $\theta_r(\mathbb{P}_r|_Z)$ of $S^{n-1}$. We may then extend $\eta$ arbitrarily to obtain $\eta \in \Omega^{j-1}(S^{n-1}, \wedge^j (\R^n)^\ast)$, defining a valuation $\Psi \in \Val_{n-j}^\infty(\R^n)$. Applying the inverse Alesker--Fourier transform, finally yields the sought after valuation $\psi = \mathbb{F}^{-1}\Psi$.
\end{proof}

\medskip

As an application of Theorem \ref{mthm:nashThm2}, we deduce Theorem~\ref{mthm:nashThm}. The key point of the proof of Theorem~\ref{mthm:nashThm} is to find an embedding $M \hookrightarrow \R^n$ for which the tangent spaces of $M$ form a perfectly self-avoiding compact submanifold of $\Gr_m(\R^n)$, leading to the following definition.

\begin{Definition}\label{def:perfskew}
	Suppose that $M$ is a compact smooth manifold. An embedding $e: M \hookrightarrow \R^n$ is \emph{perfectly non-parallel}, if the induced map $\mathbb{P}_+(TM)\to \mathbb{P}_+(\R^n)$ given by $(x,[v])\mapsto [d_xe(v)]$ is an embedding.
\end{Definition}

Note that Definition~\ref{def:perfskew} is similar to the notion of \emph{totally skew}, introduced in \cite{Ghomi2008}, and strengthens the notion of \emph{totally non-parallel}, defined in \cite{Harrison2020}. It is an application of Thom's transversality theorem that perfectly non-parallel embeddings into $\R^n$ exist whenever $n$ is large enough. 
\begin{Lemma}\label{lem:exPerfNonParEmbed}
	For any compact manifold $M^m$ and $n\geq \max(4m+1, {m+1 \choose 2}+2m)$, there exists a perfectly non-parallel embedding $e:M\hookrightarrow \R^n$.
\end{Lemma}
\proof
First note that an embedding $e: M \hookrightarrow \R^n$ is perfectly non-parallel if and only if the map $de: TM\setminus \underline{0} \to \R^n, (x,v) \mapsto de_x(v),$ is an injective immersion. We will show that any totally skew embedding can be perturbed so that $de$ is an injective immersion.

Indeed, by \cite[Prop.~2.3]{Ghomi2008}, there exists a totally skew embedding $e_0: M \hookrightarrow \R^n$ for $n\geq 4m+1$, and $e_0$ remains totally skew under $C^1$-small perturbations. Writing $f = de_0$, in local coordinates, one has
\begin{align*}
	f(x,v)=\sum_{j=1}^m \left.\frac{\partial e_0}{\partial x_j}\right|_xv_j,
\end{align*}
and so
\begin{align*}
	d_{(x,v)}f\left(\frac{\partial}{\partial x_k}\right)=\sum_{j=1}^m\left.\frac{\partial^2e_0}{\partial x_j\partial x_k}\right|_x v_j, \quad \text{ and }\quad 
	d_{(x,v)}f\left(\frac{\partial}{\partial v_k}\right)=\left.\frac{\partial e_0}{\partial x_k}\right|_x.
\end{align*}
It follows that a $C^1$-small perturbation of $e_0$ is perfectly non-parallel if the vectors $\frac{\partial e_0}{\partial x_k}, \frac{\partial^2 e_0}{\partial x_k\partial x_j}$ are all linearly independent. Denoting by $J^2_x$ the space of $2$-jets over a point $x \in M$, by Thom's transversality theorem, we may find a $C^2$-generic perturbation of $e_0$ such that its $2$-jet does not intersect the subset $C$ of linearly dependent vectors, if $m$ is smaller than the co-dimension of $C$. One easily computes that $C$ has codimension $n-(m + {m+1 \choose 2}) + 1$ inside $J^2_x$. As $n-(m+{m+1 \choose 2})+1>m\iff n\geq 2m+{m+1 \choose 2}$ holds, we are done.
\endproof

\medskip

Next, denote by $V_j^\infty(M)$ the Fr\'echet space of smooth global sections of the Fr\'echet bundle over $M$ with fiber $\Val_j^\infty(T_xM)$ over $x \in M$. By \cite{Alesker2006}, the space of smooth valuations on $M$ admits a natural filtration, $\mathcal V^\infty(M)=\mathcal W_0^\infty(M)\supset \mathcal W_1^\infty(M)\supset\dots\supset\mathcal W_{\dim M}^\infty(M)$, such that $V_j^\infty(M)$ is isomorphic to the quotient $\mathcal{W}_j^{\infty}(M)/\mathcal{W}_{j+1}^\infty(M)$.  We write $\phi \mapsto [\phi]_j$ for the composition of the quotient map $\mathcal{W}_j^\infty(M) \to \mathcal{W}_j^{\infty}(M)/\mathcal{W}_{j+1}^\infty(M)$ with this isomorphism. 

\medskip

We are now ready to prove Theorem \ref{mthm:nashThm} in the following precise form.
\begin{Theorem}\label{thm:nash_firstproof}
	Suppose that $M^m$ is a compact manifold, and let $e:M\hookrightarrow \R^n$ be a perfectly non-parallel embedding. Then the image of the restriction map $e^*:\Val^\infty(\R^n)\to \mathcal V^\infty(M)$ is given by $\mathcal W_1^\infty(M)\oplus \Span(\chi)$.
\end{Theorem}
\proof
First note that $e^*\chi=\chi$, and that for any $\psi\in\Val^\infty(\R^n)$, $[e^*\psi]_0\in C^\infty(M)$ is a constant function by translation-invariance. Consequently, for any $\phi\in \mathcal V^\infty(M)$ in the image of $e^*$ we may find $c\in \R$ such that $\phi-c\chi \in\mathcal W_1^{\infty}(M)$.

Next, let $\phi\in \mathcal W_j^\infty(M)$ with $j\geq 1$. By assumption, the collection of tangent spaces $Z=\{T_xM\}_{x\in M}\subset\Gr_m(\R^n)$ is a perfectly self-avoiding compact submanifold, and $[\phi]_j$ defines a smooth assignment of a $j$-homogeneous valuation to every $E \in Z$. By Theorem \ref{mthm:nashThm2} we may therefore choose $\psi_j\in \Val^\infty_j(\R^n)$ such that $[e^\ast\psi_j]_j=[\phi]_j$. Hence, $\phi-e^*\psi_j\in \mathcal W_{j+1}^\infty(M)$.
Starting with $j=1$, and using the above argument repeatedly, we deduce that any $\phi\in \mathcal W_1^\infty(M)$ equals $e^*(\psi_1+\psi_2+\dots+\psi_m)$, concluding the proof.
\endproof

\begin{Remark}
	In light of the results of \cite{Faifman2021}, it seems possible that one can find an embedding $e:M\hookrightarrow \R^N$ such that $\mathcal W_1^\infty(M)=e^*\Val_1^\infty(\R^N)$.
\end{Remark}

	\section{Crofton formulas}\label{sec:crofton}
	In this section, we apply Theorem~\ref{mthm:nashThm} to deduce the existence of Crofton formulas for all smooth valuations on manifolds. To this end, we first construct Crofton formulas for all translation-invariant valuations in a linear space.
	
	\subsection{Translation-invariant valuations}\label{sec:croftFormTransInvVal}
	
	Recall that for every valuation $\phi \in \Val_k^{+,\infty}(V)$ there exists a signed, translation-invariant and smooth measure $m_\phi$ on the affine grassmannian $\AGr_{n-k}(V)$, such that
	\begin{equation}\label{eq:applCroftFormEven_1}
	 \phi(K) = \int_{\AGr_{n-k}(V)} \!\!\!\chi(K \cap E) d m_\phi(E), \quad K \in \mathcal{K}(V).
	\end{equation}
	Such a measure is called a \emph{Crofton measure} for $\phi$. The existence of a Crofton measure is a consequence Alesker's irreducibility theorem~\cite[Thm.~1.3]{Alesker2001} and the Cassel\-man--Wallach theorem~\cite{Casselman1989}. Let us note here that eq.~\eqref{eq:applCroftFormEven_1} also holds with $K$ replaced by $A\in \mathcal{P}(V)$ (see \cite[Lemma 2.4.7]{Alesker2006}).
		
	\bigskip
	
	We next construct a Crofton formula for odd valuations; a similar construction was used in \cite{Alesker2011b} to define the Alesker--Fourier transform for odd valuations. Let us first outline the idea behind the construction. We start with a representation of $1$-homogeneous valuations going back to Goodey and Weil~\cite{Goodey1984}, which is extended to arbitrary degrees of homogeneity by taking Alesker products. Namely, any valuation $\phi_1 \in \Val_1^{\infty}(V)$ can be written as
	\begin{align}\label{eq:GWrepOnehomg}
		\phi_1(K) = \int_{S(V)}h_K(\theta)d\mu(\theta), \quad K \in \mathcal{K}(V),
	\end{align}
	for some unique, signed, smooth measure $\mu$ on $S(V)$ satisfying $\int_{S(V)} \theta d\mu(\theta) = 0$. Here  $h_K$ is the support function of $K \in \mathcal{K}(V)$. The condition on $\mu$ guarantees that $\phi_1$ is translation-invariant. Evidently, $\phi_1$ is odd if and only if $\mu$ is odd. 
	By taking the linear combinations of Alesker products of a $1$-homogeneous odd valuations with $(k-1)$-homogeneous even valuations, we can by the irreducibility theorem approximate any odd $k$-homogeneous smooth valuation $\phi_k \in \Val_k^{\infty}(V)$ (similarly to the proof of Lemma~\ref{lem:product_onto}). To arrive at a general Crofton formula, it is therefore natural to study representations for products of valuations of the form \eqref{eq:GWrepOnehomg}. 
	
	To this end, we note that for every $K \in \mathcal{K}(V)$, the support function of $K$ evaluated at $u \in S(V)$ can be represented as
	\begin{align}\label{eq:sptFctAsChi}
		h_K(u) = \int_0^\infty \chi(K \cap (u^+ + tu)) dt + \int_{-\infty}^0\left( \chi(K \cap (u^+ + tu)) - \chi(K)\right) dt,
	\end{align}
	where $u^+$ denotes the positive half-space defined by $u$. The substracted term $\chi(K)$ in the second integral makes the integrand compactly supported, and so plays a regularizing role. Thus, by \eqref{eq:GWrepOnehomg}, an odd $1$-homogeneous valuation is essentially given as an integral of $\chi(K \cap F^+)$ over all (affine) half-spaces $F^+ \subset V$.
	
	We then take the Alesker product with an even valuation of degree $k-1$, given by the Crofton formula \eqref{eq:applCroftFormEven_1}. The latter is an integral over the affine grassmannian of $(n-k+1)$-dimensional flats, and the Alesker product corresponds to the intersection of those subspaces with the half-spaces in the Crofton formula for $\phi_1$ (see \cite[Appendix B]{Faifman2017}), resulting in the representation of our $k$-homogeneous valuation by the integration over the grassmannian of half-subspaces of dimension $n-k$.

	We will next describe rigorously the resulting assignment of valuation to a measure on the grassmannian of half-subspaces. It will be necessary to make the construction $\GL_n$-equivariant as we will then apply the Casselman--Wallach theorem to deduce the surjectivity of this map. We achieve this by utilizing appropriate $\GL_n$-equivariant line bundles.
	
	\bigskip
	
	Denote by $\widehat{\Gr}_{k}(E)$ the grassmanian of co-oriented $k$-dimensional linear subspaces of a linear space $E$. Denote by $\widehat{\mathcal F}_{i}^{i+1}(V)$ the partial flag manifold of pairs of linear subspaces of $V$, $F = (F_0, F_1)$ with $F_0\subset F_1$, where $\dim F_0 = i$, $\dim F_1 = i+1$, and $F_0$ is co-oriented in $F_1$. For $F=(F_0,F_1)\in \widehat{\mathcal F}_{i}^{i+1}(V)$, let $F_0^{\pm} \subseteq F_1$, respectively $(F_1/F_0)^\pm\subset F_1/F_0$, denote the positively (negatively) oriented half-space, respectively half-line, given by the co-orientation. Consider the space
	\begin{align}\label{eq:croftOddSpacesMeas}
		\widetilde{\mathcal{M}}_j^-=\Gamma^-(\widehat{\mathcal F}_{n-j}^{n-j+1}(V), \Dens( T_F\widehat{\mathcal F}_{n-j}^{n-j+1}(V))\otimes \Dens(V/F_0)) 
	\end{align}
	of all sections that are odd with respect to switching the co-orientation. We will need a condition similar to the one needed for \eqref{eq:GWrepOnehomg} ensuring the resulting valuation is translation-invariant. To write it down, note that for $F=(F_0, F_1) \in \widehat{\mathcal F}_{i}^{i+1}(V)$ one has a natural identification
	\begin{align*}
		\Dens(T_F\widehat{\mathcal F}_{n-j}^{n-j+1}(V))=\Dens(T_{F_1}\Gr_{n-j+1}(V))\otimes \Dens(T_{F_0}\widehat{\Gr}_{n-j}(F_1)),
	\end{align*}
	so for $\mu \in \widetilde{\mathcal{M}}_j^-$, we can write 
	\begin{align*}
		\mu(F_0, F_1)\in \Dens(T_{F_1}\Gr_{n-j+1}(V))\otimes\Dens(T_{F_0}\widehat{\Gr}_{n-j}(F_1))\qquad\\\otimes \Dens(V)\otimes\Dens(F_0)^*.
	\end{align*}
	Given any $z\in F_1$, we consider its projection $\proj_{F_1/F_0}(z)\in F_1/F_0$, and set
	\begin{align*}
		L_z(F_0,F_1):=\sign(z, F_1/F_0)\proj_{F_1/F_0}(z)\in F_1/F_0,
	\end{align*}
	where $\sign(z, F_1/F_0)\in\{\pm1\}$ according to whether $\proj_{F_1/F_0}(z)$ points in the positive or negative direction of $F_1/F_0$, and $\sign(z, F_1/F_0)=0$ if $z\in F_0$. 
	Using the co-orientation, we may identify
	\begin{align*}
		F_1/F_0=\Dens(F_1/F_0)^*=\Dens(F_0)\otimes \Dens(F_1)^*,		
	\end{align*}
	that is, $L_z(F_0,F_1) \in \Dens(F_0)\otimes \Dens(F_1)^*$. Then we can pair $L_z(F_0,F_1)$ with $\mu(F_0, F_1)$ and
	\begin{align*}
		\langle \mu(F_0, F_1),  L_z(F_0,F_1)\rangle \in  \Dens(T_{F_1}\Gr_{n-j+1}(V))\otimes\Dens(T_{F_0}\widehat{\Gr}_{n-j}(F_1))\qquad\\\otimes \Dens(V/F_1).
	\end{align*}
	Setting
	\begin{align*}
		\langle \mu(\bullet, F_1), z\rangle_{ \widehat{\Gr}_{n-j}(F_1)}&:=\int_{\widehat{\Gr}_{n-j}(F_1)}\langle \mu(F_0, F_1),  L_z(F_0,F_1)\rangle\\
		&\in \Dens(T_{F_1}\Gr_{n-j+1}(V))\otimes \Dens(V/F_1),
	\end{align*}
	we define $\mathcal{M}_j^-$ as the subspace of all $\mu$ in $\widetilde{\mathcal{M}}_j^-$ such that for every $F_1\in\Gr_{n-j+1}(V)$ and $z\in F_1$ one has
	\begin{align}\label{eq:croftOddCondMeasMoment}
		\langle \mu(\bullet, F_1), z\rangle_{ \widehat{\Gr}_{n-j}(F_1)}=0.
	\end{align}
	Note that for $j=1$, \eqref{eq:croftOddCondMeasMoment} reduces to the condition after \eqref{eq:GWrepOnehomg}. 
	
	Next, we define $\eta(K,F,y) \in \Dens(F_1/F_0)^\ast$, for a convex body $K$, $F = (F_0, F_1) \in \widehat{\mathcal F}_{n-j}^{n-j+1}(V)$ and $y \in V$, by setting for any $dx\in \Dens(F_1/F_0)$
	\begin{align*}
		\langle \eta(K, F, y), dx\rangle &= \int_{x\in (F_1/F_0)^+}\chi((F_0^++y+x)\cap K)dx\\
		&+\int_{x\in (F_1/F_0)^-} \left(\chi((F_0^++y+x)\cap K)-\chi((F_1+y)\cap K)\right)dx.
	\end{align*}
	Note that $\eta(K,F,y)$ is well-defined, as $F_0^+ + x$ does not depend on the choice of $x \in F_1/F_0$, and the integrands are compactly supported. By \eqref{eq:sptFctAsChi}, and after choosing a Euclidean structure, $\langle \eta(K,F,y), dx\rangle$ coincides with $h_{(K-y) \cap F_1}(z)$, where $F_1 = F_0 \oplus \Span z$ and $z\in F_0^\perp$ is a positively oriented unit vector.
\begin{Lemma}\label{lem:croftOddEx}
	Let $0 < j < n$, and suppose that $\mu \in \mathcal{M}_j^-$. Then 
	\begin{equation}\label{eq:croftFormOdd}
		\phi_{j}^{\mu}(K)=\int_{F_1\in\Gr_{n-j+1}(V)} \int_{y\in V/F_1} \int_{F_0 \in \widehat{\Gr}_{n-j}(F_1)}\!\!\!\!\!\!\!\!\eta(K, (F_0, F_1), \overline y) d\mu(F_0, F_1),
	\end{equation}
	where $\overline y\in V$ is an arbitrary lift of $y\in V/F_1$ that is chosen independently of $F_0$, defines a valuation in $\Val_j^{-, \infty}(V)$.
\end{Lemma}
\proof
First note that \eqref{eq:croftFormOdd} is well-defined. Indeed, since
\begin{align*}
	&\eta(K, (F_0, F_1), \overline{y})\otimes \mu(F_0, F_1)\\\in &\Dens(T_{(F_0,F_1)}\widehat{\mathcal F}_{n-j}^{n-j+1}(V))\otimes\Dens(V/F_0)\otimes\Dens(F_1/F_0)^*\\
	&=\Dens(T_{F_1}\Gr_{n-j+1}(V))\otimes\Dens(T_{F_0}\widehat{\Gr}_{n-j}(V))\otimes\Dens(V/F_1),
\end{align*}
the integral makes formal sense. As for $z \in F_1$,
\begin{align*}
	\eta(K, (F_0, F_1), y+z) = \eta(K, (F_0, F_1), y) + L_z(F_0,F_1),
\end{align*}
whenever $K \cap (F_1+y) \neq \emptyset$, the innermost integral in eq. \eqref{eq:croftFormOdd} does not depend on the choice of $\overline y$ lifting $y\in V/F_1$ by eq. \eqref{eq:croftOddCondMeasMoment} and the assumption that $\overline{y}$ is independent of $F_0$. Moreover, since $K$ is compact, the support of $\eta(K, (F_0, F_1),\cdot)$ is bounded uniformly in $(F_0, F_1)$, so that the integral is finite.

Next, suppose that $K_k \to K \in \mathcal{K}(V)$ converges in the Hausdorff metric. Then, $(F_0^++y+x)\cap K_k \to (F_0^++y+x)\cap K$ for all $y,x$ such that either $(F_0^++y+x)\cap K = \emptyset$ or $(F_0^++y+x)\cap \mathrm{int}\, K \neq \emptyset$. In the first case, $K$ and $F_0^++y+x$ can be separated by two hyperplanes with positive distance by \cite[Thm.~1.3.7]{Schneider2014}, so $(F_0^++y+x)\cap K_k = \emptyset$ for $k$ large enough, while in the second case $K$ and $F_0^++y+x$ cannot be separated by a hyperplane and the claim follows by \cite[Thm.~1.8.10]{Schneider2014}. Consequently, as $\chi$ is continuous, the integrand in $\eta(K,F,y)$ converges pointwise almost everywhere, so by, dominated convergence, $\phi_{j}^\mu$ is continuous.

It follows at once that $\phi_j^{\mu}\in \Val(V)$, and it is straightforward to verify that $\phi_j^\mu$ is $j$-homogeneous and odd.

It remains to show that $\phi_j^{\mu}$ is a smooth valuation. To this end, we observe that since \eqref{eq:croftOddCondMeasMoment} is a $\GL(V)$-invariant condition, the space $\mathcal{M}_j^-$ is a closed $\GL(V)$-invariant subspace of the Fr\'echet space of sections 	$\widetilde{\mathcal{M}}_j^-$ \eqref{eq:croftOddSpacesMeas}. As the map $\mu \mapsto \phi_j^{\mu}$ is $\GL(V)$-equivariant, linear and continuous, it maps smooth vectors to smooth vectors, that is, $\phi_j^{\mu} \in \Val_{j}^{-,\infty}(V)$.
\endproof

In the following, we will call any $\mu\in\mathcal M_j^-$ for which $\phi=\phi_j^\mu \in \Val_{j}^{-,\infty}(V)$ a \emph{Crofton measure} for $\phi$, in analogy with the even case. We will next show that every $\phi \in \Val_j^{-,\infty}(V)$ admits a Crofton measure.

\begin{Proposition}\label{prop:exCroftMeasOddValonR}
	Suppose that $0<j<n$ and let $\phi \in \Val_j^{-,\infty}(V)$. Then there exists $\mu \in \mathcal{M}_j^-$ such that $\phi = \phi_{j}^{\mu}$.
\end{Proposition}
\proof
Observe that the map $\mathcal M_j^-\to \Val_j^{-,\infty}(V)$, $\mu \mapsto \phi_j^\mu$ is $\GL(V)$-equi\-variant, and both spaces are Fr\'echet representations of moderate growth. Provide we can show it is nonzero, then Alesker's irreducibility theorem~\cite[Thm.~1.3]{Alesker2001} implies that its image is dense in $\Val_j^{-,\infty}(V)$, and, by the theorem of Casselman--Wallach~\cite{Casselman1989}, the image is also closed, whence the map is surjective and the claim follows.

We will describe a delta measure $\mu$ for which $\phi_j^\mu \neq 0$. Convolving with an approximate identity on $\GL(V)$, one readily obtains a smooth measure with the same property. Start by fixing a Euclidean structure on $V=\R^n$ and a subspace $F_1=\R^{n-j+1}$. We may now identify $\widehat{\Gr}_{n-j}(F_1)=S^{n-j}$ and $\Dens(V/F_0)=\R$.

Let $K\subset \R^n$ be the regular $n$-dimensional simplex, rotated in such a way that $F_1$ is parallel to one of its $(n-j+1)$-dimensional faces. It follows that, whenever $(F_1+y)\cap K$ has non-empty interior, it is homothetic to a fixed regular $(n-j+1)$-dimensional simplex $S_{n-j+1}$ in $F_1$. Choosing $\Delta\subset S^{n-j}$ to be the vertices of that simplex, we may set 
\begin{align*}
	\mu=\sum_{p\in \Delta}\delta_{(p, F_1)}-\sum_{p\in \Delta}\delta_{(-p, F_1)}.
\end{align*}
One readily checks that for $F_1' \in \Gr_{n-j+1}(\R^n)$,  $\langle \mu(\bullet, F_1'), z\rangle_{ \widehat{\Gr}_{n-j}(F_1')}=0$, either trivially if $F_1'\neq F_1$, or since 
\begin{align*}
	\langle \mu(\bullet, F_1), z\rangle_{ \widehat{\Gr}_{n-j}(F_1)}=2\sum_{p \in \Delta}\langle p, z\rangle = 0,
\end{align*}
that is, $\mu\in\mathcal M_j^{-}$. In order to see that $\phi_j^\mu \neq 0$, observe that for $p\in\Delta$ and $\mathrm{int}\, K \cap (F_1+y) \neq \emptyset$,
\begin{align*}
	\eta(K, (p, F_1), y)-\eta(K, (-p, F_1),y)=h_{(K-y)\cap F_1}(p)-h_{(K-y)\cap F_1}(-p)\\
=c(y)\left(h_{S_{n-j+1}}(p) - h_{S_{n-j+1}}(-p)\right) + 2 \langle z(y),p\rangle = c(y) c_{n-j+1}  + 2 \langle z(y),p\rangle,
\end{align*}
where $c_{n-j+1}>0$ is a constant independent of $p$ as $S_{n-j+1}$ is regular. Here, we used that $(K-y)\cap F_1 = c(y)S_{n-j+1} + z(y)$, for some $c(y)>0$ and $z(y) \in \R^n$. Consequently, setting $c(y) = 0$ if $(K-y)\cap F_1$ has empty interior, and using that $\sum_{p \in \Delta}\langle p, z\rangle = 0$ for all $z \in \R^n$,
\begin{align*}
	\phi_j^\mu(K)&= \sum_{p\in \Delta} c_{n-j+1} \int_{F_1^\perp}c(y)dy +2\int_{F_1^\perp}\sum_{p\in \Delta} \langle z(y),p\rangle dy\\
	 &=c_{n-j+1}(n-j+2)\int_{F_1^\perp}c(y)dy>0,
\end{align*}
which concludes the proof.
\endproof

In the following, we will denote by $\mathrm{HGr}_{n-k+1}(V)$ the manifold of $(n-k+1)$-dimensional affine half-spaces in $V$, $\dim V = n$. Note that $\mathrm{HGr}_{n-k+1}(V)$ can naturally be para\-metri\-zed by quadruplets $(F_0, F_1,y,x)\mapsto F_0^++y+x$, where $(F_0, F_1) \in \widehat{\mathcal F}_{n-k}^{n-k+1}(V)$, $y \in V/F_1$ and $x \in F_1/F_0$.

\begin{Corollary}\label{cor:odd_crofton_compact}
	Take $\phi \in \Val_k^{-,\infty}(V)$, and let $U\subset V$ be a precompact open set. Then there exists a signed, smooth, odd, compactly supported measure $m$ on the grassmannian of half-spaces $\mathrm{HGr}_{n-k+1}(V)$ such that the Crofton integral 
	\begin{equation}\label{eq:applCroftFormEven_2}
		\psi(A): = \int_{\mathrm{HGr}_{n-k+1}(V)} \!\!\!\chi(A \cap H) d m(H)
	\end{equation}
	defines a valuation  $\psi\in\mathcal V^\infty_c(V)$ satisfying $\psi|_U=\phi|_U$.
\end{Corollary}
\proof
 By Proposition~\ref{prop:exCroftMeasOddValonR}, there exists $\mu \in \mathcal{M}_j^{-}$ such that $\phi = \phi_j^{\mu}$.
Using a Euclidean structure on $V$, we identify $\mu$ with a measure on $\widehat{\mathcal F}_{n-j}^{n-j+1}(V)$. We identify $F_0$ with the positively oriented unit vector in $F_1$ orthogonal to $F_0$.

Assume that $U$ lies inside a ball of radius $R$ around the origin, and choose a compactly supported even function $\zeta\in C^\infty(\R)$ such that $\zeta(t)=1$ for $|t|\leq 2R$. Define a functional $\psi$ on $\mathcal P(V)$ by setting for every $A\in\mathcal P(V)$
\begin{align*}
	\displaystyle&\psi(A)=\int_{F_1\in\Gr_{n-k+1}(V)}\int_{F_0\in S(F_1)} \int_{y\in F_1^\perp} \Bigg(\int_{t=0}^\infty\chi((F_0^++y+tF_0)\cap A)\zeta(t)dt \\
	&+\int_{t=-\infty}^0 \!\!\!\!\!\!\!\!\! \left(\chi((F_0^++y+tF_0)\cap A)-\chi((F_1+y)\cap A)\right)\zeta(t)dt\Bigg)\zeta(|y|)dyd\mu(F_0, F_1)
	\\&= \int_{F_1\in\Gr_{n-k+1}(V)}\int_{F_0\in S(F_1)}\int_{y\in F_1^\perp}\\&\qquad\qquad\qquad\qquad  \int_{t=-\infty}^\infty\!\!\!\!\!\!\!\!\!\chi((F_0^++y+tF_0)\cap A)\zeta(t)\zeta(|y|)dtdyd\mu(F_0, F_1),
\end{align*}
where the last equality holds since $\mu$ is odd in $F_0$ and thus $\chi((F_1+y)\cap A)$ integrates to zero. As $\mathrm{HGr}_{n-j+1}(V)$ is a homogeneous space under the affine special orthogonal group acting isotropically on $V$, it follows from \cite[Thm.~A.1]{Bernig2014} that $\psi\in\mathcal{V}_c^{\infty}(V)$. Using the natural parametrization of $\mathrm{HGr}_{n-j+1}(V)$ and letting $m$ be the measure on $\mathrm{HGr}_{n-k+1}(V)$ given by $\zeta(t)\zeta(|y|)dtdyd\mu(F_0, F_1)$, we may write
\[ \psi(A)= \int_{\mathrm{HGr}_{n-k+1}(V)}\chi(A\cap H)dm(H).\]
Observe that $\psi$ coincides with $\phi$ on convex bodies lying inside $U$. By \cite[Lemma 2.4.7]{Alesker2006}, it follows that $\phi$ and $\psi$ coincide on $U$.
\endproof

\subsection{Valuations on manifolds}
We turn now to Crofton formulas for valuations on a compact smooth manifold $M$, combining the results of the previous Section~\ref{sec:croftFormTransInvVal} with Theorem~\ref{mthm:nashThm}.
Recall that for $P \in \mathcal{P}(M)$ (resp. $\mathcal{P}(V)$), $\chi_P$ is the generalized valuation defined by $\psi \mapsto \psi(P)$, $\psi \in \mathcal{V}^\infty(M)$ (resp. $\mathcal{V}^\infty(V)$).

\begin{Proposition}\label{prop:CroftMnfldW1}
	Suppose that $M$ is a compact smooth manifold of dimension $m$, and let $e: M \hookrightarrow V$ be a perfectly non-parallel embedding into an $n$-dimensional space $V$. If $\phi \in \mathcal{W}_1^{\infty}(M)$, then there exist compactly supported measures $m_j$ on $\AGr_{n-j}(V)$ and $\mu_j$ on $\mathrm{HGr}_{n-j+1}(V)$ such that the identity
	\begin{align}\label{eq:CroftMnfldW1}
		\phi=\sum_{j=1}^{m} \int_{\AGr_{n-j}(V)}\!\!\!\!\!\!\!\! \chi_{E\cap M}\,dm_j(E) + \sum_{j=1}^{m-1} \int_{\mathrm{HGr}_{n-j+1}(V)}\!\!\!\!\!\!\!\!\chi_{H\cap M} d\mu_j(H)
	\end{align}
holds in two senses: in the Gelfand-Pettis (weak) sense, and in the integral-geometric sense of functionals on $\mathcal P(M)$.
\end{Proposition}
Let us first clarify the statement. For simplicity, consider the case
\begin{align}\label{eq:CroftMnfldW1Explain}
	\phi = \int_{\AGr_{n-j}(V)}\!\!\!\!\!\!\!\! \chi_{E\cap M}\,dm_j(E).
\end{align}
The Gelfand-Pettis sense of the integral is an equality in $\mathcal V^{-\infty}(M)$. Namely, for all $\psi\in\mathcal V^\infty(M)$, the $m_j$-a.e.\ defined function $E\mapsto \psi(E \cap M)$ is $m_j$-integrable over $\AGr_{n-j}(V)$, and we can write
\begin{align}\label{eq:crftFormMnfldGelPet}
	\langle \phi, \psi\rangle=\int_{\AGr_{n-j}(V)}\!\!\!\!\!\!\!\! \psi(E\cap M)\,dm_j(E),
\end{align}
for the Poincar\'e paring $\langle \cdot,\cdot\rangle$ in  $\mathcal{V}^{\infty}(M)$.

The integral-geometric sense of the equality means that for any $A\in\mathcal P(M)$, the $m_j$-a.e.\ defined function $E\mapsto \chi(E\cap A)$ is $m_j$-integrable over $\AGr_{n-j}(V)$, and
\begin{align}\label{eq:crftFormMnfldInsertA}
	\phi(A)=\int_{\AGr_{n-j}(V)}\!\!\!\!\!\!\!\! \chi(A \cap E)\,dm_j(E).
\end{align}

\begin{proof}[Proof of Proposition~\ref{prop:CroftMnfldW1}]
Let $\phi \in \mathcal{W}_1^{\infty}(M)$. By Theorem~\ref{mthm:nashThm}, there exists a smooth valuation $\Phi\in\Val^{\infty}(V)$ such that $e^*\Phi=\phi$. We will often identify $M$ with $e(M)$, and write $e^\ast \Phi = \Phi|_M$ accordingly. Decomposing $\Phi$ into its homogeneous components and further by parity, $\Phi=\sum_{j=0}^n \Phi_j^+ + \Phi_j^{-}$, we note that we can restrict ourselves to the case $\phi = \Phi_j^{\pm}|_M$, $0 \leq j \leq \dim M$, since $\Phi_j^{\pm}|_M = 0$ for $j > \dim M$. Moreover, $\Phi_0^- = 0$, $\Phi_{\dim M}^-|_{M} = 0$, and $\Phi_0^{+}|_M = 0$ as $\phi$ vanishes on points.
	
Consider first $\phi = \Phi_j^+|_M$, and choose a smooth Crofton measure $m_{j}$ for $\Phi_j^+$. Observe that by the transversality theorem (see, e.g., \cite{Guillemin1974}*{Ch.~2, §3}), $E\cap M \subset M$ is an embedded smooth submanifold for $m_j$-almost every $E \in \AGr_{n-j}(V)$.

Next, note that, since $m_j$ is a smooth Crofton measure for $\Phi_j^+ \in \mathcal{V}^\infty(V)$, also
\begin{align*}
	\Phi_j^+ = \int_{\AGr_{n-j}(V)} \chi_E dm_j(E)
\end{align*}
as elements of $\mathcal{V}^{-\infty}(V)$, in the Gelfand--Pettis sense. To see this, denote the integral on the right by $\widetilde \Phi_j^+\in \mathcal{V}^{-\infty}(V)$, which is well-defined as we now proceed to show. Recall that $V$ is an isotropic space under the action of $G=\SO(n)\ltimes V$. Fixing $E_0\in\AGr_{n-j}(V)$, it follows by \cite[Lemma A.2]{Bernig2014} that $g\mapsto \psi(gE_0)$ is locally integrable on $G$ for every $\psi \in \mathcal{V}_c^\infty(V)$. As the action of $G$ on $\AGr_{n-j}(V)$ is transitive and isotropic, it follows that $E \mapsto \psi(E)$ is locally integrable. Thus, $\widetilde\Phi_j^+ \in \mathcal{V}^{-\infty}(V)$ is well-defined. By \cite[Theorem.~A.1]{Bernig2014}, the assignment \[\mathcal P(V)\ni A\mapsto \int_{\AGr_{n-j}(V)}\chi(A\cap E)dm_j(E)\] is a well-defined smooth valuation.  Now by \cite[Lemma 2.4.7]{Alesker2006}, a smooth valuation is uniquely determined by its value on polytopes, and so this valuation coincides with $\Phi_j^+$. Moreover, by the definition of Alesker product and the second part of \cite[Thm.~A.1]{Bernig2014}, see also \cite[Thm.~1.1]{Fu2016},
\begin{align*}
	\langle \Phi_j^+, \psi \rangle = (\Phi_j^+ \cdot \psi)(V) = \int_{\AGr_{n-j}(V)} \psi(E) dm_j(E) = \langle \widetilde{\Phi}_j^+, \psi\rangle,
\end{align*}
for every $\psi \in \mathcal{V}^\infty_c(V)$, that is, $\Phi_j^+ = \widetilde{\Phi_j^+}$ as elements of $\mathcal{V}^{-\infty}(V)$.

We have seen in particular that for all $A\in\mathcal P(M)$, one has \[\phi(A)=\Phi_j^+(A)=\int_{\AGr_{n-j}(V)}\chi(A\cap E)dm_j(E).\]

For $\epsilon > 0$, let $Z_\epsilon\subset \AGr_{n-j}(V)$ be a precompact open subset of Euclidean measure $\sigma(Z_\epsilon)<\epsilon$, $Z_\epsilon'\subset Z_\epsilon$ an open subset with closure $\overline{Z_\epsilon'}\subset Z_\epsilon$ such that each $E\in \AGr_{n-j}(V)\setminus Z_\epsilon'$ intersects $M$ transversally. Such sets exist by the transversality theorem, and since $M$ is assumed to be compact. Let further $\eta_\epsilon(E)\in C^\infty(\AGr_{n-j}(V))$ satisfy $0\leq \eta_\epsilon\leq 1$, $\eta_\epsilon=1$ outside $Z_\epsilon$ and  $\eta_\epsilon=0$ inside $Z_\epsilon'$. Then the Gelfand--Pettis integral
\[\Phi_{j,\epsilon}^+:=\int_{\AGr_{n-j}(V)} \eta_\epsilon(E)\chi_Edm_j(E)\]
defines a valuation $\Phi_{j,\epsilon}^+ \in\mathcal V^{-\infty}(V)$, and $\Phi_{j,\epsilon}^+\to \Phi_j^+$ in $\mathcal V^{-\infty}(V)$ as $\epsilon\to 0$. Indeed, taking $\psi\in\mathcal V^\infty_c(V)$ it holds that $E\mapsto \psi(E)$ is a continuous, compactly supported function on $\AGr_{n-j}(V)$. The dominated convergence theorem then implies that $\langle \Phi_{j, \epsilon}^+, \psi\rangle \to \langle \Phi_{j}^+, \psi\rangle$. Since $\mathcal{V}^{-\infty}(V)$ is endowed with the weak-* topology, the assertion follows.

We may now use \cite[Claim~3.5.4 and Prop.~3.5.12]{Alesker2009} to write
\begin{align*}
	 \phi=e^*\Phi_j^+& = e^\ast \left(\lim_{\epsilon\to 0}\Phi_{j,\epsilon}^+\right)=\lim_{\epsilon\to 0} \left(\int_{\AGr_{n-j}(V)}\!\!\!\!\!\!\!\! \eta_\epsilon(E) \chi_{E\cap M}\, d m_j(E)\right).
\end{align*}

By the definition of the Gelfand--Pettis integral of a function with values in $\mathcal{V}^{-\infty}(M)$, the limit on the right-hand side satisfies
\begin{align*}
	\left\langle \lim_{\epsilon\to 0} \left(\int_{\AGr_{n-j}(V)}\!\!\!\!\!\!\!\!\!\!\!\!\!\!\!\! \eta_\epsilon(E) \chi_{E\cap M}\, d m_j(E)\right), \psi \right\rangle &= \lim_{\epsilon\to 0} \int_{\AGr_{n-j}(V)}\!\!\!\!\!\!\!\! \!\!\!\!\!\!\!\!\eta_\epsilon(E) \langle \chi_{E\cap M}, \psi \rangle\, d m_j(E)
\end{align*}
for every $\psi \in \mathcal{V}^\infty(M)$. 

Observe that $E\mapsto \psi(E\cap M)$ is a locally integrable function on $\AGr_{n-j}(V)$. To see this, first extend $p \mapsto \psi(\{p\})\in C^\infty(M)$ to a smooth function $f$ on $V$. Next, utilize the Euclidean structure on $V$ to define $\omega \in \Omega^{n-1}(\mathbb{P}_V)$, by $\omega_{(x,\xi)} = f(x) \vol_{S_xV}$, where $\vol_{S_x  V}$ is the normalized volume form on $S_xV=S^{n-1}$. Define the valuation $\widetilde\Psi \in \mathcal{V}^\infty(V)$ by
\begin{align*}
	\widetilde\Psi(A) = \int_{\nc(A)} \omega, \quad A \in \mathcal{P}(V).
\end{align*}
Its restriction $\widetilde \psi=\widetilde \Psi|_M$ satisfies $\widetilde\psi(\{p\}) = \psi(\{p\})$ for all $p \in M$, and so $\psi-\widetilde \psi\in\mathcal W_1^\infty(M)$. Now by Theorem~\ref{mthm:nashThm} we may choose $\Psi_1\in\Val^\infty(V)$ with $\Psi_1|_M=\psi - \widetilde{\psi} \in \mathcal{W}_1^\infty(M)$, so that $\psi(E\cap M)=\Psi_1(E\cap M) + \widetilde{\Psi}(E \cap M)$, whenever $E\in\AGr_{n-j}(V)$ intersects $M$ transversally. As before, since $V$ is an isotropic space under the action of $G=\SO(n)\ltimes V$, and fixing $E_0\in\AGr_{n-j}(V)$, it follows by \cite[Lemma A.2]{Bernig2014} that both $g\mapsto \Psi_1(gE_0\cap M)$ and $g\mapsto \widetilde \Psi(gE_0\cap M)$ are locally integrable on $G$. It then immediately follows that $\psi(E\cap M)$ is integrable on $\AGr_{n-j}(V)$.

By dominated convergence, using that $0 \leq \eta_\epsilon \leq 1$, we find that
\[ \lim_{\epsilon\to 0} \int_{\AGr_{n-j}(V)}\!\!\!\!\!\!\!\! \eta_\epsilon(E) \psi(E\cap M)d m_j(E) = \int_{\AGr_{n-j}(V)}\!\!\!\!\!\!\!\! \psi(E\cap M)d m_j(E),\]
and so
\begin{align*}
	\phi=\lim_{\epsilon\to 0} \int_{\AGr_{n-j}(V)}\!\!\!\!\!\!\!\! \eta_\epsilon(E) \chi_{E\cap M}\, d m_j(E) = \int_{\AGr_{n-j}(V)}\!\!\!\!\!\!\!\! \chi_{E\cap M}\,dm_j(E)
\end{align*}
in the sense of Gelfand--Pettis. This concludes the proof of the assertion in the even case, noting that we can replace $m_j$ by a compactly supported measure without changing the values for $A \in \mathcal{P}(M)$ as $M$ is compact.

\medskip

Next we consider the odd case. Let $\phi = \Phi_j^-|_M$, $1\leq j \leq \dim M-1$. Using Corollary \ref{cor:odd_crofton_compact} we choose a smooth, odd, compactly supported measure $\mu_j$ on $\mathrm{HGr}_{n-j+1}(V)$ such that 
\[ \Psi(A):=\int_{\mathrm{HGr}_{n-j+1}(V)}\chi(A\cap H)d\mu_j(H), \quad A\in\mathcal P(V)\]
is a smooth valuation on $V$, and $\Psi|_U=\Phi_j^-|_U$ for some neighborhood $U$ of $M$.
In particular, $\phi=\Psi|_M$, hence $\phi$ satisfies 
\[ \phi(A)=\int_{\mathrm{HGr}_{n-j+1}(V)}\chi(A\cap H)d\mu_j(H), \quad A\in\mathcal P(M).\]

By the second part of \cite[Thm.~A.1]{Bernig2014}, see also \cite[Thm.~1.1]{Fu2016},
\begin{align*}
	\langle \Psi, \psi\rangle = (\Psi \cdot \psi)(V) = \int_{\mathrm{HGr}_{n-j+1}(V)} \psi(H)d\mu_j(H),
\end{align*}
for all $\psi \in \mathcal{V}_c^\infty(V)$, and by the definition of Gelfand--Pettis integral,
\begin{align*}
	\Psi= \int_{\mathrm{HGr}_{n-j+1}(V)}\chi_Hd\mu_j(H).
\end{align*}

We proceed as in the even case to conclude that
\[ \phi= \int_{\mathrm{HGr}_{n-j+1}(V)}\chi_{H\cap M} d\mu_j(H)\]
in the sense of Gelfand--Pettis in $\mathcal V^{-\infty}(M)$, completing the proof.

\end{proof}

	\bigskip
	\noindent
	Theorem~\ref{mthm:croftonFormulaMf} will follow from Proposition~\ref{prop:CroftMnfldW1} combined with the following.
	\begin{Proposition}\label{prop:CroftMnfldW0}
		Suppose that $M$ is a compact smooth manifold and let $e: M \hookrightarrow V$ be an embedding into a linear space $V$ of odd dimension $n$. Let $\phi \in \mathcal{V}^{\infty}(M)$. Then there exists $R>0$, depending only on $e$, a compactly supported measure $\mu$ on $V \times [0,R]$ and a constant $C \in \R$ such that
		\begin{align*}
			\phi - C \chi - \int_{V \times [0,R]} \chi_{B_\rho(y) \cap M}\, d\mu(y,\rho) \in \mathcal{W}_1^\infty(M),
		\end{align*}
	where $B_d(x)$ is the closed Euclidean ball of radius $d>0$ centered at $x \in V$.
	\end{Proposition}
	\proof
	As $\mathcal{W}_1^\infty(M)$ consists of all valuations that vanish on points, we need to find $\mu$ and $C$, such that
	\begin{align}\label{eq:prfCroftFormW0DefPhiTilde}
		\widetilde{\phi} = C \chi + \int_{V \times [0,R]} \chi_{B_\rho(y) \cap M}\, d\mu(y,\rho)
	\end{align}
	defines a smooth valuation and $\phi(\{p\})=\widetilde \phi(\{p\})$ for all $p\in M$. To this end, define $f_\phi \in C^\infty(M)$ by $f_\phi(p) = \phi(\{p\})$. Identifying $e(M)$ with $M$ and choosing a Euclidean structure on $V$, there exists $r>0$ such that $M \subset B_r(0)$.
	
	Let now $f\in C^\infty_c(V)$ be an arbitrarily smooth function, such that $f|_M=f_\phi$ and $\Supp f \subset B_{2r}(0)$. We claim that there exists a smooth function $g\in C^\infty_c(V)$, $\Supp g \subset B_{2r}(0)$, such that
	\begin{align*}
		f(x) = \left(g \ast e^{-|\cdot|}\right)(x) = \int_{V} g(y) e^{-|x-y|}dy, \quad x \in V.
	\end{align*}
	Indeed, denoting by $\widehat{\bullet}$ the Fourier transform on $V$, we define $g$ by
	\begin{align*}
		\widehat{g} = \frac{\widehat{f}}{\widehat{e^{-|\cdot|}}} = c_n (1+|\cdot|^2)^{\frac{n+1}{2}} \widehat{f},
	\end{align*}
	where $c_n \in \R$ is some constant depending on $n$. Note that, since $\widehat{f}$ is a Schwartz function and $(1+|\cdot|^2)^{\frac{n+1}{2}}$ is smooth and all its derivatives have at most polynomial growth, $\widehat{g}$ is again a Schwartz function, that is $g$ is well-defined. Since $n+1$ is even, we can write
	\begin{align*}
		\widehat{g}(y) = c_n \sum_{i=0}^{\frac{n+1}{2}} \binom{\frac{n+1}{2}}{i} |y|^{2i} \widehat{f}(y) = c_n \sum_{i=0}^{\frac{n+1}{2}} (-1)^i\binom{\frac{n+1}{2}}{i} \widehat{\Delta^i f}(y),
	\end{align*}
	and, consequently,
	\begin{align*}
		g(x) = c_n \sum_{i=0}^{\frac{n+1}{2}} (-1)^i\binom{\frac{n+1}{2}}{i} \Delta^i f(x).
	\end{align*}
	We conclude that the support of $g$ is contained in the support of $f$, and thus in $B_{2r}(0)$. Moreover, using that $e^{-|\cdot|}$ defines a tempered distribution, we obtain
	\begin{align*}
		\langle f, h \rangle = \langle \widehat{f}, \widehat{h} \rangle = \langle \widehat{e^{-|\cdot|}}, \widehat{g} \, \widehat{h}\rangle = \langle \widehat{e^{-|\cdot|}}, \widehat{g \ast h}\rangle = \langle e^{-|\cdot|}, g \ast h\rangle = \langle g \ast e^{-|\cdot|}, h\rangle,
	\end{align*}
	for every Schwartz function $h$ on $V$, which implies $f = g \ast e^{-|\cdot|}$.
	
	Next, let $R=3r$, so that $\Supp g + B_r(0) \subset B_R(0)$ and note that $R$ only depends on the embedding of $M$ into $V$. We define $\widetilde{\Phi} \in \mathcal{V}^\infty(V)$ by
	\begin{align*}
		\widetilde{\Phi}(A) = \int_{0}^{R}\!\!\!\!\!e^{-\rho}d\rho\int_{V} \!\!\chi(A \cap B_\rho(y)) g(y) dy  + \left(e^{-R} \!\int_{V}\!\! g(y) dy \right) \chi(A), A \in \mathcal{P}(V).
	\end{align*}
	Note that, by \cite[Thm.~A.1]{Bernig2014} applied to the isotropic space $V$ with $G = \SO(n) \ltimes V$, the inner integral of the first term defines a smooth valuation on $V$ for every $\rho>0$. Since these valuations are continuous in $\rho$, and $g$ has compact support, the integral exists, and $\widetilde{\Phi} \in \mathcal{V}^\infty(V)$ is well-defined. As in the proof of Proposition \ref{prop:CroftMnfldW1}, using the second part of \cite[Thm.~A.1]{Bernig2014}, one deduces
	\begin{align*}
		\widetilde{\Phi} = \int_{0}^{R}e^{-\rho} d\rho\int_{V} \chi_{B_\rho(y)} g(y) dy  + \left(e^{-R} \int_{V} g(y) dy \right) \chi.
	\end{align*}
	Evaluating $\widetilde{\Phi}$ on $x \in M \subset B_r(0)$, we find
	\begin{align*}
		\widetilde{\Phi}(\{x\}) &= \int_{0}^{R}e^{-\rho} d\rho\int_{V} \mathbbm{1}_{B_\rho(y)}(x) g(y) dy  + e^{-R} \int_{V} g(y) dy  \\
		&= \int_{V} g(y) \left(\int_0^R \mathbbm{1}_{B_\rho(y)}(x) e^{-\rho} d\rho + e^{-R} \right) dy.
	\end{align*}
	For $\rho>R$, it holds for $y\in\Supp g$ that $|y|\leq 2r$ and so $\mathbbm 1_{B_\rho(y)}(x)=1$. Hence by the layer-cake formula, we have
	\begin{align*}
		\widetilde{\Phi}(\{x\}) = \int_{V} g(y) dy\int_0^\infty \mathbbm{1}_{B_\rho(y)}(x) e^{-\rho} d\rho  = \int_{V} g(y) e^{-|x-y|}dy = f(x).
	\end{align*}
	We are thus left to verify that $\widetilde{\phi} = \widetilde{\Phi}|_M$ is of the form \eqref{eq:prfCroftFormW0DefPhiTilde} with $\mu = (e^{-\rho} d\rho)\otimes(g(y)dy)$ and $C = e^{-R}\int_V g(y) dy$.
	
	Indeed, by Sard's theorem, almost every $e^{-\rho}$, with $\rho\in[0,R]$, is a regular value of $(e^{-|\cdot - y|})|_M$ and, by the transversality theorem, $M \cap B_\rho(y)$ is an embedded submanifold as superlevel set of $(e^{-|\cdot - y|})|_M$ for almost every $(y,\rho) \in V \times [0,R]$. Repeating the argument with $M$ replaced by $A \in \mathcal{P}(M)$, the same holds true for $A \cap B_\rho(y)$.	We can now use a similar argument as in the proof of Proposition~\ref{prop:CroftMnfldW1}, using again \cite[Claim~3.5.4 and Prop.~3.5.12]{Alesker2009}, to conclude the proof.
	\endproof

\begin{Remark}
	Fu \cite{Fu2016} introduced the notion of a \emph{principal kinematic valuation}, which is very similar to the integral representations from Theorem~\ref{mthm:croftonFormulaMf}, with the family of submanifolds being the orbit of a subset under the isotropic action of a group. The main difference with our approach is that we allow integration with families of subsets, which are only almost always submanifolds.
	
	Another closely related notion is the Radon transform with respect to the Euler characteristic, introduced by Alesker \cite{Alesker2009}. That setting is also too restrictive for our purposes, as it requires the submanifolds to come from a double fibration, which need not be the case (particularly for the case handled in Proposition~\ref{prop:CroftMnfldW0}).
\end{Remark}

As a direct consequence of Theorem~\ref{mthm:croftonFormulaMf}, the Alesker product of two smooth valuations on $M$, one given by a Crofton formula, takes the following simple form.
\begin{Corollary}
	Suppose that $M$ is a compact smooth manifold, embedded into $\R^N$ as in Theorem \ref{mthm:croftonFormulaMf}, and let $\phi, \psi \in \mathcal{V}^\infty(M)$, where $\phi$ is given by a Crofton formula as in eq. \eqref{eq:crofton}. Then 	for all $A \in \mathcal{P}(M)$,
	\begin{align*}
		(\phi \cdot \psi)(A) = \, &  C \psi(A) + \int_{\R^N \times [0,R]}\!\!\!\!\!\!\!\!\!\!\!\!\!\!\! \psi(B_\rho(y) \cap e(A))\, d\mu(y,\rho) \\
		&+ \sum_{j=1}^{n} \int_{\AGr_{N-j}(\R^N)}\!\!\!\!\!\!\!\!\!\!\!\!\!\!\!\!\!\!\!\!\!\! \psi(E\cap e(A))\,dm_j(E) + \sum_{j=1}^{n-1} \int_{\mathrm{HGr}_{N-j+1}(\R^N)}\!\!\!\!\!\!\!\!\!\!\!\!\!\!\!\!\!\!\!\!\!\!\!\!\!\!\!\psi(H\cap e(A)) d\mu_j(H).
	\end{align*}

\end{Corollary}
\begin{proof}
 	Let $\Phi\in\mathcal V^\infty(\R^n)$ be given by the same Crofton formula \eqref{eq:crofton} as $\phi$, and let $\Psi\in\mathcal V^\infty(\R^n)$ be an arbitrary valuation such that $\Psi|_M=\psi$ (which is easily shown to exist, but also follows from Theorem~\ref{mthm:nashThm} and Proposition~\ref{prop:CroftMnfldW0}). Then for $A\in\mathcal P(M)$ we have
	$(\phi\cdot\psi)(A)=(\Phi\cdot \Psi)(A)$, and \cite[Theorem.~A.1]{Bernig2014} completes the proof.
\end{proof}

\appendix

\section{A representation-theoretic proof of Theorem \ref{mthm:resProperty} for even valuations}\label{app:repProofResThm}
In this appendix, we present an alternative proof of Theorem \ref{mthm:resProperty} in the case of even valuations. Its interest stems from the very different toolset utilized in the proof, namely the Klain embedding and the representation theory of the general linear group. It would be interesting to extend this approach also to the odd case.

Before we give the proof, let us recall the composition series from \cite{Howe1999}. It is well-known that the isomorphism classes of irreducible representations of $\SO(n)$ can be parametrized by their highest weights $\lambda_1 \geq \dots \geq \lambda_{\lfloor n/2 \rfloor}$. Denoting by $\Lambda_0^+(k)$ the set of all highest weights with $\lambda_l = 0$ for $l \geq k+1$, and setting 
\begin{align*}
	R^+(1,s) = \begin{cases}
					\{ \mu \in \Lambda_0^+(k): \, \mu_i \in 2\N, \, 1+s \geq \mu_{s+1}\}, & \text{ if } s=0,1\\
					\{ \mu \in \Lambda_0^+(k): \, \mu_i \in 2\N, \, \mu_{s-1} \geq 1+s \geq \mu_{s+1}\}, & \text{ if } 2\leq s \leq k-1,
				\end{cases}
\end{align*}
a composition series of $\Gamma(\Gr_j(\R^n), \Dens(E))$ can then be written as
\begin{align}\label{eq:compSerKlainSec}
	\{0\} \subset R^+(1,1) \subset R^+(1,1) \oplus R^+(1,3) \subset \dots \subset \bigoplus_{s=0}^{\lfloor k/2 \rfloor} R^+(1,2s+1),
\end{align}
where $k = \min\{j, n-j\}$, see \cite[Thm.~3.4.4]{Howe1999}. In particular, the composition series is of length less or equal $2$, when $k \leq 3$. Moreover, by Alesker's irreducibility theorem~\cite[Thm.~1.3]{Alesker2001}, the space $\Val_j^{+,\infty}(\R^n)$ is an irreducible $\GL_n(\R)$-module with $\SO(n)$-types given by $R^+(1,1)$ (see also \cite{Alesker2011}), that is, it is the first module in above's series.

\medskip

\noindent
We have now all ingredients to prove
\begin{Proposition}\label{prop:resPropertyEvenRepTheoryPrf}
	Suppose that $1 \leq j \leq n-3$ and let $j+2 \leq r \leq n-1$.	Then for every $\psi \in U_j^{+,\infty}(r,\R^n)$ there exists a valuation $\psi \in \Val_j^{+,\infty}(\R^n)$ such that $\phi|_E = \psi_E$ for all $E \in \Gr_r(\R^n)$. If $j=1$, the conclusion holds for all $2 \leq r \leq n-1$.
\end{Proposition}
\proof
The assertion for $j=1$ follows directly, since in this case the Klain embedding is an isomorphism $\Val_1^{+,\infty}(\R^n)\cong \Gamma(\Gr_1(\R^n), \Dens(E))$. Hence, assume $2\leq j\leq r-2$ in the following.

Next, note that the restriction map $\res_r: \Val_j^{+,\infty}(\R^n)\to U_j^{+, \infty}(r,\R^{n})$ is injective by Klain's part of Theorem~\ref{thm:KlSchnInjective}, and by the Casselman--Wallach theorem~\cite{Casselman1989}, it has closed image. Using the Klain map, we can further map $U_j^{+, \infty}(r,\R^{n})$ injectively into $\Gamma^\infty(\Gr_j(\R^n), \Dens(E))$, where this map is well-defined by the compatibility of sections in $U_j^{+, \infty}(r,\R^{n})$. We therefore have the following chain of closed $\GL_n(\R)$-invariant subspaces
$$ \Val_j^{+,\infty}(\R^n)\subset U_j^{+,\infty}(r,\R^{n})\subset \Gamma^\infty(\Gr_j(\R^n), \Dens(E)),$$
where the last inclusion is strict. Indeed, fixing $\R^r\subset\R^n$, the Klain section of $\phi\in\Val_j^+(\R^r)$ has to be contained in the image of the cosine transform and, hence, cannot be arbitrary.

By \eqref{eq:compSerKlainSec}, the composition series has length $2$, when $j=2, 3, n-3$. We conclude that $\Val_j^+(\R^n)= U_j^{+,\infty}(r, \R^{n})$, which yields the claim for any $r$ in these cases.

Now, fix $4 \leq j \leq n-4$. As in the (general) proof in Section~\ref{sec:extFullGrass}, an inductive argument with basis $n=j+3$ and $r=j+2$ shows that we can restrict to $r=n-1$. Indeed, assume that the proposition is proved for $\dots, n-2, n-1$ for all $r$. By restricting the valuations, there is a natural map
\begin{align*}
	U_j^{+, \infty}(n-1,\R^{n})\to U_j^{+,\infty}(r,\R^{n}),
\end{align*}
which is injective by Theorem~\ref{thm:KlSchnInjective}. Let us show it is an isomorphism. To show surjectivity, let $\phi \in U_j^{+,\infty}(r,\R^n)$ be given and let $E \in \Gr_{n-1}(\R^n)$. Restricting, we find $\phi|_E \in U_j^{+, \infty}(r,E)$, and by the induction hypothesis there exists $\psi_E \in \Val_j^{+,\infty}(E)$ such that $\phi|_E = \res_r(\psi_E)$. Now $\psi:= (\psi_E)_{E \in \Gr_{n-1}(\R^n)}$ defines a smooth section over $\Gr_{n-1}(\R^n)$, as the corresponding Klain sections depend smoothly on $E$ (see Lemma~\ref{lem:charSmoothCompSecEven}). Thus $\psi\in  U_j^{+, \infty}(n-1,\R^{n})$ is the desired preimage of $\phi$.

It remains to show that the inclusion $\Val_j^{+,\infty}(\R^n)\subset U_j^{+, \infty}(n-1,\R^{n})$ is an equality. Assume in contradiction it is a proper inclusion, and fix a Euclidean structure on $\R^n$. By \eqref{eq:compSerKlainSec}, there exists $s \geq 1$ and an $\SO(n)$-submodule of $U_j^{+,\infty}(n-1,\R^{n})$ with $\SO(n)$-type given by $R^+(1,2s+1)$. In particular, the $\SO(n)$-module $V^n_{((2s+2)^k)}$ of highest weight $((2s+2)^k) = (2s+2,\dots, 2s+2, 0,\dots, 0)$, with $2s+2$ taken $k$ times and $0$ taken $\lfloor \frac n 2\rfloor -k$ times, for $k= \min\{j, n-j\} \geq 4$, appears in the restriction of $U_j^{+,\infty}(n-1,\R^{n})$ to $\SO(n)$.

Take now any hyperplane $H\in \Gr_{n-1}(\R^n)$ and consider the evaluation map
\begin{align*}
	\mathrm{ev}_H: U_j^{+,\infty}(n-1,\R^{n})\to \Val_j^{+,\infty}(H),
\end{align*}
where $\mathrm{ev}_H(\phi)=\phi(H)$. Clearly, $\mathrm{ev}_H$ is $\SO(H)\cong\SO(n-1)$-equivariant. 

If $n$ is even then the $\SO(n)$-module $V^n_{((2s+2)^k)}$, restricted to $\SO(n-1)$, is just the irreducible module $V^{n-1}_{((2s+2)^k)}$, which does not appear in $\Val_{j}^{+,\infty}(H)$ by Alesker's irreducibility theorem. Similarly if $n$ is odd then under the action of $\SO(n-1)$, 
$$V^n_{((2s+2)^k)}=\bigoplus_{i=-s-1}^{s+1}V^{n-1}_{((2s+2)^{k-1}, 2i)},$$
none of which can appear in $\Val_j^+(H) $ since $k\geq 4$.

Consequently, $\mathrm{ev}_H(V^n_{((2s+2)^k)})$ must vanish. Since this holds for all $H$, any $\phi\in V^n_{((2s+2)^k)}$ must be $0$, again by Theorem~\ref{thm:KlSchnInjective}. This is a contradiction and we conclude that $\Val_j^{+,\infty}(\R^n)$ coincides with $U_j^{+, \infty}(n-1,\R^{n})$, which finishes the proof.
\endproof

\section{The nuclear proof of the valuation Nash theorem}\label{app:NuclearNash}

In this appendix, we give an alternative proof of Theorem~\ref{mthm:nashThm}. As before, we will make use of a perfectly non-parallel embedding of the compact manifold $M$ into some $\R^n$. We then show, using a different method, that any section in $V_j^\infty(M)$, where $V_j^\infty(M)$ is the Fr\'echet space of smooth global sections of the Fr\'echet bundle over $M$ with values in $\Val_j^\infty(T_xM)$, is given by restricting a valuation in $\Val_j^\infty(\R^n)$. For $j=1$, this extension is rather straightforward and explicit. The key step is passing from $j=1$ to general $j$ using the Alesker product, the irreducibility theorem, and the nuclearity of the various spaces involved.

\medskip

Fix a smooth compact manifold $M$ and a perfectly non-parallel embedding $e: M \hookrightarrow \R^n$. In the following, we will write $\res_M:\Val_j^\infty(\R^n)\to V_j^\infty(M)$ for the map $\psi \mapsto [e^\ast \psi]_j$, where $j$ will be clear from context, and denote by $V_j^{\pm, \infty}(M)$ the Fr\'echet spaces of sections of even/odd valuations, defined similarly to $V_j^\infty(M)$. 

\begin{Lemma}\label{lem:nash1d}
	There is a continuous map $E:V^\infty_1(M)\to \Val_1^\infty(\R^n)$ such that $\mathrm{res}_M\circ E=\id$.
\end{Lemma}
\proof
The map is defined separately on even and odd valuations. Let us first consider $V^{+,\infty}_1(M)$. We write $V=\R^n$.
By assumption, the oriented projectivized tangent bundle $Z:=de(\mathbb P_+(TM))\subset \mathbb P_+(V)$ is an embedded submanifold.

Fix a Euclidean structure on $V$, which in turn induces a Riemannian metric on $\mathbb P_+(V) \cong S(V)$. 
Choose $\epsilon>0$ such that $Z_\epsilon:=\cup_{x\in Z} B_{\epsilon}(x)$, where $B_{\epsilon}(x)$ is the open ball of radius $\epsilon$ centered at $x$, is a tubular neighborhood of $Z$ on which the least distance projection $\pi_Z:Z_\epsilon\to Z$ is well-defined and smooth.

Fix a function $\alpha\in C^\infty (\R)$ which is identically $1$ near $0$, and identically $0$ outside $[-\frac 12,\frac 12]$. Using the Klain embedding, we may identify $\phi\in V_1^{+, \infty}(M)$ with its Klain section, considered as an even function $f\in C^\infty(Z)$. For $\theta\in\mathbb P_+(V)$, set $$\widetilde f(\theta)=\alpha\left(\frac{d(\theta, Z)}{\epsilon}\right) f(\pi_Z \theta).$$
The first factor is indeed a smooth function on $Z_\epsilon$, and can be smoothly extended by $0$ outside $Z_\epsilon$. Clearly, $\widetilde \psi$ is even. We then define $E\phi$ to be the unique valuation in $\Val^{+,\infty}_1(V)$ with Klain section $\widetilde f$.

Now consider $\phi\in V_1^{-,\infty}(M)$, which we identify with a smooth choice of valuations $\phi_E\in\Val_1^{-,\infty}(E)$ for $E\in \mathrm{Image}(de: M\to \Gr_m(V))$. Applying the Alesker--Fourier transform, we obtain a family of (twisted) valuations $\psi_E=\mathbb F \phi_E\in \Val^{-,\infty}_{m-1}(E^*)\otimes \Dens(E)$.

The Euclidean structure on $V$ trivializes $\Dens(E)\cong\R$ for all $E$, and identifies $E^*\cong  E$, $V^*\cong V$.
By Schneider's theorem we can find a smooth family of odd functions $f_E\in C^\infty(S(E))$ such that 
$$\psi_E(K)=\int_{S(E)} f_E(\theta)dS_{m-1}(K;\theta).$$
Furthermore, $f_E$ can be chosen to be orthogonal in $L^2(S(E))$ to all restrictions of linear functions to $S(E)$. Such a choice is unique, and $E\mapsto f_E$ is then smooth.

As before, since $e$ is perfectly non-parallel, $S(E)\cap S(E')=\emptyset$ whenever $E\neq E'$ is a tangent plane of $e(M)$, and we may define $f:Z\to \R$ by setting $f|_{S(E)}=f_E$. The map $f$ is smooth by construction. 

Now extend $f$ to an odd function $\widetilde f$ on $S(V)$, similarly to the even case. 
Define the valuation $\psi\in\Val_{n-1}^{-,\infty}(V)$ by 
$$\psi(K)=\int_{S(V)} \widetilde f(\theta)dS_{n-1}(K;\theta).$$
Denoting by $\pi_{E}:V\to E$ the orthogonal projection, it now follows that $\psi_E=(\pi_{E})_*\psi$ for all $E$. It remains to set 
$E\phi:=\mathbb F^{-1}\psi$, since $\mathbb F^{-1}(\pi_{E})_*=i_E^*\mathbb F^{-1}$, where $i_E:E\hookrightarrow V$ is the inclusion map.
\endproof

To extend the result from Lemma~\ref{lem:nash1d} to other degrees of homogeneity, we will makes use of the nuclearity of the corresponding spaces, established below.

\begin{Lemma}\label{lem:nuclear}
	The Fr\'echet spaces $\Val_j^\infty(\R^n)$ and $V_j^\infty(M)$ are nuclear.
\end{Lemma}

\proof 
We embed $M$ into some $\R^n$ and endow $M$ with the Riemannian structure inherited from $\R^n$. If $j=0$ or $j=\dim M$, then $V_j^\infty(M) \cong C^\infty(M)$, which is nuclear. We may therefore assume $1 \leq j \leq \dim M - 1$ in the following.

By Section~\ref{sec:bgValDiffForm}, it holds for $x \in M$ that $\Val_j^\infty(T_xM)$ is naturally a subspace of
\begin{align*}
	\Gamma(S(T_x^\ast M), \wedge^{m-j}\xi^\perp \otimes \wedge^{m-j} T_x M),
\end{align*}
which, by composing with the embeddings $\xi^\perp \subset T_x M$ and $T_x M \subset \R^n$, is a subspace of $C^\infty(S(T_x^\ast M), \wedge^{m-j} \R^n \otimes \wedge^{m-j}\R^n)$. We conclude that $V_j^\infty(M) = \Gamma(M, \Val_j^\infty(T_x M))$ is a linear subspace of $C^\infty(S^\ast M, \wedge^{m-j} \R^n \otimes \wedge^{m-j} \R^n)$. As the latter space is isomorphic to $C^\infty(S^\ast M)\otimes (\wedge^{m-j} \R^n \otimes \wedge^{m-j} \R^n)$, which is nuclear by \cite[Prop.~50.1, (50.9)]{Treves1967}, we see that $V_j^\infty(M)$ is nuclear by \cite[Prop.~50.1]{Treves1967}.

As the space of translation-invariant forms is a subspace of the space of forms, the same reasoning applies for $\Val_j^\infty(\R^n)$, completing the proof.
\endproof 

By nuclearity, the injective and projective topologies on the tensor product coincide and we denote the (unique) completed tensor product by $\widehat{\otimes}$.
\begin{Corollary}
	The Alesker product $$m:\Val_i^\infty(\R^n)\times \Val_j^\infty(\R^n)\to \Val_{i+j}^\infty(\R^n)$$ admits a continuous extension to the completed tensor product, denoted $$\widehat m:\Val_i^\infty(\R^n)\widehat \otimes \Val_j^\infty(\R^n)\to \Val_{i+j}^\infty(\R^n).$$ Similarly, the Alesker product $m:V_i^\infty(M)\times V_j^\infty(M)\to V_{i+j}^\infty(M)$ extends to
	$$\widehat m:V_i^\infty(M)\widehat \otimes V_j^\infty(M)\to V_{i+j}^\infty(M).$$ 
\end{Corollary}
\proof
This follows from  the universal property of the projective topology.
\endproof

Denoting the respective extensions of the $k$-fold product maps to the completed tensor $k$-powers by $$\widehat m_k:(\Val_1^\infty(\R^n))^{\widehat\otimes k}\to \Val_k^\infty(\R^n),\qquad \widehat m_k^M:(V_1^\infty(M))^{\widehat\otimes k}\to V_k^\infty(M),$$ 
and similarly extending the map $E$ of Lemma~\ref{lem:nash1d} to a continuous operator
$$E^k:=E^{\widehat\otimes k}: V_1^\infty(M)^{\widehat\otimes k}\to \Val_1^\infty(\R^n)^{\widehat\otimes k},$$
the product maps on $M$ and on $\R^n$ can be related.

\begin{Lemma}\label{lem:extend_multiply_restrict}
	$\widehat m_k^M:(V_1^\infty(M))^{\widehat\otimes k}\to V_k^\infty(M)$ coincides with $ \mathrm{res}_M\circ \widehat m_k \circ E^k$.
\end{Lemma}
\proof
It suffices by continuity and linearity to check the identity on elements of the form $\phi_1\otimes\dots\otimes\phi_k$. This amounts to the verification
$$  \mathrm{res}_M( E(\phi_1)\cdots E(\phi_k))= \phi_1\cdots\phi_k,$$
which follows from Lemma \ref{lem:nash1d} since the Alesker product commutes with restriction.
\endproof

It remains to see that all $\phi \in V_k^\infty(M)$ lie in the image of $\widehat{m}_k^M$.

\begin{Lemma}\label{lem:product_onto}
	The map $\widehat m_k^M: (V_1^\infty(M))^{\widehat\otimes k}\to V_k^\infty(M)$ is onto.
\end{Lemma}
\proof
First note that both $(\Val_1^\infty(\R^m))^{\widehat\otimes k}$ and $\Val_k^\infty(\R^m)$ are Fr\'echet spaces that are nuclear by Lemma~\ref{lem:nuclear} and \cite[Prop.~50.1, (50.9)]{Treves1967}. By Alesker's irreducibility theorem and the Casselman-Wallach theorem~\cite{Casselman1989}, the map $\widehat m_k:(\Val_1^\infty(\R^m))^{\widehat\otimes k}\to \Val_k^\infty(\R^m)=\Val_k^{+,\infty}(\R^m)\oplus \Val_k^{-,\infty}(\R^m)$ is surjective.

Next, consider the case $M = \R^m$. Then $V_k^\infty(\R^m) = C^\infty(\R^m, \Val_k^\infty(\R^m))$, and, by \cite[Rem.~1.1.9]{Alesker2006b}, $\widehat{m}_k$ extends to a continuous and surjective map
\begin{align*}
	C^\infty(\R^m, (\Val_1^\infty(\R^m))^{\widehat\otimes k}) \to V_k^\infty(\R^m).
\end{align*}
Observe that the natural map $C^\infty(\R^m)^{\widehat{\otimes}k} \to C^\infty(\R^m)$ is surjective and, as all spaces below are Fr\'echet, the map
\begin{align*}
	C^\infty(\R^m)^{\widehat{\otimes} k} \widehat\otimes \Val_1^\infty(\R^m)^{\widehat\otimes k} \to C^\infty(\R^m) \widehat\otimes \Val_1^\infty(\R^m)^{\widehat\otimes k} 
\end{align*}
is also surjective. The claim follows in this case as, by nuclearity, \begin{align*}V_1^\infty(\R^m)^{\widehat{\otimes}k} &\cong C^\infty(\R^m)^{\widehat{\otimes} k} \widehat\otimes \Val_1^\infty(\R^m)^{\widehat\otimes k} \\
	C^\infty(\R^m, (\Val_1^\infty(\R^m))^{\widehat\otimes k}) &\cong C^\infty(\R^m) \widehat\otimes \Val_1^\infty(\R^m)^{\widehat\otimes k}.
\end{align*}

Next, consider a general compact manifold $M$, and let $U\subset M$ be a neighborhood that is diffeomorphic to $\R^m$. Observe that a $k$-tuple of compactly supported functions $f_1,\dots,f_k\in C^\infty_c(U)$ defines a continuous multilinear map
\begin{align*}
	(V_1^\infty(U))^k\to (V_1^\infty(M))^{\widehat\otimes k},\qquad (\psi_1,\dots,\psi_k)\mapsto f_1\psi_1\otimes\dots\otimes f_k\psi_k,
\end{align*}
and so extends to a continuous linear map $S_{f_1,\dots,f_k}: (V_1^\infty(U))^{\widehat\otimes k}\to  (V_1^\infty(M))^{\widehat\otimes k}$. It holds that
$$\widehat m^M_k(S_{f_1,\dots, f_k}\psi)=f_1\dots f_k\widehat m^M_k(\psi),$$
as this is clearly true for decomposable $\psi$, which is sufficient.

Fix a finite open cover of $M$ by precompact, contractible sets $U_i$. For $\phi\in V^\infty_k(M)$, we can use a partition of unity to write $\phi=\sum\phi_i$ with $\Supp(\phi_i)\subset U_i$. As every $U_i$ is diffeomorphic to $\R^m$, by the first step, we may choose $\psi_i\in(V^\infty_{1}(U_i))^{\widehat \otimes k}$ with $\widehat m_k^M(\psi_i)=\phi_i$. Fixing smooth functions $w_i\in C^\infty(M)$ with $\Supp(w_i)\subset U_i$ such that $w_i=1$ on $\Supp(\phi_i)$, it follows that $$\widehat m^M_k(S_{w_i,\dots, w_i}(\psi_i))=w_i^k\widehat m^M_k(\psi_i)=w_i^k\phi_i=\phi_i.$$ 
Setting $\psi_i':=S_{w_i,\dots, w_i}(\psi_i)\in V_1^\infty(M)^{\widehat\otimes k}$, we have $\Supp(\psi_i')\subset U_i$ and so may write $\phi=\widehat m_k^M(\sum\psi_i)$ since the cover is locally finite, concluding the proof.

\endproof 

\begin{Corollary}\label{cor:product_onto}
	$\mathrm{res}_M:\Val^\infty_k(\R^n) \to V_k^\infty(M)$ is onto.
\end{Corollary}
\proof
For $\phi\in V_k^\infty(M)$, we may by Lemma \ref{lem:product_onto} choose $\psi \in (V_1^\infty(M))^{\widehat\otimes k}$ such that $\widehat m_k^M(\psi)=\phi$. By Lemma \ref{lem:extend_multiply_restrict}, $\mathrm{res}_M(\widehat m_k( E_1^k(\psi)))=\phi$.
\endproof

Corollary \ref{cor:product_onto} can now be utilized in the proof of Theorem \ref{thm:nash_firstproof} instead of Theorem~\ref{mthm:nashThm2}, completing the alternative proof of Theorem~\ref{mthm:nashThm}.

\bigskip

\subsection*{Acknowledgements}
The authors would like to thank Semyon Alesker, Andreas Bernig, Joe Fu, Gil Solanes and Thomas Wannerer for valuable discussions and suggestions.

\bibliographystyle{abbrv}
\bibliography{./books}
\end{document}